\documentclass[11pt]{article}
\usepackage[T1]{fontenc}
\usepackage[utf8]{inputenc}
\usepackage{lmodern}
\usepackage[a4paper,margin=27mm,headheight=14pt]{geometry}
\usepackage{microtype}
\usepackage{amsmath,amssymb,amsthm,mathtools,mathrsfs}
\usepackage{graphicx,booktabs,longtable,array,multirow}
\usepackage{tikz}
\usetikzlibrary{arrows.meta,calc,positioning,decorations.pathreplacing,angles,quotes}
\usepackage{enumitem}
\usepackage{caption,subcaption}
\usepackage{placeins}
\usepackage{fancyhdr}
\usepackage{aliascnt}
\usepackage[hidelinks,pdfencoding=auto]{hyperref}
\usepackage[nameinlink,noabbrev]{cleveref}
\usepackage{listings}
\numberwithin{equation}{section}
\newtheorem{theorem}{Theorem}[section]
\newaliascnt{proposition}{theorem}
\newtheorem{proposition}[proposition]{Proposition}
\aliascntresetthe{proposition}
\crefname{proposition}{proposition}{propositions}
\Crefname{proposition}{Proposition}{Propositions}
\newaliascnt{lemma}{theorem}
\newtheorem{lemma}[lemma]{Lemma}
\aliascntresetthe{lemma}
\crefname{lemma}{lemma}{lemmas}
\Crefname{lemma}{Lemma}{Lemmas}
\newaliascnt{corollary}{theorem}
\newtheorem{corollary}[corollary]{Corollary}
\aliascntresetthe{corollary}
\crefname{corollary}{corollary}{corollaries}
\Crefname{corollary}{Corollary}{Corollaries}
\newtheorem{classicalthm}{Theorem}

\crefname{classicalthm}{Theorem}{Theorems}
\Crefname{classicalthm}{Theorem}{Theorems}
\theoremstyle{definition}
\newaliascnt{definition}{theorem}
\newtheorem{definition}[definition]{Definition}
\aliascntresetthe{definition}
\crefname{definition}{definition}{definitions}
\Crefname{definition}{Definition}{Definitions}
\newaliascnt{example}{theorem}
\newtheorem{example}[example]{Example}
\aliascntresetthe{example}
\crefname{example}{example}{examples}
\Crefname{example}{Example}{Examples}
\theoremstyle{remark}
\newaliascnt{remark}{theorem}
\newtheorem{remark}[remark]{Remark}
\aliascntresetthe{remark}
\crefname{remark}{remark}{remarks}
\Crefname{remark}{Remark}{Remarks}
\newcommand{\Dg}{\widetilde D}

\newcommand{\ZZ}{\mathbb Z}
\newcommand{\RR}{\mathbb R}
\newcommand{\CC}{\mathbb C}
\newcommand{\Def}{\mathrel{\vcentcolon=}}
\newcommand{\ind}{\mathbf 1}
\newcommand{\ee}{\mathrm e}
\newcommand{\ii}{\mathrm i}
\newcommand{\dif}{\,\mathrm d}

\newcommand{\Hcal}{\mathcal H}
\newcommand{\Ncal}{\mathcal N}
\newcommand{\Gcal}{\mathscr G}
\newcommand{\Lphase}{L_{\mathrm{phase}}}
\newcommand{\Lfar}{L_{\mathrm{far}}}
\newcommand{\Lhigh}{L_{\mathrm{high}}}
\DeclareMathOperator{\des}{des}

\DeclareMathOperator{\sgn}{sgn}
\DeclareMathOperator{\Res}{Res}
\DeclareMathOperator{\diag}{diag}
\DeclareMathOperator{\sech}{sech}
\DeclareMathOperator{\csch}{csch}

\graphicspath{{figures/}}
\setlist{topsep=4pt,itemsep=2pt,parsep=0pt}
\hypersetup{pdftitle={Real-Rootedness and Log-Concavity of Zigzag Eulerian Polynomials},pdfauthor={Ethan Torres}}
\title{Real-Rootedness and Log-Concavity of\protect\\ Zigzag Eulerian Polynomials}
\author{Ethan Torres\\\small\texttt{ethanjt2@illinois.edu}}
\date{November 2026}
\begin{document}
\maketitle
\begin{abstract}
We prove that the zigzag Eulerian polynomials have only real, negative,
simple zeros, and hence strictly log-concave coefficients, settling the
first two parts of a conjecture of Petersen and Zhuang. Real-rootedness in
combinatorics is usually proved inside a single algebraic framework, such as
interlacing or stability. Our proof rests on a different premise: different
parts of a zero set become visible in different representations of the same
polynomial, and classical tools suffice once each part is viewed in the
right one. The bridge between representations is a local rule. Forgetting
the order in which a bounded-jump path visits positions leaves only a
nearest-neighbor condition on sorted gaps, which makes an encoding via a transfer
matrix much simpler to see. From there the argument passes to complex analysis. The sizes of the
matrix's eigenmodes locate the zeros at the two ends of the negative axis,
through Rouch\'e's theorem. In a logarithmic coordinate the generating
function becomes a combination of derivatives of the hyperbolic secant,
akin to the classical derivative polynomials, and the arguments of its
nearest poles produce a sine that fixes the signs between consecutive
central zeros. The zeros are thus read off from the poles of an associated
kernel. More broadly, the proof suggests that combinatorial polynomials
defined by local constraints carry analytic structure, in their spectra and
in the poles of their kernels, that is complementary to the usual algebraic
methods and that classical complex analysis can exploit directly.
\end{abstract}
\tableofcontents

\section*{Statement on the use of AI}

I want to be fully transparent about how AI was used in this work, and
equally clear about what it was not used for: this is not an AI-generated
proof. The core of the argument is my own. I constructed the recombining tree and
its position-visit quotient, proved the locality theorem, and developed the
Rouch\'e argument for the edge zeros through my own thinking and intuition.
While that argument was my only line of attack, I used AI to probe it for
mathematical weaknesses, which helped me see the central sine problem: the
gap that the Rouch\'e argument alone leaves in the central region. I then
built the central sine lemma with substantial help from AI. The geometric
concepts behind it, and their intuition and description throughout the
paper, are mine (and were subsequently verified and made more clear to a reader by AI).

AI assisted in three further ways. Once the main ideas were in place and I
had restated the theorems more concisely, AI helped check and verify the
parts of the proofs I was not working out myself, and I reviewed those
checks. All code used to verify the results, including the coefficient
checkers and validators, was designed by AI. AI also helped with LaTeX
formatting and with keeping the notation consistent with conventions in
both complex analysis and algebraic combinatorics, so that the paper can
reach the broadest possible audience.

AI played its largest role in two places: producing this particular draft,
and an extensive literature review. The paper itself is a compression of roughly 150
pages of mathematics that I wrote out myself. AI carried out the
compression, selecting the core elements and pruning long computations the
argument does not need, and I reviewed the result. I chose this route so
that readers can use these techniques now, rather than work through 150
pages of notationally inconsistent mathematics in which the core ideas, and
the steps leading to them, were hard to follow. I will keep condensing the
paper in later arXiv revisions toward a journal-ready version, and the
longer drafts will remain available to anyone who wants more complete
proofs.

In its current form, I am comfortable releasing this paper and claiming it
as my own work. It has been through extensive revision but it is not yet in its
final form, and I will continue to revise it for journal submission. Please, notify the author of any errata.

\section{Results and proof architecture}\label{sec:introduction}
 
Long before anyone asked about the zeros of the zigzag Eulerian
polynomials, these polynomials had already shown up in a number of
constructions, and it is worth seeing how the question of their zeros
arose. The common thread is the order polynomial of the fence, the poset
whose elements alternate up and down.
The order-polynomial and descent-enumerator connection comes from
Stanley's theory of $P$-partitions
\cite[Propositions~13.2--13.3]{StanleyPartitions}, in the modern formulation
of \cite[Section~3.15]{StanleyEC1}. Its reciprocity is the order-polynomial
reciprocity theorem \cite[Proposition~2.1]{StanleyReciprocity}.
Stanley's order and chain polytopes provide the general lattice-point
framework \cite[Definition~2.1, Theorem~4.1]{StanleyPolytopes}.
Chen and Zhang proved unimodality of the fence chain-polytope numerator
\cite[Section~3]{ChenZhang}. Coons and Sullivant identified the same
numerator through the swap statistic on alternating permutations and proved
symmetry and unimodality geometrically
\cite[Theorems~1.9 and~4.1]{CoonsSullivant}.
 
The same count appears in several other guises: as a magic labeling of
a path, as a symmetric tridiagonal matrix with constant row and column
sums, and through the alternating heights of Kekul\'e structures. These
graph-labeling, tridiagonal-matrix, and Kekul\'e
interpretations belong to established enumerative frameworks
\cite{StanleyMagic,DiaconisWood,CyvinGutman}; they help identify the object,
but are not assertions of new equinumeracy. We collect these realizations in Remark~\ref{rem:classical-models}.
 
Petersen and Zhuang give the systematic zigzag Eulerian treatment, including
the big-return interpretation, refined recurrences, and
$\gamma$-nonnegativity \cite[Section~3.3, Theorems~1.1 and~2.5]{PZ}.
Their last conclusion is an application of Br\"and\'en's theorem on
sign-graded posets \cite[Theorem~4.2]{BrandenSign}; we state the specialization
in Section~\ref{sec:dictionary}. Thus symmetry and unimodality precede the
real-rootedness question.
 
To say precisely what we prove, we fix the objects. The fence on $\Dg$
elements is the poset $x_1<x_2>x_3<x_4>\cdots$, and an order-preserving map
from it to the chain $1<2<\cdots<M$ is the same thing as a word whose
letters alternately weakly rise and weakly fall.
For integers $\Dg\ge1$ and $M\ge1$, let $\Omega_{\Dg}(M)$ count these weakly
alternating words
\[
(y_1,\ldots,y_{\Dg})\in\{1,\ldots,M\}^{\Dg},
\qquad y_1\le y_2\ge y_3\le y_4\ge\cdots.
\]
These are the order polynomials of fence, or zigzag, posets. The generating functions of these order polynomials
define the normalized numerators $h_{\Dg}$ by
\begin{equation}\label{eq:intro-numerator}
 \sum_{M\ge1}\Omega_{\Dg}(M)t^{M-1}
 =\frac{h_{\Dg}(t)}{(1-t)^{\Dg+1}},
 \qquad Z_{\Dg}(t)=t h_{\Dg}(t).
\end{equation}
The convention $Z_{\Dg}=t h_{\Dg}$ is the one used to compare with the zigzag
Eulerian polynomials of Petersen and Zhuang \cite[Eq.~(1), Lemma~2.2]{PZ}.
Section~\ref{sec:dictionary} proves the polynomiality and records the precise
normalization. In particular, $h_1=h_2=1$, and $h_{\Dg}$ is a monic reciprocal
polynomial of degree $\Dg-2$ for $\Dg\ge2$.
 
Our goal in this paper is to prove the following two statements about
these numerators.
 
\begin{theorem}[Real-rootedness and simplicity]\label{thm:main}
For every integer $\Dg\ge3$, the polynomial $h_{\Dg}$ has exactly $\Dg-2$
distinct strictly negative zeros. Consequently, $Z_{\Dg}=t h_{\Dg}$ has only
simple real zeros: one at the origin and all others strictly negative.
For $\Dg=1,2$, one has $Z_{\Dg}(t)=t$.
\end{theorem}
 
\begin{corollary}[Coefficient log-concavity]\label{cor:main-lc}
Write $h_{\Dg}(t)=\sum_{j=0}^{m}\mathsf h_jt^j$, where $m=\Dg-2$ and $\Dg\ge2$.
All $\mathsf h_j$ are positive, and
\begin{equation}\label{eq:intro-newton}
 \mathsf h_j^2\ge
 \frac{(j+1)(m-j+1)}{j(m-j)}\mathsf h_{j-1}\mathsf h_{j+1}
 >\mathsf h_{j-1}\mathsf h_{j+1},\qquad 1\le j<m.
\end{equation}
Thus the coefficient sequence is strictly log-concave at every interior
index and is unimodal. The coefficient sequence of $Z_{\Dg}$ is log-concave,
with its initial zero understood in the usual way.
\end{corollary}
 
Together, these two statements are exactly parts~(a) and~(b) of the
conjecture that Petersen and Zhuang posed for the zigzag Eulerian
polynomials \cite[Conjecture~6.1]{PZ}, which they had checked by computer
for every size up to~$35$. Their conjecture has a third part. Part~(c)
asks that the $\gamma$-vector of $h_{\Dg}=Z_{\Dg}/t$, that is, the vector of
coefficients $\gamma_{\Dg,j}$ in the expansion \eqref{eq:known-gamma-expansion}
recorded in Section~\ref{sec:dictionary}, satisfy the
Frankl--F\"uredi--Kalai inequalities; equivalently, that it be the
$f$-vector of a balanced simplicial complex. Petersen and Zhuang place this
part in the line of Nevo and Petersen, who conjectured that the
$\gamma$-vectors of flag homology spheres satisfy these inequalities
\cite{NevoPetersen}. Part~(c) looks as though it should follow from the
real-rootedness proved here, but the author is not yet experienced enough
with these inequalities to claim it, and we leave it open.
Nor is this the only open direction around these polynomials: Petersen
and Zhuang also ask whether real-rootedness persists for the relaxed
generalizations they introduce \cite[Question~6.8]{PZ}, a question we do
not take up here.
 
The order in which the earlier results arrived, with symmetry and
unimodality established well before any statement about zeros, is the
natural one. For a polynomial with positive
coefficients, the familiar properties form a ladder \cite{Branden}. Real
zeros give log-concave coefficients, by Newton's inequalities, and positive
log-concave coefficients are unimodal. If the polynomial is also symmetric,
real zeros give $\gamma$-nonnegativity as well. Each earlier result climbs
part of this ladder from below, while Conjecture~6.1 of \cite{PZ} asks for
the top rung directly. Theorem~\ref{thm:main} supplies it, and
Corollary~\ref{cor:main-lc} is the first step back down. The methods that
give symmetry, unimodality, and $\gamma$-nonnegativity do not by themselves
reach the zeros, so the next question is how one gets there.
 
\subsection{The proof in three parts}\label{sec:three-parts}
 
The real zeros of an enumerative polynomial can be difficult to see in its
original combinatorial definition. The difficulty is not necessarily a lack
of structure. Rather, different parts of the zero set may become visible in
different representations. In this paper, a local constraint on the geometry
of path visits first gives a finite matrix, and then a meromorphic kernel.
The matrix controls the extreme zeros through the sizes of its eigenmodes;
the kernel controls the central zeros through the arguments of its poles.
The proof is organized around this change of viewpoint.
 
The proof gives more than a root count.
For every size in its uniform range, it constructs explicit negative
intervals containing the zeros, one per interval. Their endpoints come from
two different geometric mechanisms, and the argument verifies that these
intervals fit together without an overlap or a missing root.
 
Those two mechanisms divide the zeros of $h_{\Dg}$ into two groups, and it
helps to name them now. The \emph{edge zeros} are the zeros nearest the
origin, together with their reciprocal partners at the far negative end;
the matrix sees them. The \emph{central zeros} are the rest. They sit
around $t=-1$, the fixed point of the reciprocity $t\mapsto1/t$, and the
kernel sees them.
The argument has three logically separate parts. Each ends with a statement
that is used, without further reference to the method, by the next.
\begin{enumerate}[label=(\roman*)]
\item \emph{Near-zero and far roots} (Sections~\ref{sec:locality}--\ref{sec:edges}).
Section~\ref{sec:locality} shows that a path's visit profile is
governed by one local rule: two consecutive sorted gaps never add up to
more than the jump radius. Read as a compatibility condition between
neighboring gaps, this rule gives the threshold matrix; its spectrum gives nested
balanced circles on which one term dominates; Rouch\'e's theorem then
counts the zeros nearest the origin one at a time, and reciprocity places
their partners at the far negative end.
\item \emph{Explicit endpoint verification}
(Sections~\ref{sec:finite} and~\ref{sec:joining}). A short exact
computation handles $3\le\Dg\le10$. For larger sizes the outer ends of the
central grid are compared with the edge intervals in the logarithmic
coordinate, in Sections~\ref{sec:boundary-separation}
and~\ref{sec:boundary-intervals}, so that the
two families of intervals abut with no overlap and no gap.
\item \emph{The central sine lemma} (Sections~\ref{sec:geometry}
and~\ref{sec:central}). This part handles the central zeros. A conjugate pole pair of the hyperbolic secant
produces an oscillation whose extrema are the sign grid; the three-loss
decomposition shows that the discrepancy from the full polynomial cannot
reverse any sign.
\end{enumerate}
Why the zeros divide in this way, and why the edge mechanism cannot
simply be pushed into the middle, is explained next.
 
\subsection{Why the proof has an edge and a central part}
 
Start with the edges, where the matrix does the work. The threshold matrix is symmetric and has an explicit sine spectrum. At a
fixed alphabet size, one eigenmode is strictly larger in modulus than all the
others. A uniform remainder estimate makes that observation useful while
the alphabet size varies. The logarithm of the leading eigenvalue is
strictly concave. Consequently, one can balance two neighboring growth
rates and choose a circle on which a single term of
\eqref{eq:intro-numerator} dominates the sum of the rest. Rouch\'e's theorem
then counts zeros inside that circle. At its intersection with the negative
axis the same domination gives a sign. Nested circles and alternating signs
therefore locate the zeros closest to zero; reciprocity locates their far
negative counterparts.
 
This mechanism eventually stops giving the required domination. In the
central region, several alphabet-size terms are comparable. An estimate
that replaces all of them by their absolute values loses the cancellation
that determines the sign. We therefore change representations rather than
force the edge estimate farther than it can go. The substitution $t=-\ee^{-z}$ turns the negative real axis into the real
$z$-axis and turns reciprocity into reflection. A centered expansion of the
order polynomial expresses the transformed numerator as a finite linear
combination of derivatives of
\[
 f(z)=\frac{1}{2\cosh(z/2)}.
\]
 
Its poles are known: they lie at the odd imaginary multiples of $\pi$.
A conjugate pair of simple poles contributes a sine after differentiation.
The closest pair therefore suggests a grid of alternating extrema on the
real axis. This is the geometric origin of the central sign argument.
 
The sine is a model, not an identity for the whole polynomial. Indeed, it
has two more interior zeros than the degree permits, and its endpoint signs
are eventually wrong at every fixed size. The proof must retain the
corrections. We keep the first three derivative terms together, where their
phases can reinforce the principal sine, and divide the discrepancy into
three parts: the phase deficit of that retained group, its farther poles,
and the remaining coefficient corrections. The resulting bounds are
\begin{equation}\label{eq:intro-budget}
 \Lphase<\frac{19}{100},\qquad
 \Lfar<\frac{13}{1000},\qquad
 \Lhigh<\frac34.
\end{equation}
Their sum is less than one by $47/1000$. Section~\ref{sec:geometry} explains
the exact decomposition before Section~\ref{sec:central} proves the bounds. There is a finite-size issue as well as a spatial one. Put
$N=\Dg+1$ and $K_N=\lfloor\sqrt{8N/25}-1/2\rfloor$. The uniform central
window begins at $\Dg=11$, whereas the nonempty edge construction begins at
$\Dg=19$. The division of work is therefore as follows.
\begin{center}
\begin{tabular}{@{}cl@{}}\toprule
Size & Root certification\\\midrule
$\Dg=1,2$ & Constant normalized numerator.\\
$3\le\Dg\le10$ & Exact reciprocal reductions and integer sign brackets.\\
$11\le\Dg\le18$ & Central grid and its two exterior intervals.\\
$\Dg\ge19$ & Near-zero intervals, central intervals, reciprocal far intervals.\\
\bottomrule\end{tabular}
\end{center}
The finite cases are not a third spatial group of zeros. They are the sizes
for which a short exact calculation is preferable to extending a uniform
analytic estimate.
 
One picture runs through all three parts: an argument-principle count on a
contour. In part~(i) the contours are the circles $|t|=a_{\Dg,j}$ in the
polynomial coordinate and the boundaries of the annuli between them.
Rouch\'e's theorem (Theorem~\ref{thm:rouche}) fixes the count inside each circle,
so each annulus holds one zero, and the intermediate value theorem
(Theorem~\ref{thm:ivt}) places that zero on the negative axis. In
part~(iii) they are the squares of Proposition~\ref{prop:poles} around the poles of the
kernel, which show that a finite sum of pole contributions approximates the
kernel derivative to controlled accuracy. In both cases an explicit
dominant term fixes the count, and the work of the proof is to
show that the dominance is quantitative and uniform in the size.
 
\subsection{The literature behind the method}
 
The three parts above draw on several established theories, and since each
returns later in the paper, we say here what is borrowed, what is new, and
where each piece is used.
 
The construction starts from the quotient of paths by position-visit
counts, the organizing object of the enumeration framework of Torres,
Sreenivas, and Sowers \cite{TSS}. At the start of
Section~\ref{sec:locality} we explain why this quotient is not the object of
classical lattice-path enumeration, with its continued fractions and
weighted paths \cite{FlajoletCF,Viennot,FlajoletGuillemin,BanderierFlajolet}:
there, different temporal trajectories are counted separately, while here
they may represent a single object.
The target of the construction, the fence chain polytope, is also
classical \cite{StanleyPolytopes}. The explicit source-side
construction in Section~\ref{sec:locality} explains how it arises from
position-visit profiles.
 
The threshold matrix and its trigonometric eigenvalues are part of the
unit-primitive-matrix and Kekul\'e-number theory studied by Xin and Zhong
\cite[Theorem~1.6, Lemma~4.9]{XinZhong}; we recover it in
Proposition~\ref{prop:spectrum}. We give a short derivation because
both the boundary weights and the signs of the eigenvalues are needed, and
because the cumulative-sum structure explains the choice of eigenvectors.
The contribution of the spectral calculation here is its uniform use in a
root-localization argument, not a claim that the underlying spectrum is new.
 
There is also a continuous spectral antecedent. Elkies diagonalizes an
integral operator with triangular support in his study of alternating
permutations and reciprocal-power sums \cite[Section~3, Lemma~3]{Elkies}.
Section~\ref{sec:edges} (Remark~\ref{rem:elkies-kernel}) makes its rescaling and its relation to our finite
spectrum explicit.
 
The two classical results that turn spectral domination into root counts,
Rouch\'e's theorem \cite[Chapter~4, Section~5.2]{Ahlfors} and the
intermediate value theorem \cite[Theorem~4.23]{RudinPMA}, are stated in
Section~\ref{sec:classical-tools} in the form in which they are used,
together with the annulus contours to which they are applied.
At the transition to the central part, the argument uses Petersen
and Zhuang's fixed-alphabet recurrence \cite[Theorem~3.15]{PZ}; its
conversion into a centered differential expansion is carried
out explicitly in Section~\ref{sec:centered}.
 
The differential representation also has a substantial literature behind
it. Hoffman's derivative polynomials for tangent and secant
\cite{Hoffman1995,Hoffman1999}, Boyadzhiev's hyperbolic formulas
\cite[Section~4]{Boyadzhiev}, and the real-zero results associated with
Hetyei's transforms \cite{Hetyei}, recorded in the derivative-polynomial
notation of \cite[Lemma~3.1]{Ma}, motivate the use of a hyperbolic kernel.
They do not imply that an arbitrary size-dependent combination of its
derivatives is real-rooted. We explain this distinction when the exact
combination is obtained, in Section~\ref{sec:kernel-transform}.
 
The pole viewpoint belongs to the classical relation between the zeros of
high derivatives and the geometry of singularities. Its coefficient-extraction
form is the meromorphic expansion method of
\cite[Theorem~IV.10 and the ensuing effective error bound]{FlajoletSedgewick},
and Section~\ref{sec:logarithm} and Proposition~\ref{prop:poles} use it in
that form. A recent geometric
treatment is \cite[Section~1.3]{BSTW}. Here one needs a finite-size statement
for a \emph{linear combination} of derivatives, with coefficients depending
on the size. The central estimates of Section~\ref{sec:central} address precisely that additional
requirement: neither an asymptotic nearest-pole principle nor numerical
agreement with a sine suffices for the sign certification used below.
 
\subsection{Important note on notation}
 
We reserve $d$ for a chronological level and $D$ for path depth. The
polynomial index and gap dimension are $\Dg$; in the path correspondence,
$\Dg=D-1$. The jump radius is $r$, and its associated alphabet size is
$M=r+1$. Outside that correspondence, $M$ varies independently. Position-visit
counts are denoted by $\nu_q$, leaving $c_k(\Dg)$ for the correction
coefficients of the analytic expansion. The variable $t$ is the polynomial
coordinate, $z$ the complex logarithmic coordinate, and $s$ its real
restriction. The indices $j,k,p$ refer, respectively, to pole pairs,
correction terms, and selected extrema.
Positions of a path are $a_i$ and its gaps $g_i$; the coefficients of
$h_{\Dg}$ are $\mathsf h_j$; the Euler-ratio constants are $\varkappa_\ell$; the
reduction polynomial of Section~\ref{sec:finite} is $\Upsilon$; the orthogonal
eigenvector matrix is $\Xi_M$; $\theta_j$ is always the argument of the vector from the
$j$th conjugate pole pair to the evaluation point, with $\theta_0$ the
principal pair and $\theta_p$ its value at the $p$th selected extremum; and
$\Dg^{\star}_k$ is the first size at which the order-$k$ correction cap is used.
The number of zeros of $h$ in $|t|<a$, written as a contour integral, is
$\mathcal I_h(a)$; a zero is $\rho$ and its multiplicity $m_\rho$.
A symbol local to one proof is introduced there and not reused.
 
\section{Position-visit profiles, locality, and chain polytopes}\label{sec:locality}
 
Our change of viewpoint, described at the start of
Section~\ref{sec:three-parts}, also tells us that our starting question is
not initially about zeros. Instead, we want to ask:
\begin{quote}
Which visit profiles can a path with bounded integer jumps actually
produce once we forget the order of its steps, and how can we count them?
\end{quote}
\noindent Consider a path with
bounded integer jumps in a recombining tree, and record how many times it
visits each position. If a node's value depends only on position, these visit
counts determine the sum of the node values along the path. Different orders
of the same visits then become indistinguishable. This quotient of paths by
position-visit counts is the organizing object of the enumeration framework
of Torres, Sreenivas, and Sowers \cite{TSS}.
 
The distinction from ordinary walk enumeration is essential. Continued
fractions and weighted lattice paths are classical companions in
Flajolet's theory \cite{FlajoletCF} and in Viennot's combinatorial theory of
orthogonal polynomials \cite{Viennot}; related path expansions underlie the
birth-and-death analysis of Flajolet and Guillemin
\cite{FlajoletGuillemin}. Banderier and Flajolet develop the analytic
enumeration of directed lattice paths for general finite step sets
\cite{BanderierFlajolet}. Those theories provide the natural background,
but here different temporal trajectories may represent a single object.
The local criterion below must therefore characterize the quotient by visit
counts, rather than count trajectories with their multiplicities.
 
The temporal order of a path is a large object, since every level offers a
choice of jump. A position-visit profile is much smaller, but it is not an
arbitrary vector of nonnegative integers, so the task in this section is to
identify exactly what geometric information survives once temporal order has
been discarded. Sorting the visits exposes it as a nearest-neighbor
constraint. For an origin-to-origin path, remove the final occurrence of the
origin and sort the remaining occurrences. What survives is local:
\emph{the sum of two consecutive sorted gaps cannot exceed the jump radius}.
Necessity is a cut-crossing argument and sufficiency a constructive cyclic
ordering of the sorted occurrences, so a global existence problem for a path
reduces to inequalities involving only neighboring gaps.

This observation has two simultaneous consequences. The gap vectors are
lattice points in a dilation of a fence chain polytope, and the compatibility
of two successive gaps is encoded by a $0$--$1$ matrix. The polytope
identifies the enumerative polynomial and explains why the path quotient
leads to this particular family, so it is not an ornamental interpretation
added after the proof; the matrix supplies an analytic way to study that
polynomial. We keep both descriptions because their roles are different.
We begin our proof by formalizing these observations.

\subsection{The path quotient and its cyclic form}

Fix a jump radius $r\ge0$ and a depth $D\ge1$. The recombining
$(2r+1)$-nomial tree has vertices $(x,d)$, where $d\in\{0,\ldots,D\}$ and
$x\in\ZZ$ is reachable from the origin, and directed edges
\[
 ((x,d-1),(x+u,d)),\qquad 1\le d\le D,\quad -r\le u\le r.
\]
Only the position sequence is needed below. An origin-to-origin path is
\begin{equation}\label{eq:path-definition}
 \mathsf p=(x_0,\ldots,x_D),\qquad
 x_0=x_D=0,\qquad |x_d-x_{d-1}|\le r.
\end{equation}
Its position-visit profile, also called its cardinality tuple, is
\begin{equation}\label{eq:profile}
 \nu(\mathsf p)=(\nu_q)_{q\in\ZZ},\qquad
 \nu_q=\#\{d:0\le d\le D,\ x_d=q\}.
\end{equation}
It has finite support, total mass $D+1$, and $\nu_0\ge2$. Two paths are equivalent when their profiles agree. Thus profiles count classes, not temporal realizations. In particular,
$\sum_{d=0}^{D}v(x_d)=\sum_q\nu_qv(q)$ for every function $v$ of position; this is the reason for taking the quotient in the first place. Removing the final zero gives the cyclic profile
\begin{equation}\label{eq:cyclic-profile}
 \nu_q^\circ=\nu_q-\ind_{\{q=0\}},\qquad
 \sum_q\nu_q^\circ=D.
\end{equation}
The occurrences $x_0,\ldots,x_{D-1}$ form a cycle, with the final edge $x_{D-1}\to x_0$ supplied by the last step of the path. Occurrences at the same position remain distinct occurrences. This distinction matters when we remove one occurrence from a cycle in the necessity proof. List the multiset specified by $\nu^\circ$ in nondecreasing order:
\begin{equation}\label{eq:sorted-occurrences}
 a_1\le a_2\le\cdots\le a_D,
 \qquad
 g_i=a_{i+1}-a_i\quad(1\le i<D),\qquad g_0=g_D=0.
\end{equation}
The zeros at the ends of the gap sequence are boundary conventions, not
additional occurrences. They allow the endpoint inequalities to be written
in the same form as the interior ones.

\begin{theorem}[Locality criterion]\label{thm:locality}
Let $\nu$ be a finitely supported nonnegative integer vector with
$\sum_q\nu_q=D+1$ and $\nu_0\ge2$. Form \eqref{eq:cyclic-profile} and
\eqref{eq:sorted-occurrences}. Then $\nu$ is the profile of a path
\eqref{eq:path-definition} if and only if
\begin{equation}\label{eq:locality}
 g_i+g_{i+1}\le r,\qquad 0\le i\le D-1.
\end{equation}
\end{theorem}

\begin{proof}
Suppose first that a path realizes $\nu$. Label its $D$ cyclic occurrences
by $1,\ldots,D$ so that their positions are $a_1,\ldots,a_D$ in sorted
order; ties can be resolved arbitrarily. Every edge of this labeled cycle
joins positions at distance at most $r$.

Each single gap satisfies $g_i\le r$. Indeed, if $g_i>r$, no allowed edge
could join an occurrence among $1,\ldots,i$ to an occurrence among
$i+1,\ldots,D$, whereas the cyclic ordering must cross that cut. This proves
the two endpoint inequalities in \eqref{eq:locality}. For an interior index
$1\le i\le D-2$, remove occurrence $i+1$ and its two incident edges. The
remaining occurrences, in their inherited cyclic order, form a path. That
path contains both the nonempty set $\{1,\ldots,i\}$ and the nonempty set
$\{i+2,\ldots,D\}$. Somewhere along it an edge crosses between these two
sets. The distance across that edge is at least
\[
 a_{i+2}-a_i=g_i+g_{i+1},
\]
and at most $r$. This proves every interior inequality. When $D=1$ there are no nonzero gaps; when $D=2$, the single-gap argument is the whole necessity proof. Conversely, assume \eqref{eq:locality}. Visit the odd-numbered sorted occurrences in increasing order, then the even-numbered ones in decreasing order, and close the cycle. For example, when $D=6$ the order is $1,3,5,6,4,2,1,$ and when $D=5$ it is $1,3,5,4,2,1$. Every consecutive pair of indices in this cyclic order differs by one or two. A difference of two gives a position distance $a_{i+2}-a_i=g_i+g_{i+1}\le r$. A difference of one gives $g_i\le r$, which follows from \eqref{eq:locality} because all gaps are nonnegative. For $D=1$, the cycle is the single zero loop. For $D=2$, it traverses the same allowed pair in both directions. At least one cyclic occurrence is at zero because $\nu_0^\circ\ge1$. Rotate the constructed cyclic order so that this occurrence is first, then append it as the final occurrence. The resulting path starts and ends at
zero, has $D$ allowed jumps, and has profile $\nu^\circ+\ind_{\{q=0\}}=\nu$. This proves sufficiency.
\end{proof}

The proof explains why adjacent sums, rather than individual gap bounds, are the right condition. A single-gap bound controls a cut. Removing one
occurrence makes two neighboring gaps visible to the same cut, and the
cyclic ordering forces an edge across it. The converse construction uses
exactly those one- and two-gap distances. Necessity and sufficiency
therefore identify the same local obstruction.

\subsection{Anchoring the sorted gaps}

For general profiles the gaps do not record a translation of the support. To obtain a bijection without an additional coordinate, restrict to profiles whose support is nonnegative. Because zero occurs, their first sorted occurrence is then $a_1=0$. The translation is fixed, and the sorted list can be reconstructed by cumulative sums. Set $\Dg=D-1$. For $M\ge1$ define
\begin{equation}\label{eq:gap-set}
 \Gcal_{\Dg}(M)=
 \left\{(g_1,\ldots,g_{\Dg})\in\ZZ_{\ge0}^{\Dg}:
 g_0=g_{\Dg+1}=0,\quad
 g_i+g_{i+1}\le M-1\ (0\le i\le\Dg)\right\}.
\end{equation}
The empty vector is the unique member of $\Gcal_0(M)$. The boundary
inequalities imply $0\le g_i\le M-1$ for all coordinates.

\begin{proposition}[Profile--gap correspondence]\label{prop:profile-gap}
For $D\ge1$ and $r\ge0$, the distinct nonnegative-support profiles of
origin-to-origin paths with jump radius $r$ are in bijection with
$\Gcal_{D-1}(r+1)$.
\end{proposition}

\begin{proof}
The forward map is \eqref{eq:cyclic-profile}--\eqref{eq:sorted-occurrences};
Theorem~\ref{thm:locality} puts its image in the indicated gap set. Conversely,
for $g\in\Gcal_{D-1}(r+1)$ define
\begin{equation}\label{eq:gap-inverse}
 a_1=0,\qquad a_i=\sum_{j=1}^{i-1}g_j\quad(2\le i\le D),
 \qquad
 \nu_q=\#\{i:a_i=q\}+\ind_{\{q=0\}}.
\end{equation}
This profile has the required mass and at least two occurrences of zero.
Its support is nonnegative, and Theorem~\ref{thm:locality} constructs a
realizing path. Sorting its cyclic occurrences recovers the same list
$a_1,\ldots,a_D$, hence the same gaps. In the other direction, $a_1=0$
ensures that cumulative sums reconstruct the original sorted list.
Thus the two maps are inverse on profiles.
\end{proof}

\begin{example}\label{ex:path-gap}
At $D=4$ and $r=2$, the path $(0,1,3,2,0)$ has cyclic occurrences
$(0,1,3,2)$, sorted occurrences $(0,1,2,3)$, and gap vector $(1,1,1)$.
The two interior inequalities are $1+1\le2$, and the endpoint inequalities
are $1\le2$. The profile is $(\nu_0,\nu_1,\nu_2,\nu_3)=(2,1,1,1)$.
The correspondence records this profile once, regardless of how many
allowed temporal orderings realize it.
\end{example}

\begin{figure}[h]
 \centering
 \includegraphics[width=\linewidth]{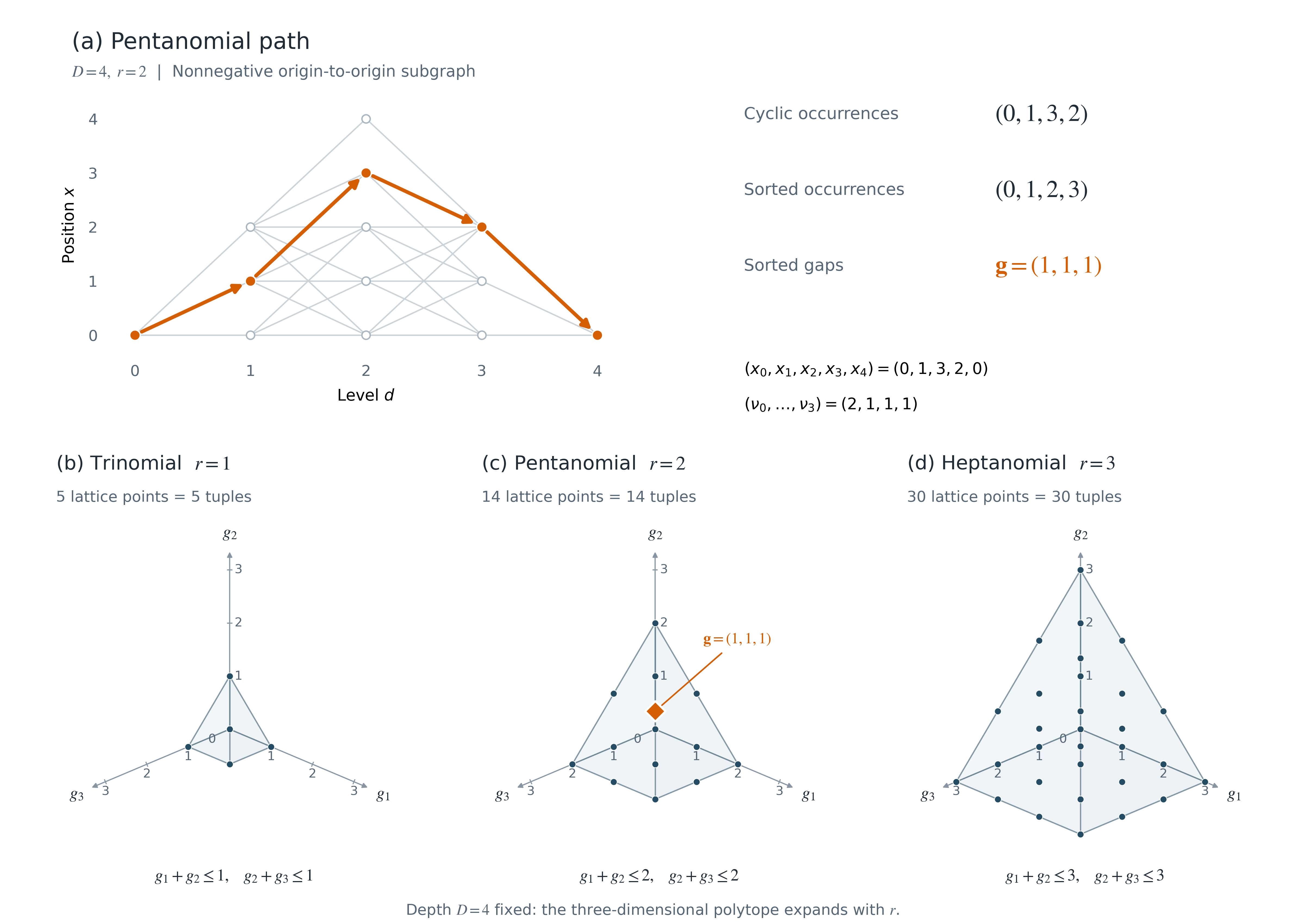}
 \caption{A nonnegative origin-to-origin path and its sorted-gap image.
 The highlighted pentanomial path in Example~\ref{ex:path-gap} represents
 the profile associated with $(1,1,1)$. At fixed depth $D=4$, the three
 lower panels show the same three-dimensional chain polytope dilated by
 $r=1,2,3$. Their $5,14,30$ lattice points count distinct profiles, not
 individual paths. The displayed path is one representative of its class.}
 \label{fig:path-polytope}
\end{figure}

\subsection{The fence chain polytope}

We now interpret the inequalities themselves. This identifies the classical
object reached by the path construction, while preserving the distinction
between that object and the profiles that led to it.

\begin{definition}\label{def:fence}
The fence poset $P_{\Dg}$ has elements $p_1,\ldots,p_{\Dg}$ and cover
relations $p_i<p_{i+1}$ for odd $i$, and $p_i>p_{i+1}$ for even $i$.
The chain polytope of a finite poset $P$ is
\begin{equation}\label{eq:chain-polytope}
 \mathcal C(P)=\left\{x\in\RR_{\ge0}^{P}:
 \sum_{p\in C}x_p\le1\text{ for every chain }C\subseteq P\right\}.
\end{equation}
\end{definition}

Definition~\ref{def:fence} uses Stanley's chain-polytope convention
\cite[Definition~2.1]{StanleyPolytopes}. Every nontrivial chain in a fence
has two elements, at neighboring indices. For $\Dg\ge2$, each element is
in such a chain, so the inequalities become
\[
 x_i\ge0,\qquad x_i+x_{i+1}\le1\quad(1\le i<\Dg).
\]
The singleton bounds follow from nonnegativity and a neighboring pair.
For $\Dg=1$, the singleton condition $x_1\le1$ must be included explicitly.
Equivalently, for all $\Dg\ge1$ we may pad with $x_0=x_{\Dg+1}=0$ and use
the consecutive inequalities for $0\le i\le\Dg$.
Consequently,
\begin{equation}\label{eq:gap-polytope}
 \Gcal_{\Dg}(M)
 =(M-1)\mathcal C(P_{\Dg})\cap\ZZ^{\Dg}.
\end{equation}
There is no change of lattice in this identity: the gap coordinates are the polytope coordinates. For completeness, the direct passage to order maps is particularly simple for a fence. Given $g\in\Gcal_{\Dg}(M)$, put
\begin{equation}\label{eq:gap-order-map}
 y_i=\begin{cases}
       g_i+1,&i\text{ odd},\\
       M-g_i,&i\text{ even}.
      \end{cases}
\end{equation}
Then $1\le y_i\le M$. When $i$ is odd, $g_i+g_{i+1}\le M-1$ is equivalent
to $y_i\le y_{i+1}$; when $i$ is even, it is equivalent to
$y_i\ge y_{i+1}$. The inverse formulas are $g_i=y_i-1$ at odd indices
and $g_i=M-y_i$ at even ones. Thus
\begin{equation}\label{eq:omega-gap}
 \#\Gcal_{\Dg}(M)=\Omega_{\Dg}(M).
\end{equation}
This explicitly realizes, in the fence case and with the shift $M-1$, the
Ehrhart/order-polynomial relationship in
\cite[Theorem~4.1]{StanleyPolytopes}. It also identifies the count in
Proposition~\ref{prop:profile-gap} with the numerator family in
\eqref{eq:intro-numerator}.

\subsection{Why a transfer matrix is already present}

A gap word can be built one coordinate at a time. Once its current gap is $i$, the next gap $j$ is allowed precisely when $i+j\le M-1$. Nothing else about the earlier word is needed. This locality is the reason a matrix with states $0,\ldots,M-1$ will count the words: its entries are the indicators of these allowed pairs. The zero gaps at the two boundaries specify the
initial and final state. Before using the matrix analytically, we establish the numerator's algebraic properties. Those properties perform two later tasks. Reciprocity reflects the near-zero roots to the far negative side, and the degree supplies the final check that the interval count is complete. Neither fact should enter the proof as an unexplained feature of the numerical examples.

% ===== 03_dictionary =====
\section{The zigzag Eulerian dictionary}\label{sec:dictionary}

The gap model gives a count at every alphabet size. To identify its numerator,
we use the linear extensions of the fence. This section records the exact
normalizations and proves the properties later used in the analytic argument.
It is not necessary to develop the full theory of $P$-partitions for this
purpose.

\subsection{Linear extensions and the numerator}

A labeling $\omega:P\to\{1,\ldots,|P|\}$ is \emph{natural} when it is a
bijection and $p<q$ implies $\omega(p)<\omega(q)$. For the fence, one choice
labels all odd-indexed elements first, followed by the even-indexed ones:
\[
 \omega(p_i)=\begin{cases}
 (i+1)/2,&i\text{ odd},\\
 \lceil\Dg/2\rceil+i/2,&i\text{ even}.
 \end{cases}
\]
Every cover runs from an odd index to an even one, so this is natural.
A linear extension is a listing of the elements consistent with the poset
order. Read its labels as a permutation $w=w_1\cdots w_{\Dg}$, and let
$\des(w)=\#\{i:w_i>w_{i+1}\}$.

The following proposition is the fundamental $P$-partition decomposition
specialized to a naturally labeled fence. The general decomposition is
Stanley's \cite[Lemma~3.15.3 and Theorem~3.15.8]{StanleyEC1}; its original
order-polynomial form is \cite[Proposition~13.3]{StanleyPartitions}.
We give the specialization because the direction of the order maps and the
factor of $t$ in the numerator determine the later root normalization.

\begin{proposition}[Order count and descent numerator]\label{prop:numerator}
For every $\Dg\ge1$,
\begin{align}
 \Omega_{\Dg}(M)
 &=\sum_{w\in\mathcal L(P_{\Dg},\omega)}
   \binom{M+\Dg-1-\des(w)}{\Dg},\label{eq:order-binomial}\\
 h_{\Dg}(t)&=\sum_{w\in\mathcal L(P_{\Dg},\omega)}t^{\des(w)},
 \qquad
 F_{\Dg}(t)\Def\sum_{M\ge1}\Omega_{\Dg}(M)t^{M-1}
 =\frac{h_{\Dg}(t)}{(1-t)^{\Dg+1}}.\label{eq:numerator}
\end{align}
In particular, $h_{\Dg}$ has nonnegative integer coefficients and
$h_{\Dg}(0)=1$. The polynomial is independent of the chosen natural labeling.
\end{proposition}

\begin{proof}
Take a weak order-preserving map from the fence to $\{1,\ldots,M\}$.
Sort the elements first by their values and then, for ties, by increasing
natural label. The resulting word is a linear extension: a comparable pair
appears in the right order whether its values differ or agree. Along a fixed
extension $w$, the sorted values satisfy
\[
 1\le b_1\le\cdots\le b_{\Dg}\le M,
 \qquad b_i<b_{i+1}\text{ whenever }w_i>w_{i+1}.
\]
Conversely, these conditions recover the same tie-breaking order and define a weak order-preserving map. If $a=\des(w)$, subtract from $b_i$ the number of descents preceding position $i$. This gives a weakly increasing sequence of length $\Dg$ in $\{1,\ldots,M-a\}$, counted by
$\binom{M-a+\Dg-1}{\Dg}$. This proves \eqref{eq:order-binomial}, including zero counts when $M\le a$. For $0\le a\le\Dg-1$, the binomial generating-function identity gives
\[
 \sum_{M\ge1}\binom{M+\Dg-1-a}{\Dg}t^{M-1}
 =\frac{t^a}{(1-t)^{\Dg+1}}.
\]
Indeed, the terms $M\le a$ vanish, and putting $q=M-1-a\ge0$ reduces the
sum to $t^a\sum_{q\ge0}\binom{q+\Dg}{\Dg}t^q$.
Summing over extensions proves \eqref{eq:numerator}. The only descent-free
permutation is $1\cdots\Dg$, and it is an extension under a natural
labeling, so the constant coefficient is one. Finally, the left side of
\eqref{eq:numerator} is an unlabeled order-map count; multiplication by
$(1-t)^{\Dg+1}$ proves independence of the labeling.
\end{proof}

This is the normalization of the fundamental order-polynomial identity in \cite[Eq.~(1), Lemma~2.2, Eqs.~(5)--(6)]{PZ}: their generating series uses $t^M$, whereas ours uses $t^{M-1}$. Thus their $Z_{\Dg}(t)$ is exactly $t h_{\Dg}(t)$, not $h_{\Dg}$ itself. This factor accounts for the zero at the origin and is kept separate throughout the root analysis. The polyhedral normalization is equally explicit. The order polytope is
\[
 \mathcal O(P_{\Dg})=\{x\in[0,1]^{\Dg}:x_i\le x_j
                         \text{ whenever }p_i<p_j\}.
\]
Adding one to every coordinate gives a bijection from its lattice points
at dilation $r$ to the order maps counted by $\Omega_{\Dg}(r+1)$.
Together with \cref{eq:gap-polytope,eq:omega-gap,eq:numerator}, this gives
\begin{equation}\label{eq:ehrhart-numerator}
 \sum_{r\ge0}\#\bigl(r\mathcal C(P_{\Dg})\cap\ZZ^{\Dg}\bigr)t^r
 =\sum_{r\ge0}\#\bigl(r\mathcal O(P_{\Dg})\cap\ZZ^{\Dg}\bigr)t^r
 =\frac{h_{\Dg}(t)}{(1-t)^{\Dg+1}}.
\end{equation}
The numerator of this lattice-point series is called the Ehrhart
$h^*$-polynomial. Thus our $h_{\Dg}$ is exactly the $h^*$-polynomial studied
by Coons and Sullivant, with no further power of $t$ or change of variable.
Their swap enumerator \cite[Theorem~1.9 and Corollary~3.5]{CoonsSullivant}
and the descent enumerator in \cref{eq:numerator} therefore describe the
same polynomial. We use the latter because it makes the degree and
reciprocity arguments particularly direct.

\subsection{Reciprocity and the exact degree}

Let $\Omega^\circ_{\Dg}(M)$ count maps that are strictly increasing on
comparable pairs. Sorting such a map with \emph{decreasing} labels within
ties gives a linear extension. Along that extension, strict increases are
required at its ascents, of which there are $\Dg-1-\des(w)$. Conversely,
these requirements enforce strictness on comparable pairs: their labels
increase, so the intervening segment cannot be entirely descending. The
same subtraction argument as above gives
\begin{equation}\label{eq:strict-binomial}
 \Omega^\circ_{\Dg}(M)
 =\sum_w\binom{M+\des(w)}{\Dg},
 \qquad
 \Omega_{\Dg}(-M)=(-1)^{\Dg}\Omega^\circ_{\Dg}(M).
\end{equation}
The second identity uses
$\binom{-M+\Dg-1-a}{\Dg}=(-1)^{\Dg}\binom{M+a}{\Dg}$ and is a polynomial identity. This is precisely Stanley's order-polynomial reciprocity
\cite[Proposition~13.2(i)]{StanleyPartitions},
\cite[Proposition~2.1]{StanleyReciprocity}; see also
\cite[Corollary~3.15.12]{StanleyEC1}. The direct binomial verification here
fixes the strict-map convention rather than replacing that attribution.

For $\Dg\ge2$, every element of the fence is incident with a cover and its
rank is zero or one. Subtracting one from the values at all maximal
elements takes a strict map into $\{1,\ldots,M\}$ to a weak map into
$\{1,\ldots,M-1\}$. The bounds are legitimate: each minimal value is at
most $M-1$, and each maximal value is at least two. Adding one back at
maximal elements is the inverse. Hence
\begin{equation}\label{eq:order-reciprocity}
 \Omega^\circ_{\Dg}(M)=\Omega_{\Dg}(M-1),
 \qquad
 \Omega_{\Dg}(-1-M)=(-1)^{\Dg}\Omega_{\Dg}(M)
 \quad(\Dg\ge2).
\end{equation}
The latter identity is the rank-one case of the graded-poset identity in
\cite[Corollary~3.15.18(b)(iii)]{StanleyEC1}. In that result the longest
chain has one edge, so the shift is one. This accounts for the center
$M=-1/2$ of the later expansion; it is not a centering chosen from a fit to
the coefficients.

\begin{proposition}[Reciprocal numerator]\label{prop:reciprocity}
For $\Dg\ge2$, $h_{\Dg}$ is monic of degree $\Dg-2$ and
\begin{equation}\label{eq:h-reciprocity}
 h_{\Dg}(t)=t^{\Dg-2}h_{\Dg}(1/t).
\end{equation}
Moreover, $h_1=h_2=1$.
\end{proposition}

\begin{proof}
Using \eqref{eq:strict-binomial} and the same binomial summation as in
Proposition~\ref{prop:numerator},
\[
 \sum_{M\ge1}\Omega^\circ_{\Dg}(M)t^M
 =\frac{t^{\Dg}h_{\Dg}(1/t)}{(1-t)^{\Dg+1}}.
\]
On the other hand, \eqref{eq:order-reciprocity}, together with
$\Omega_{\Dg}(0)=0$, makes this series
$t^2h_{\Dg}(t)/(1-t)^{\Dg+1}$. Equating the two expressions gives
\eqref{eq:h-reciprocity}. The numerator initially has degree at most
$\Dg-1$ by its descent definition. The identity forces the coefficient
at degree $\Dg-1$ to vanish, and pairs its constant coefficient one with
its coefficient at degree $\Dg-2$. Thus its degree is exactly $\Dg-2$
and its leading coefficient is one. For $\Dg=1$, the count is $M$, so
$F_1(t)=(1-t)^{-2}$ and $h_1=1$; the case $\Dg=2$ follows either from the
argument or from the unique extension of a two-element chain.
\end{proof}

Reciprocity preserves multiplicities at nonzero roots. More explicitly,
if $(t-\rho)^q$ is the exact vanishing factor at $\rho\ne0$, substitution
$t\mapsto1/t$ in \eqref{eq:h-reciprocity} gives the same order at
$\rho^{-1}$. We will use this to reflect the edge root count, not to assume
any reality before it has been proved.

\subsection{Euler numbers and the coefficient consequences to come}

Let $E_{\Dg}$ denote the number of alternating permutations of size $\Dg$,
with $E_0=E_1=1$. Assign to each element of the fence its position in a
linear extension. These ranks, read along the fence, form an alternating
permutation, and the inverse assignment reconstructs the extension. Thus
\begin{equation}\label{eq:euler-normalization}
 h_{\Dg}(1)=E_{\Dg},
 \qquad
 \sum_{q\ge0}E_q\frac{u^q}{q!}=\sec u+\tan u.
\end{equation}
The second identity is Andr\'e's classical exponential generating function;
we use the normalization in \cite[Theorem~1.1, Eqs.~(1.2)--(1.3)]{StanleyAlt}.
Complementing each entry $i\mapsto\Dg+1-i$ interchanges the up--down and
down--up conventions, so their counts and this normalization agree.

There is a second convention for ``Euler numbers,'' in which the exponential
generating function is $\sech u$ rather than $\sec u+\tan u$.
Sagan \cite[Section~1, Eqs.~(1)--(3)]{Sagan} explicitly distinguishes these
conventions in his treatment of generalized Euler numbers. If
$\mathcal E_q=[u^q]\sech(u)\,q!$, then
$\mathcal E_{2j}=(-1)^jE_{2j}$ and $\mathcal E_{2j+1}=0$.
We always use the positive $E_q$ in \eqref{eq:euler-normalization},
including odd indices. This matters in the ratio estimates: a sign from a
kernel residue must not inadvertently be included a second time in $E_q$.

Petersen and Zhuang prove the stronger structural statement that
$Z_{\Dg}(t)/t=h_{\Dg}(t)$ is $\gamma$-nonnegative
\cite[Theorem~2.5]{PZ}, by applying Br\"and\'en's sign-graded-poset theorem
\cite[Theorem~4.2]{BrandenSign}. To make the parameters explicit, a natural
labeling assigns positive sign to each cover, and every maximal chain of
our fence has one edge. The rank in that theorem is therefore $r=1$ and
the number of elements is $p=\Dg$. Its exponent $p-r-1$ becomes $\Dg-2$,
so for $\Dg\ge2$ the conclusion is
\begin{equation}\label{eq:known-gamma-expansion}
 h_{\Dg}(t)=\sum_{j=0}^{\lfloor(\Dg-2)/2\rfloor}
 \gamma_{\Dg,j}t^j(1+t)^{\Dg-2-2j},\qquad \gamma_{\Dg,j}\ge0.
\end{equation}
This explains the already-known symmetry and unimodality. We do not use it as a substitute for a root argument:
nonnegative coefficients in a $\gamma$-expansion do not alone establish
real-rootedness. Instead, the exact degree in
Proposition~\ref{prop:reciprocity} will be exhausted by real intervals.
Once that is done, Newton's inequalities will strengthen the coefficient
conclusion to \eqref{eq:intro-newton}.

\begin{remark}[Classical interpretations of the same local inequalities]
\label{rem:classical-models}
The inequalities in the gap model also give a transparent graph-labeling
interpretation. On the path with vertices $0,\ldots,\Dg$, label the edge
$\{i-1,i\}$ by $g_i$ and the loop at vertex $i$ by
$r-g_i-g_{i+1}$, where $g_0=g_{\Dg+1}=0$. The labels are nonnegative
and the sum incident to each vertex, counting a loop once, is $r$.
Conversely, the non-loop labels recover the gaps. This is a magic labeling
in the convention of Stanley \cite[Section~1]{StanleyMagic}; the specific
zigzag correspondence is recorded in \cite[Section~5]{PZ}.

Equivalently, put $A_{i-1,i}=A_{i,i-1}=g_i$ and
$A_{ii}=r-g_i-g_{i+1}$, with all other entries zero. Then $A$ is a
symmetric nonnegative integer tridiagonal matrix with every row and column
sum equal to $r$. For $r>0$, $A/r$ lies in the tridiagonal doubly stochastic
polytope studied by Diaconis and Wood \cite{DiaconisWood}. The converse is
also exact for integer matrices: equality of the row and column sums at
vertex zero forces $A_{01}=A_{10}$, and the same comparison successively
at vertices $1,2,\ldots$ forces symmetry on every adjacent pair.
Finally, the alternating-height interpretation of Kekul\'e structures in
\cite[Section~5]{PZ} connects this count with the chemical-graph setting of
Cyvin and Gutman \cite{CyvinGutman} and Xin and Zhong \cite{XinZhong}.
These are existing realizations of the count. The profile-locality theorem
supplies the source-side map needed here, rather than a new claim about the
older target objects.
\end{remark}

The dictionary is now complete. The same count describes position-visit
profiles, lattice points of a chain polytope, and order-preserving maps.
We next return to its local gap description, since this is the representation
in which spectral information is directly accessible.

% ===== 04_edges =====
\section{Threshold spectra and the extreme zeros}\label{sec:edges}

The local inequality of Theorem~\ref{thm:locality} now becomes a linear
operator. Its spectrum will first describe the growth of
$\Omega_{\Dg}(M)$ at fixed $M$. A separate comparison across $M$ will then
turn that growth information into circles on which one term of the
generating function dominates. Keeping these two steps separate is
important: a dominant matrix eigenvalue does not, by itself, prove
real-rootedness of the numerator.

\subsection{The compatibility matrix and its boundary weights}

For $M\ge1$, index rows and columns by $0,\ldots,M-1$ and define
\begin{equation}\label{eq:threshold}
 T_{M-1}=(\ind_{\{i+j\le M-1\}})_{0\le i,j<M}.
\end{equation}
This is a symmetric matrix. A product of its entries is one exactly when every consecutive pair of the corresponding gap word is allowed. Padding with the two zero boundary gaps therefore gives
\begin{equation}\label{eq:matrix-count}
 \Omega_{\Dg}(M)
 =(T_{M-1}^{\Dg+1})_{00}
 =\sum_{g_1,\ldots,g_{\Dg}=0}^{M-1}
   \prod_{i=0}^{\Dg}(T_{M-1})_{g_i,g_{i+1}}.
\end{equation}
There are $\Dg+1$ compatibility checks, not $\Dg$. This is the transfer-matrix principle in the form of \cite[Section~4.7]{StanleyEC1}: the state is the current gap, and multiplication sums over successive allowed states. Unlike a walk matrix on spatial positions, this matrix operates after temporal paths have been identified by their visit counts. This change of state space is what makes the local criterion an enumeration of profiles. For $\Dg=0$ the formula gives $(T_{M-1})_{00}=1$, consistent with the empty gap word. For example,
\[
 T_2=\begin{pmatrix}1&1&1\\1&1&0\\1&0&0\end{pmatrix},
 \qquad \Omega_2(3)=(T_2^3)_{00}=3+2+1=6.
\]
The three summands count the possible second gaps after a first gap of
$0,1,2$, respectively.

\begin{proposition}[Spectrum and spectral count]\label{prop:spectrum}
Put
\begin{equation}\label{eq:spectral-data}
 \beta_{M,\ell}=\frac{(2\ell-1)\pi}{2M+1},\qquad
 \lambda_{M,\ell}=\frac{(-1)^{\ell-1}}{2\sin(\beta_{M,\ell}/2)},\qquad
 \alpha_{M,\ell}=\frac4{2M+1}\cos^2(\beta_{M,\ell}/2).
\end{equation}
For $1\le\ell\le M$, these are the simple eigenvalues of $T_{M-1}$ and
their weights at the zero boundary state. More precisely,
\begin{equation}\label{eq:spectral-count}
 \Omega_{\Dg}(M)=\sum_{\ell=1}^M\alpha_{M,\ell}\lambda_{M,\ell}^{\Dg+1},
 \qquad \alpha_{M,\ell}>0,\qquad\sum_{\ell=1}^M\alpha_{M,\ell}=1.
\end{equation}
An orthogonal diagonalizing matrix has entries
\begin{equation}\label{eq:eigenvectors}
 (\Xi_M)_{i\ell}=\frac2{\sqrt{2M+1}}
 \sin((M-i)\beta_{M,\ell}),\qquad 0\le i<M,
 \quad1\le\ell\le M.
\end{equation}
\end{proposition}

\begin{proof}
Let $J$ reverse the coordinates. Then $U=T_{M-1}J$ has entries
$U_{ij}=\ind_{\{i\le j\}}$. It is the upper-triangular cumulative-sum
matrix, and its inverse $B$ has ones on the diagonal and minus ones on the
superdiagonal. Thus
\begin{equation}\label{eq:inverse-difference}
 T_{M-1}^{-1}=JB,\qquad
 T_{M-1}^{-2}=B^{\mathsf T}B.
\end{equation}
The latter matrix has diagonal $(1,2,\ldots,2)$ and adjacent off-diagonal entries $-1$; for $M=1$ it is $(1)$. This is the structural benefit of reversing columns: cumulative summation becomes a first difference, and squaring the inverse gives a second difference. Column reversal alone is
not a similarity transformation; diagonalizability of $T_{M-1}$ follows from its symmetry, not from the triangular form of $U$. With auxiliary boundary values $x_{-1}=x_0$ and $x_M=0$, the eigenvalue equation for $B^{\mathsf T}B$ is
\[
 -x_{i-1}+2x_i-x_{i+1}=\mu x_i,\qquad 0\le i<M.
\]
The trial solution $x_i=\cos((i+1/2)\beta)$ satisfies the first boundary
condition and has $\mu=4\sin^2(\beta/2)$. The second condition is
$\cos((M+1/2)\beta)=0$, giving exactly the angles in
\eqref{eq:spectral-data}. These $M$ eigenvalues are positive and distinct,
since their half-angles increase strictly in $(0,\pi/2)$. We have therefore
constructed a full eigenbasis, without needing an additional classification
of solutions to the recurrence.

The matrix $T_{M-1}$ commutes with its inverse square, so it preserves each
of these one-dimensional eigenspaces. Its last row is $(1,0,\ldots,0)$.
The corresponding eigenvalue is consequently
\[
 \frac{x_0}{x_{M-1}}
 =\frac{\cos(\beta/2)}{(-1)^{\ell-1}\sin\beta}
 =\frac{(-1)^{\ell-1}}{2\sin(\beta/2)}.
\]
Here $(M+1/2)\beta=(\ell-1/2)\pi$, so
$x_{M-1}=\cos((\ell-1/2)\pi-\beta)=(-1)^{\ell-1}\sin\beta$.
The equivalent eigenvector
$v_i=\sin((M-i)\beta)=(-1)^{\ell-1}x_i$ has norm
\[
 \sum_{q=1}^M\sin^2(q\beta)
 =\frac M2-\frac12\sum_{q=1}^M\cos(2q\beta)
 =\frac{2M+1}{4}.
\]
For the last equality, the finite geometric sum gives
$\sum_{q=-M}^M\ee^{2\ii q\beta}
 =\sin((2M+1)\beta)/\sin\beta=0$; hence
$\sum_{q=1}^M\cos(2q\beta)=-1/2$.
Distinct eigenvalues of the real symmetric matrix have orthogonal eigenvectors, so normalization gives \eqref{eq:eigenvectors} and
$T_{M-1}=\Xi_M\diag(\lambda_{M,\ell})\Xi_M^{\mathsf T}$.
Taking the $00$ entry of its $(\Dg+1)$st power proves
\eqref{eq:spectral-count}, because
\[
 (\Xi_M)_{0\ell}^2=\frac4{2M+1}\sin^2(M\beta_{M,\ell})
 =\frac4{2M+1}\cos^2(\beta_{M,\ell}/2).
\]
The weights sum to the squared norm of row zero of $\Xi_M$, namely one.
\end{proof}

This is the spectrum in \cite[Theorem~1.6, Lemma~4.9]{XinZhong}, written in
decreasing modulus. In their indexing the eigenvalues are
$(-1)^{M+1}/[2\cos((2j-1)\pi/(2M+1))]$. The exact permutation is
\[
 j=\begin{cases}
 (M-\ell+2)/2,&M\equiv\ell\pmod2,\\
 (M+\ell+1)/2,&M\not\equiv\ell\pmod2.
 \end{cases}
\]
The cosine then equals $(-1)^{M-\ell}\sin(\beta_{M,\ell}/2)$.
The boundary weights in \eqref{eq:spectral-count} are included here because
the count uses a specified initial and final state, rather than a trace.

\begin{remark}[The continuous triangular kernel]
\label{rem:elkies-kernel}
The cumulative-sum calculation has a continuous analogue in Elkies's
operator \cite[Eq.~(37) and Lemma~3]{Elkies},
\[
 (\mathcal T F)(v)=\int_0^{\pi/2-v}F(u)\dif u,
 \qquad 0<v<\pi/2.
\]
Set $v=\pi x/2$, $u=\pi y/2$, and $F(u)=f(2u/\pi)$.
Then
\[
 (\mathcal T F)(\pi x/2)=\frac\pi2(\mathcal K f)(x),
 \qquad (\mathcal K f)(x)=\int_0^{1-x}f(y)\dif y.
\]
Thus the eigenvalues $1/(4q+1)$, $q\in\ZZ$, of Elkies's operator
become $2/[\pi(4q+1)]$ for $\mathcal K$. In decreasing modulus these are
\begin{equation}\label{eq:elkies-scaled-spectrum}
 \kappa_\ell=\frac{2(-1)^{\ell-1}}{(2\ell-1)\pi},\qquad \ell\ge1.
\end{equation}
One can check the correspondence directly: integration of
$\cos((2\ell-1)\pi y/2)$ gives $\kappa_\ell
\cos((2\ell-1)\pi x/2)$. The finite formula has exactly this limit:
for each fixed $\ell$, \eqref{eq:spectral-data} gives
$\lambda_{M,\ell}/M\to\kappa_\ell$ as $M\to\infty$.
The boundary conditions also agree. Differentiating
$\mathcal Kf=\kappa f$ twice gives $f''=-f/\kappa^2$, with
$f'(0)=0$ and $f(1)=0$, the continuous counterparts of the two discrete
conditions above. This comparison explains the half-integer angles and
alternating eigenvalue signs. We use the exact finite spectrum and weights,
not an operator-limit assertion, in all subsequent bounds.
\end{remark}

\subsection{Uniform isolation of the leading mode}

Write $N=\Dg+1$ and abbreviate
\begin{equation}\label{eq:leading-mode}
 \vartheta_M=\frac{\pi}{4M+2},\qquad
 \lambda_M=\frac1{2\sin\vartheta_M},\qquad
 \alpha_M=\frac4{2M+1}\cos^2\vartheta_M.
\end{equation}
The first mode is positive and strictly largest in modulus. We need an
error bound independent of $M$, because the generating function sums over
all alphabet sizes.

\begin{lemma}[Uniform spectral remainder]\label{lem:spectral-remainder}
For $M\ge1$ and $N=\Dg+1\ge4$,
\begin{equation}\label{eq:spectral-remainder}
 \Omega_{\Dg}(M)=\alpha_M\lambda_M^N(1+\varepsilon_{M,N}),
 \qquad |\varepsilon_{M,N}|\le e_N\Def2(2/5)^N.
\end{equation}
\end{lemma}

\begin{proof}
The error is zero for $M=1$. For $M\ge2$, divide
\eqref{eq:spectral-count} by its first summand. Since cosine decreases and
sine increases on $(0,\pi/2)$,
\[
 \frac{\alpha_{M,\ell}}{\alpha_M}\le1,\qquad
 \frac{|\lambda_{M,\ell}|}{\lambda_M}
 =\frac{\sin\vartheta_M}{\sin((2\ell-1)\vartheta_M)}.
\]
For $\ell=2$, the triple-angle identity and $\vartheta_M\le\pi/10$ give $\frac{\sin\vartheta_M}{\sin3\vartheta_M} =\frac1{3-4\sin^2\vartheta_M}<\frac5{13}<\frac25.$ Indeed, $\sin\vartheta_M\le\vartheta_M$ and $\pi^2<10$ imply $4\sin^2\vartheta_M<2/5$. For $\ell\ge3$, use
$\sin u\ge2u/\pi$ on $[0,\pi/2]$ and $\sin u\le u$ to obtain $|\lambda_{M,\ell}|/\lambda_M \le\pi/[2(2\ell-1)]$. Therefore $|\varepsilon_{M,N}|\le(2/5)^N+(\pi/10)^N\left(1+\frac5{2(N-1)}\right).$ Here the sum over odd denominators was extended to infinity and bounded by
its first term plus the integral of the decreasing function $(2x-1)^{-N}$ from $3$ to infinity.
More explicitly,
\[
 \sum_{\ell=3}^{\infty}(2\ell-1)^{-N}
 \le 5^{-N}+\int_3^\infty(2x-1)^{-N}\dif x
 =5^{-N}\left(1+\frac5{2(N-1)}\right).
\]
The bound therefore remains valid when $M=2$, in which case this part
of the finite spectral sum is empty, as well as when $M$ is arbitrarily
large. It is the summability in the mode index, rather than the number of
modes, that supplies uniformity. Relative to $(2/5)^N$, the second term
is at most
\[
 (\pi/4)^N\left(1+\frac5{2(N-1)}\right)
 <(4/5)^4\left(1+\frac56\right)=\frac{1408}{1875}<1.
\]
Both factors decrease with $N\ge4$, and $\pi/4<4/5$ follows from
$\pi<22/7$. This proves the claimed bound.
\end{proof}

The estimate controls the entire coefficient, including the boundary
weight. It is deliberately uniform rather than asymptotically sharp. The
next comparison requires the ratio of any two coefficients, for which
\eqref{eq:spectral-remainder} gives the common multiplicative factor
$(1+e_N)/(1-e_N)$.

\subsection{Curvature across alphabet sizes}

Extend the leading quantities to real $x\ge1$:
\begin{equation}\label{eq:continuous-leading}
 \vartheta(x)=\frac{\pi}{4x+2},\qquad
 \lambda(x)=\frac1{2\sin\vartheta(x)},\qquad
 \alpha(x)=\frac{2\cos^2\vartheta(x)}{x+1/2},\qquad
 \phi(x)=\log\lambda(x).
\end{equation}
The relevant curvature is that of $\phi$, not that of the polynomial
coefficients themselves. It measures how rapidly neighboring exponential
growth rates separate.

\begin{lemma}[Curvature and weight ratios]\label{lem:curvature}
The function $\phi$ is strictly increasing and strictly concave on
$[1,\infty)$. For every integer $j\ge2$,
\begin{equation}\label{eq:discrete-curvature}
 \delta_j\Def2\phi(j)-\phi(j-1)-\phi(j+1)
 \ge\frac9{10(j+1/2)^2}.
\end{equation}
The function $\alpha$ is decreasing, and for integers $q\ge1$,
\begin{equation}\label{eq:weight-ratios}
 \frac{\alpha_{j+q}}{\alpha_j}\le1,
 \qquad
 \frac{\alpha_{j-q}}{\alpha_j}\le1+\frac{2q}{3}
 \quad\text{when }j-q\ge1.
\end{equation}
\end{lemma}

\begin{proof}
Put $b=x+1/2$ and $\theta=\pi/(4b)$. Differentiation gives
\begin{equation}\label{eq:phi-derivatives}
 \phi'(x)=\frac{\theta\cot\theta}{b}>0,\qquad
 -\phi''(x)=\frac{\theta(\sin2\theta-\theta)}{b^2\sin^2\theta}.
\end{equation}
For $0<\theta\le\pi/6$, the Taylor bounds
$\sin2\theta\ge2\theta-(2\theta)^3/6$ and $\sin\theta\le\theta$ imply
\begin{equation}\label{eq:phi-second-lower}
 -\phi''(x)\ge\frac1{b^2}-\frac{\pi^2}{12b^4}>0.
\end{equation}
The numerator used in dividing by $\sin^2\theta$ is positive, since $1-4\theta^2/3\ge1-\pi^2/27>0$. Twice integrating a second derivative gives $\delta_j=\int_{-1}^1(1-|u|)(-\phi''(j+u))\dif u.$ To see the weight in this identity, write $\delta_j=\int_0^1[\phi'(j-1+t)-\phi'(j+t)]\dif t$ and integrate $-\phi''$ over $[j-1+t,j+t]$ inside the brackets. A point $j+u$
belongs to these intervals for a set of $t$ of length $1-|u|$. Thus no pointwise replacement of the curvature at $j$ is being made: its variation over the whole neighboring interval is retained.
Set $b=j+1/2$ and $v=b^{-2}\le4/25$. The two integrals in \eqref{eq:phi-second-lower} can be evaluated exactly:
\begin{align*}
 \int_{-1}^1\frac{1-|u|}{(b+u)^2}\dif u
 =-\log(1-v), \qquad 
 \int_{-1}^1\frac{1-|u|}{(b+u)^4}\dif u
 =\frac{v^2(3-v)}{3(1-v)^2}.
\end{align*}
For example, these are the second differences of $-\log b$ and $1/(6b^2)$, whose second derivatives are $b^{-2}$ and $b^{-4}$. It follows that $\delta_j\ge-\log(1-v)-\frac{\pi^2}{36}
\frac{v^2(3-v)}{(1-v)^2}.$ The function $(3-v)/(1-v)^2$ has derivative $(5-v)/(1-v)^3>0$.
Using $-\log(1-v)\ge v+v^2/2$ and $\pi^2<10$, we obtain
\begin{align*}
\delta_j\ge v+\frac12v^2-\frac5{18}\frac{1775}{441}v^2=v-\frac{2453}{3969}v^2\ge\frac{89413}{99225}v>\frac9{10}v.
\end{align*}
This proves \eqref{eq:discrete-curvature} with all constants uniform. For the weights, $(\log\alpha)'(x)=\frac{-1+2\theta\tan\theta}{b}<0,$ because $2\theta\tan\theta\le\pi/(3\sqrt3)<1$. This gives the first ratio in \eqref{eq:weight-ratios}. For a leftward shift, the cosine factor decreases, so
\[
 \frac{\alpha_{j-q}}{\alpha_j}
 \le\frac{j+1/2}{j-q+1/2}
 =1+\frac{q}{j-q+1/2}\le1+\frac{2q}{3}.
\]
\end{proof}

\subsection{Two classical theorems and the contours on which they are used}\label{sec:classical-tools}

Two classical theorems carry every root count in this paper. The first
counts zeros inside a closed contour; the second places a zero on the real
axis. We state both, fix the two kinds of contour to which the first is
applied, and then put each into the form in which later proofs quote it.

\begin{classicalthm}[Rouch\'e {\cite[Chapter~4, Section~5.2, Corollary to Theorem~18]{Ahlfors}}]\label{thm:rouche}
Let $\Gamma$ be a simple closed contour, and let $F$ and $\psi$ be
holomorphic on and inside $\Gamma$. If
\begin{equation}\label{eq:rouche-hypothesis}
 |F(t)-\psi(t)|<|\psi(t)|\qquad\forall t\text{ on }\Gamma,
\end{equation}
then neither function vanishes on $\Gamma$, and $F$ and $\psi$ have the same
number of zeros inside $\Gamma$, counted with multiplicity.
\end{classicalthm}

\begin{classicalthm}[Intermediate value theorem {\cite[Theorem~4.23]{RudinPMA}}]\label{thm:ivt}
A continuous real function on a closed interval that takes values of
opposite sign at the two endpoints has a zero in the open interval.
\end{classicalthm}

\Cref{thm:rouche} is a statement about a contour integral, and the
integral is easier to use once one has seen where it comes from. So before
applying the theorem we pause on the question that the edge argument asks
again and again: \emph{how many zeros does $h$ have inside a closed contour
$\Gamma$, when we do not know where they are?}

The residue theorem answers a neighboring question. It evaluates a contour
integral by adding the residues at the poles inside the contour, wherever
those poles happen to be. To count zeros with it, we therefore want a
function, manufactured from $h$, that has a pole at every zero of $h$, with
residue equal to the multiplicity of that zero, and no other singularity.
Integrating such a function around $\Gamma$ would return the number of zeros
inside, counted with multiplicity, without our having located a single one.

How does one find such a function? A polynomial is a product of factors,
one for each zero. A product is awkward to integrate, whereas a sum can be
integrated one term at a time, and the classical device for turning a
product into a sum is the logarithm. Let $h$ be a polynomial of degree
$n\ge1$. By the fundamental theorem of algebra,
\begin{equation}\label{eq:h-factored}
 h(t)=\operatorname{lc}(h)\prod_\rho(t-\rho)^{m_\rho},
 \qquad \operatorname{lc}(h)\ne0,
 \qquad \sum_\rho m_\rho=n,
\end{equation}
where $\rho$ runs over the distinct complex zeros of $h$, the integer
$m_\rho\ge1$ is the multiplicity of $\rho$, and $\operatorname{lc}(h)$ is the
leading coefficient. Taking a logarithm formally, and then differentiating
to remove it again, suggests
\begin{equation}\label{eq:log-heuristic}
 \log h(t)=\log\operatorname{lc}(h)+\sum_\rho m_\rho\log(t-\rho)
 \qquad\Longrightarrow\qquad
 \frac{h'(t)}{h(t)}=\sum_\rho\frac{m_\rho}{t-\rho}.
\end{equation}
The right-hand side is exactly the function we were looking for: one simple pole for each zero, weighted by its multiplicity. This is why $h'/h$ is called the \emph{logarithmic derivative} of $h$, and why it is the integrand of every root count below. But \eqref{eq:log-heuristic} is only a
guide. Each $\log(t-\rho)$ is multivalued around $\rho$, and we do not want the argument to depend on choosing a logarithm globally. It does not have to: the logarithm suggested the identity, and the product rule proves it. The product rule differentiates a product one factor at a time and leaves the other factors untouched. Since the factor $(t-\rho)^{m_\rho}$ has derivative $m_\rho(t-\rho)^{m_\rho-1}$, applying the rule to
\eqref{eq:h-factored} gives
\begin{equation}\label{eq:h-product-rule}
 h'(t)=\operatorname{lc}(h)\sum_\rho m_\rho(t-\rho)^{m_\rho-1}
       \prod_{\rho'\ne\rho}(t-\rho')^{m_{\rho'}}.
\end{equation}
For any $t$ that is not a zero of $h$, we may divide
\eqref{eq:h-product-rule} by \eqref{eq:h-factored}. In the term belonging
to $\rho$, the leading coefficient cancels, the factors with $\rho'\ne\rho$
cancel completely, and the factor belonging to $\rho$ loses all but one
power of $t-\rho$:
\begin{equation}\label{eq:log-derivative-partial-fractions}
\begin{aligned}
 \frac{h'(t)}{h(t)}
 &=\frac{\operatorname{lc}(h)\sum_\rho m_\rho(t-\rho)^{m_\rho-1}
         \prod_{\rho'\ne\rho}(t-\rho')^{m_{\rho'}}}
        {\operatorname{lc}(h)\prod_{\rho'}(t-\rho')^{m_{\rho'}}}\\[3pt]
 &=\sum_\rho m_\rho\,
   \frac{(t-\rho)^{m_\rho-1}}{(t-\rho)^{m_\rho}}\cdot
   \frac{\prod_{\rho'\ne\rho}(t-\rho')^{m_{\rho'}}}
        {\prod_{\rho'\ne\rho}(t-\rho')^{m_{\rho'}}}
  =\sum_\rho\frac{m_\rho}{t-\rho}.
\end{aligned}
\end{equation}
This is the identity predicted by \eqref{eq:log-heuristic}, obtained
without any logarithm. Now read it near a single zero $\rho$. The term
$m_\rho/(t-\rho)$ is a simple pole, and every other term of the sum is
holomorphic there. Hence
\begin{equation}\label{eq:log-derivative-residue}
 \Res_{t=\rho}\frac{h'(t)}{h(t)}=m_\rho.
\end{equation}
So $h'/h$ is the function we asked for: its poles are the zeros of $h$, and
its residues are their multiplicities. If $\Gamma$ is a simple closed
contour, positively oriented, with no zero of $h$ on it, the residue
theorem turns \eqref{eq:log-derivative-residue} into the count
\begin{equation}\label{eq:argument-principle}
 \frac1{2\pi\ii}\oint_\Gamma\frac{h'(t)}{h(t)}\dif t
 =\sum_{\rho\text{ inside }\Gamma}\Res_{t=\rho}\frac{h'(t)}{h(t)}
 =\sum_{\rho\text{ inside }\Gamma}m_\rho.
\end{equation}
The left-hand side is computed from the values of $h$ on the contour alone;
the right-hand side is the number of zeros inside it. This is the argument
principle \cite[Chapter~4, Section~5.2, Theorem~18]{Ahlfors}.\footnote{The
residue theorem, Cauchy's theorem and coefficient formula, and the identity
theorem are used in their standard forms; background is in
\cite[Section~1.10]{DLMF}. Subsequent uses will not carry individual
citations. Specific trigonometric products and partial-fraction identities
are cited where they enter.}

We shall also need \eqref{eq:argument-principle} for functions that are not
polynomials, because the generating series $F_{\Dg}$ is a power series. For
a function $h$ holomorphic on and inside $\Gamma$, with no zero on $\Gamma$,
there is no global factorization \eqref{eq:h-factored}, but none is needed:
the residue \eqref{eq:log-derivative-residue} is a local statement. Near a
zero, $h(t)=(t-\rho)^{m_\rho}u(t)$ with $u$ holomorphic and $u(\rho)\ne0$,
and the same product-rule computation, now with two factors, gives
\begin{equation}\label{eq:log-derivative-local}
 \frac{h'(t)}{h(t)}
 =\frac{m_\rho(t-\rho)^{m_\rho-1}u(t)+(t-\rho)^{m_\rho}u'(t)}
       {(t-\rho)^{m_\rho}u(t)}
 =\frac{m_\rho}{t-\rho}+\frac{u'(t)}{u(t)},
\end{equation}
where $u'/u$ is holomorphic near $\rho$. The residue is again $m_\rho$, and
\eqref{eq:argument-principle} holds as before.

Before using this integral, two features of the choice we have just made
are worth recording, because both return later. First, the logarithm was
only scaffolding. Where a branch of $\log h$ exists its derivative
is $h'/h$, but no branch is ever chosen, and none could be chosen on a
contour that surrounds zeros: the integral \eqref{eq:argument-principle}
measures, in units of $2\pi\ii$, how far $\log h$ is from returning to its
starting value. The obstruction to defining the logarithm is the very
number we are counting. Second, the integrand is not forced. The count used
only the principal parts in \eqref{eq:log-derivative-partial-fractions}, so
any function with simple poles of residue $m_\rho$ at the zeros, and no
other singularity inside $\Gamma$, has the same integral. Adding a function
holomorphic inside $\Gamma$ changes nothing; this is how the count will be
transferred from the series $F_{\Dg}$ to the polynomial $h_{\Dg}$ in the
proof of Proposition~\ref{prop:dominance}. Likewise a holomorphic weight
$\varphi$ replaces each residue by $m_\rho\varphi(\rho)$, as in
\eqref{eq:weighted-contour}. We use $h'/h$ because it is the member of this
family that is written directly in terms of $h$.

With the integrand understood, \cref{thm:rouche} can be read as a
comparison of two such integrals, one for $\psi$ and one for $F$. Under the
hypothesis
\eqref{eq:rouche-hypothesis}, every member of the
family $F_\tau=\psi+\tau(F-\psi)$, $0\le\tau\le1$, satisfies
\begin{equation}\label{eq:circle-homotopy-gap}
 |F_\tau(t)|\ge|\psi(t)|-\tau|F(t)-\psi(t)|>0
 \qquad\text{for every }t\text{ on }\Gamma.
\end{equation}
The integral \eqref{eq:argument-principle} for $F_\tau$ is therefore a
continuous integer-valued function of $\tau$, hence constant, and its values
at $\tau=0$ and $\tau=1$ are the two zero counts. This is the mechanism
behind \cref{thm:rouche}: the hypothesis keeps the integrand finite along
the whole deformation from $\psi$ to $F$.

We integrate over two kinds of contour. The first is a circle $|t|=a$,
always traversed counterclockwise. When $h$ has no zero on it, we write
\begin{equation}\label{eq:circle-count}
 \mathcal I_h(a)\Def
 \frac1{2\pi\ii}\oint_{|t|=a}\frac{h'(t)}{h(t)}\dif t
 =\#\{\text{zeros of }h\text{ in }|t|<a\},
\end{equation}
counted with multiplicity. The second is the boundary of an annulus
$a<|t|<b$. The annulus lies to the left of its outer circle traversed
counterclockwise and of its inner circle traversed clockwise, so its
positively oriented boundary consists of these two curves. If $h$ has no
zero on either circle, \eqref{eq:argument-principle} for this boundary
reads
\begin{equation}\label{eq:annular-count-difference}
 \frac{1}{2\pi\ii}\oint_{\partial\{a<|t|<b\}}\frac{h'(t)}{h(t)}\dif t
 =\mathcal I_h(b)-\mathcal I_h(a)
 =\sum_{a<|\rho|<b}m_\rho.
\end{equation}
The inner circle enters with a minus sign because of its orientation. The
second equality is also the additivity of the count over the three disjoint
pieces $|t|<a$, $|t|=a$, and $a<|t|<b$ of the disk $|t|<b$. Thus a
difference of two circle counts \emph{is} a contour integral around an
annulus. Its value, not the width of the annulus, is what certifies one
zero at a time below.

These integrals count the zeros in a region but do not say where in the
region they lie. \Cref{thm:ivt} locates them, in the following form.

\begin{lemma}[Sign changes exhaust a contour count]\label{lem:sign-count}
Let $h$ be holomorphic on a neighborhood of the closure of an open set
$U\subseteq\CC$ and real on the real axis, with exactly $n$ zeros in $U$,
counted with multiplicity. Suppose that $U\cap\RR$ contains $n$ pairwise
disjoint open intervals, at the two endpoints of each of which $h$ takes
values of opposite sign. Then each interval contains exactly one zero of
$h$, that zero is simple, and $h$ has no other zero in $U$, real or complex.
\end{lemma}

\begin{proof}
By \cref{thm:ivt}, each interval contains a zero. The intervals are disjoint,
so these are $n$ distinct zeros in $U$, each of multiplicity at least one.
The total multiplicity in $U$ is $n$. Hence every one of them is simple, no
interval contains a second zero, and $U$ contains no further zero.
\end{proof}

Reality and simplicity therefore enter at different steps. A sign change
produces at least one real zero; the contour count leaves no room for a
second zero, a multiple zero, or a nonreal one. On the edge, $U$ is an
annulus with $n=1$. In Sections~\ref{sec:finite} and~\ref{sec:joining} the
same lemma is used with larger $n$.

\subsection{Balanced circles and a reusable dominance principle}

We now isolate what the hypothesis \eqref{eq:rouche-hypothesis} requires of
a power series with positive coefficients when $\Gamma$ is a circle and
$\psi$ is one term of the series. It will clarify what must be supplied by
the spectral estimates.

\begin{lemma}[Positive-coefficient circle certification]\label{lem:circle-tool}
Let $F(t)=\sum_{m\ge0}b_mt^m$ be holomorphic on $|t|<R$, with $b_m>0$.
If $0<a<R$ and
\begin{equation}\label{eq:circle-tool-condition}
 \sum_{m\ne q}b_ma^m<b_qa^q,
\end{equation}
then $F$ has no zero on $|t|=a$, exactly $q$ zeros in $|t|<a$ counted with
multiplicity, that is, $\mathcal I_F(a)=q$, and $\sgn F(-a)=(-1)^q$.
If this condition holds on successively larger circles with successive
indices $q=1,\ldots,q_{\max}$, then each intervening annulus, and the innermost
disk, contains exactly one zero of $F$; that zero is simple and negative.
\end{lemma}

\begin{proof}
Apply \cref{thm:rouche} on the circle $|t|=a$ with $\psi(t)=b_qt^q$. For
$|t|=a$, the triangle inequality and the positivity of the coefficients give
\begin{equation}\label{eq:rouche-on-circle}
 |F(t)-\psi(t)|=\Bigl|\sum_{m\ne q}b_mt^m\Bigr|
 \le\sum_{m\ne q}b_ma^m<b_qa^q=|\psi(t)|.
\end{equation}
The two outer bounds do not depend on the point of the circle. This is how
the pointwise hypothesis \eqref{eq:rouche-hypothesis} becomes the single
numerical condition \eqref{eq:circle-tool-condition}; since the middle
bound is attained at $t=a$, nothing is lost in the reduction. Hence $F$ has
no zero on the circle and as many zeros inside it as $\psi$. Because
$\psi'/\psi=q/t$,
\[
 \mathcal I_F(a)=\mathcal I_\psi(a)
 =\frac1{2\pi\ii}\oint_{|t|=a}\frac qt\dif t=q.
\]
At $t=-a$ every term is real: $F(-a)=(-1)^qb_qa^q+\sum_{m\ne q}b_m(-a)^m$, and by
\eqref{eq:circle-tool-condition} the sum is smaller in absolute value than the first term. This proves the sign. Now let two consecutive circles have radii $a<a'$ and indices $q$ and
$q+1$. The closed contour is the boundary of the annulus between them, and \eqref{eq:annular-count-difference} gives
\[
 \frac1{2\pi\ii}\oint_{\partial\{a<|t|<a'\}}\frac{F'(t)}{F(t)}\dif t
 =\mathcal I_F(a')-\mathcal I_F(a)=(q+1)-q=1.
\]
The values $F(-a)$ and $F(-a')$ have the opposite signs $(-1)^q$ and
$(-1)^{q+1}$. Lemma~\ref{lem:sign-count}, with $n=1$ and the interval
$(-a',-a)$, shows that the annulus contains exactly one zero, which is
simple and negative. For the innermost disk, use $\mathcal I_F(a)=1$ at
$q=1$, the interval $(-a,0)$, and $F(0)=b_0>0>F(-a)$.
\end{proof}

For our series, $1\le\Omega_{\Dg}(M)\le M^{\Dg}$ ensures local uniform
convergence in $|t|<1$. For $j\ge2$, define the balanced radius
\begin{equation}\label{eq:balanced-radius}
 a_{\Dg,j}=\left(\frac{\lambda_{j-1}}{\lambda_{j+1}}\right)^{N/2},
 \qquad q_{N,j}=\exp(-N\delta_j/2).
\end{equation}
The radius is chosen so that the leading exponential comparisons with both neighbors have the same deficit. This balances the hardest terms on the two sides of the proposed dominant index.
Indeed, if $a$ is an unknown positive radius, the logarithms of the two neighboring exponential ratios are
\[
 N(\phi(j+1)-\phi(j))+\log a,\qquad
 N(\phi(j-1)-\phi(j))-\log a.
\]
Equating these expressions gives
$2\log a=N(\phi(j-1)-\phi(j+1))$, exactly
\eqref{eq:balanced-radius}. Their common value is $-N\delta_j/2$.
This calculation balances the exponential factors only. The boundary weights $\alpha_M$ and spectral remainders still have to be included;
\cref{prop:dominance} performs that separate step.

\begin{lemma}[Geometric decay away from the balanced index]\label{lem:balanced-decay}
The radii $a_{\Dg,j}$ increase strictly with $j$ and belong to $(0,1)$.
For every $M\ge1$,
\begin{equation}\label{eq:balanced-decay}
 \frac{\lambda_M^N a_{\Dg,j}^{M-1}}
      {\lambda_j^N a_{\Dg,j}^{j-1}}
 \le q_{N,j}^{|M-j|}.
\end{equation}
\end{lemma}

\begin{proof}
Let $\Delta_i=\phi(i+1)-\phi(i)$. Strict increase and concavity of $\phi$ imply that $\Delta_i>0$ and $\Delta_i$ decreases strictly. Also $\delta_j=\Delta_{j-1}-\Delta_j$ and $\log a_{\Dg,j}=-\frac N2(\Delta_{j-1}+\Delta_j).$ Thus $0<a_{\Dg,j}<1$, and the radii increase. The logarithm of the left-hand side of \eqref{eq:balanced-decay} is $N\left(\phi(M)-\phi(j)-\frac{M-j}{2}(\Delta_{j-1}+\Delta_j)\right).$ If $M>j$, every increment from $j$ to $M$ is at most $\Delta_j$, so this is at most $-N(M-j)\delta_j/2$. If $M<j$, every increment from $M$ to $j$ is at least $\Delta_{j-1}$, giving at most $-N(j-M)\delta_j/2$. Exponentiation proves the bound.
\end{proof}

\begin{proposition}[A sufficient domination bound]\label{prop:dominance}
Suppose $N\ge6$, $j\ge2$, and $q_{N,j}\le1/4$. Then the term
$\Omega_{\Dg}(j)t^{j-1}$ strictly dominates the sum of all other terms of
$F_{\Dg}$ on $|t|=a_{\Dg,j}$. Therefore
\begin{equation}\label{eq:edge-count-sign}
 \mathcal I_{\Dg}(a_{\Dg,j})=j-1,
 \qquad \sgn h_{\Dg}(-a_{\Dg,j})=(-1)^{j-1},
\end{equation}
where $\mathcal I_{\Dg}$ abbreviates $\mathcal I_{h_{\Dg}}$ in
\eqref{eq:circle-count}, and $h_{\Dg}$ has no zero on $|t|=a_{\Dg,j}$.
\end{proposition}

\begin{proof}
Apply \eqref{eq:spectral-remainder}, \eqref{eq:weight-ratios}, and
\eqref{eq:balanced-decay}, and write $q=q_{N,j}$. After division by the
proposed dominant term, the sum of the absolute values of all other terms
is at most
\begin{align}
 \frac{1+e_N}{1-e_N}
 \left(\sum_{m\ge1}q^m+
       \sum_{m=1}^{j-1}(1+2m/3)q^m\right)
 &\le\frac{1+e_N}{1-e_N}\Psi(q),\label{eq:tail-dominance}\\
 \Psi(q)&\Def\frac{2q}{1-q}+\frac{2q}{3(1-q)^2}.
 \notag
\end{align}
All coefficients in the two geometric sums are nonnegative, so $\Psi(q)\le\Psi(1/4)=26/27$. Moreover $e_N\le e_6=128/15625$, and $\frac{1+e_6}{1-e_6}\frac{26}{27}=\frac{136526}{139473}<1.$ For clarity, the two summations that produced the constants are $\sum_{m\ge1}q^m=q/(1-q)$ and
$\sum_{m\ge1}mq^m=q/(1-q)^2$. At $q=1/4$ they give $2/3+8/27=26/27$; the leftward boundary-weight correction is the $8/27$ term. The spectral allowance is $(1+e_6)/(1-e_6)=15753/15497$, so the remaining strict dominance margin is
\begin{equation}\label{eq:rouche-explicit-margin}
 1-\frac{136526}{139473}=\frac{2947}{139473}>0.
\end{equation}
In particular, this is a bound for the sum of \emph{all} competing terms,
not merely for either nearest neighbor.
Lemma~\ref{lem:circle-tool} now applies to $F_{\Dg}$ on $|t|=a_{\Dg,j}$,
with $q=j-1$ and $b_q=\Omega_{\Dg}(j)$. To transfer its conclusions, recall
from \eqref{eq:numerator} that $h_{\Dg}(t)=(1-t)^NF_{\Dg}(t)$, so that
\[
 \frac{h_{\Dg}'(t)}{h_{\Dg}(t)}=\frac{F_{\Dg}'(t)}{F_{\Dg}(t)}-\frac{N}{1-t}.
\]
The last term is holomorphic in $|t|<1$, and its integral over the circle
vanishes by Cauchy's theorem. Hence
$\mathcal I_{\Dg}(a_{\Dg,j})=\mathcal I_{F_{\Dg}}(a_{\Dg,j})=j-1$. The factor
$(1-t)^N$ is nonzero in $|t|<1$ and positive on the negative axis, so the
absence of zeros on the circle and the sign transfer as well.
\end{proof}

\subsection{The edge index and explicit intervals}

The curvature estimate makes a convenient choice of the last certified
index possible:
\begin{equation}\label{eq:edge-index}
 K_N=\left\lfloor\sqrt{\frac{8N}{25}}-\frac12\right\rfloor.
\end{equation}
For $2\le j\le K_N$, one has
$N\ge(25/8)(K_N+1/2)^2\ge(25/8)(j+1/2)^2$, and therefore
\[
 \frac{N\delta_j}{2}\ge\frac{45}{32},\qquad
 q_{N,j}\le\ee^{-45/32}<\frac14.
\]
The final comparison follows without decimal approximations from
$\sum_{q=0}^{4}(45/32)^q/q!>4$.

The constants in the index choice can now be read backwards. The curvature
bound contributes $9/10$, the balancing divides it by two, and the size
condition contributes $25/8$. Their product is $45/32$. It is enough
that the resulting exponential deficit makes $q_{N,j}<1/4$, because
\eqref{eq:rouche-explicit-margin} has already allocated room for every
weight and spectral error. The floor in \eqref{eq:edge-index} simply
selects all the integer indices for which this one sufficient criterion
holds; it is not intended to be an optimal cutoff. The condition $K_N\ge2$ is equivalent to
$N\ge20$, so $N\ge6$ in every application.

\begin{theorem}[Extreme root localization]\label{thm:edges}
Let $\Dg\ge19$, $N=\Dg+1$, and $K=K_N$. There is exactly one simple
negative zero in each of the intervals
\begin{equation}\label{eq:edge-intervals}
 (-a_{\Dg,2},0),\qquad
 (-a_{\Dg,j+1},-a_{\Dg,j})\quad(2\le j<K).
\end{equation}
These are all the $K-1$ zeros in $|t|<a_{\Dg,K}$. Their reciprocals are
$K-1$ simple negative zeros in $|t|>a_{\Dg,K}^{-1}$, and there are no
other zeros in that exterior region. The remaining $\Dg-2K$ zeros, counted
with multiplicity, lie in the open annulus
$a_{\Dg,K}<|t|<a_{\Dg,K}^{-1}$.
\end{theorem}

\begin{proof}
Write $h=h_{\Dg}$ and $\mathcal I=\mathcal I_{\Dg}$. The whole proof is a sequence of evaluations of \eqref{eq:argument-principle}, moving outward from the origin: first on the certified circles, then on the annuli between them, then on the circles that reciprocity reflects to the far side, and finally on the one annulus that is left over. We start with the circles themselves. By the choice of $K$, Proposition~\ref{prop:dominance} holds for $2\le j\le K$. Thus $h$ has no zero on $|t|=a_{\Dg,j}$, and

\begin{equation}\label{eq:edge-circle-data}
 \mathcal I(a_{\Dg,j})=j-1,\qquad
 \sgn h(-a_{\Dg,j})=(-1)^{j-1}\qquad(2\le j\le K).
\end{equation}
The radii are strictly nested in $(0,1)$ by Lemma~\ref{lem:balanced-decay}. Now take two consecutive circles. Fix $2\le j<K$. The closed contour for the $j$th interval of \eqref{eq:edge-intervals} is the boundary of the annulus $a_{\Dg,j}<|t|<a_{\Dg,j+1}$. By \eqref{eq:annular-count-difference} and \eqref{eq:edge-circle-data},
\begin{equation}\label{eq:edge-annulus-count}
 \frac1{2\pi\ii}\oint_{|t|=a_{\Dg,j+1}}\frac{h'(t)}{h(t)}\dif t
 -\frac1{2\pi\ii}\oint_{|t|=a_{\Dg,j}}\frac{h'(t)}{h(t)}\dif t
 =j-(j-1)=1.
\end{equation}
The signs of $h$ at $-a_{\Dg,j}$ and $-a_{\Dg,j+1}$ are opposite. By
Lemma~\ref{lem:sign-count} with $n=1$, the annulus contains exactly one zero,
which is simple and lies in $(-a_{\Dg,j+1},-a_{\Dg,j})$. For the innermost
disk, $\mathcal I(a_{\Dg,2})=1$ and $h(0)=1>0>h(-a_{\Dg,2})$ give one simple
zero in $(-a_{\Dg,2},0)$ in the same way. Summing the annular integrals,
\begin{equation}\label{eq:edge-telescoping}
 \mathcal I(a_{\Dg,K})
 =\mathcal I(a_{\Dg,2})+\sum_{j=2}^{K-1}
   \bigl[\mathcal I(a_{\Dg,j+1})-\mathcal I(a_{\Dg,j})\bigr]
 =1+(K-2)=K-1,
\end{equation}
so the zeros just located are all the zeros in $|t|<a_{\Dg,K}$.

Everything so far has happened inside the unit disk. To reach the far
zeros we need no new estimate, because reciprocity reflects the circles we
already have. Reciprocity \eqref{eq:h-reciprocity} reads
$h(t)=t^{\Dg-2}h(1/t)$. Differentiating and dividing by $h(t)$,
\begin{equation}\label{eq:reflected-log-derivative}
 \frac{h'(t)}{h(t)}
 =\frac{\Dg-2}{t}-\frac1{t^2}\,\frac{h'(1/t)}{h(1/t)}.
\end{equation}
Let $0<a<1$ be a radius with no zero on $|t|=a$. Then $|t|=a^{-1}$ carries
no zero either, because the zero set is invariant under $t\mapsto1/t$. As
$t$ describes $|t|=a^{-1}$ counterclockwise, $1/t$ describes $|t|=a$
clockwise, and $\dif(1/t)=-\dif t/t^2$. The two changes of sign cancel:
\[
 \frac1{2\pi\ii}\oint_{|t|=a^{-1}}
 \frac1{t^2}\,\frac{h'(1/t)}{h(1/t)}\dif t
 =\frac1{2\pi\ii}\oint_{|t|=a}\frac{h'(t)}{h(t)}\dif t
 =\mathcal I(a).
\]
Integrating \eqref{eq:reflected-log-derivative} over $|t|=a^{-1}$, on which
$(2\pi\ii)^{-1}\oint\dif t/t=1$, gives the reflection identity
\begin{equation}\label{eq:reflection-identity}
 \mathcal I(a^{-1})+\mathcal I(a)=\Dg-2.
\end{equation}
Since $h$ has $\Dg-2$ zeros in all, $\mathcal I(a)$ is therefore also the number of zeros in $|t|>a^{-1}$. At $a=a_{\Dg,K}$ this number is $K-1$. Reciprocity preserves multiplicity and maps the negative axis to itself, so these are exactly the reciprocals of the $K-1$ zeros found above.
One closed contour remains. It is the boundary of the annulus $a_{\Dg,K}<|t|<a_{\Dg,K}^{-1}$, which is nonempty because $a_{\Dg,K}<1$. By \eqref{eq:annular-count-difference},
\eqref{eq:reflection-identity}, and \eqref{eq:edge-telescoping},
\begin{equation}\label{eq:central-annulus-count}
 \frac1{2\pi\ii}\oint_{|t|=a_{\Dg,K}^{-1}}\frac{h'(t)}{h(t)}\dif t
 -\frac1{2\pi\ii}\oint_{|t|=a_{\Dg,K}}\frac{h'(t)}{h(t)}\dif t
 =(\Dg-K-1)-(K-1)=\Dg-2K.
\end{equation}
This is the final assertion.
\end{proof}

The contour integrals \eqref{eq:edge-telescoping},
\eqref{eq:central-annulus-count}, and \eqref{eq:reflection-identity}
partition the $\Dg-2$ zeros into three blocks, of sizes $K-1$, $\Dg-2K$,
and $K-1$. In the two outer blocks every annulus has count one, and one
sign change per annulus places its zero on the negative axis. In the
central annulus the count is $\Dg-2K$, all but $O(\sqrt N)$ of the zeros,
and the spectral argument supplies no sign inside it. Lemma~\ref{lem:sign-count}
says exactly what is missing: $\Dg-2K$ disjoint negative intervals inside
the central annulus, each with a sign change at its endpoints. The next
representation is designed to supply those signs.

% ===== 05_geometry =====
\section{Logarithmic geometry and the central sine mechanism}\label{sec:geometry}

The edge proof has located the part of the zero set for which one
alphabet-size term can control the others in modulus. The remaining question
is different: how can a sum of comparable terms acquire a predictable real
sign? Contour integrals suggest a change of representation. Instead of
estimating the original alternating power series directly, we express its
values through a kernel with explicitly known poles. The arguments of the
vectors from those poles to the evaluation point then supply the oscillation
that an absolute-value estimate would discard.

This section carries that transition all the way to an exact three-loss
identity. The bounds themselves are deferred until Section~\ref{sec:central}.
The aim is that, by the time those bounds begin, every quantity has a
geometric or algebraic role, and the reason for estimating it is already
visible.

\subsection{From annular root counts to a logarithmic strip}\label{sec:logarithm}

The root counts of Theorem~\ref{thm:edges} were values of the contour
integral \eqref{eq:argument-principle} over circles and over the boundaries
of annuli. The same residue computation gives a weighted form: if $\varphi$
is holomorphic on and inside a simple positively oriented contour $\Gamma$
on which $h$ has no zero, then
\begin{equation}\label{eq:weighted-contour}
 \frac1{2\pi\ii}\oint_\Gamma
 \varphi(t)\frac{h'(t)}{h(t)}\dif t
 =\sum_{\rho\text{ inside }\Gamma}m_\rho\varphi(\rho).
\end{equation}
The logarithmic derivative therefore provides an exact contour description of what has been counted and of what remains, namely the $\Dg-2K$ zeros of \eqref{eq:central-annulus-count} in an annulus $a<|t|<a^{-1}$. Now put
\begin{equation}\label{eq:log-map}
 t=-\ee^{-z},\qquad z=s+\ii y.
\end{equation}
Its modulus is $|t|=\ee^{-s}$. Thus the annulus
$a<|t|<a^{-1}$ becomes the vertical band
\begin{equation}\label{eq:vertical-band}
 |\Re z|< -\log a.
\end{equation}
A convenient fundamental strip is $|\Im z|<\pi$; its boundary maps to
the positive real axis, while its real axis maps bijectively to the
negative real axis. The vertical band and the horizontal fundamental
strip are different objects. Their intersection is the part of a single
logarithmic sheet corresponding to the annulus. Reciprocity becomes
$s\mapsto-s$ because $(-\ee^{-s})^{-1}=-\ee^s$.

On any contour segment on which this coordinate is one-to-one, the change
of differential is also explicit:
\[
 \frac{h'(t)}{h(t)}\dif t
 =\frac{\partial_z h(-\ee^{-z})}{h(-\ee^{-z})}\dif z,
 \qquad \frac{\dif t}{\dif z}=\ee^{-z}=-t.
\]
We do not apply this identity to an entire circle using a single-valued
logarithm across its branch cut. A circle lifts to a vertical segment of
length $2\pi$, whose endpoints differ by a period. The real negative
intervals, in contrast, lift without any branch choice or ambiguity.

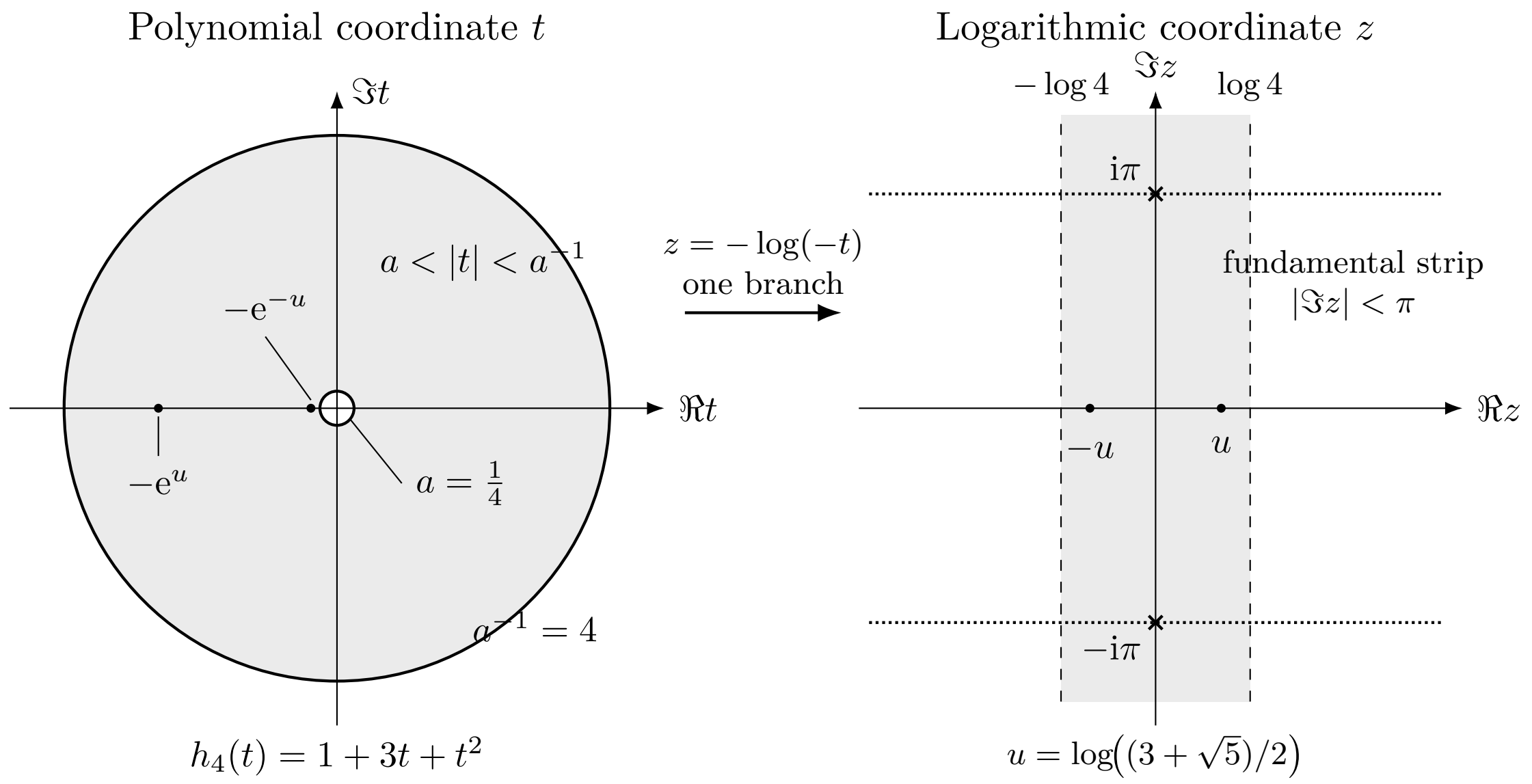
\begin{figure}[tbp]
\centering
\begin{tikzpicture}[x=.66cm,y=.66cm,>=Latex,font=\small]
\begin{scope}
 \node[font=\normalsize] at (0,5.55) {Polynomial coordinate $t$};
 \fill[black!8,even odd rule] (0,0) circle (4) (0,0) circle (.25);
 \draw[thick] (0,0) circle (4);
 \draw[thick] (0,0) circle (.25);
 \draw[->] (-4.8,0)--(4.8,0) node[right] {$\Re t$};
 \draw[->] (0,-4.65)--(0,4.65) node[right] {$\Im t$};
 \fill (-2.618034,0) circle (.065);
 \fill (-.381966,0) circle (.065);
 \draw (-2.618034,-.12)--(-2.618034,-.7)
 node[below] {$-\ee^{u}$};
 \draw (-.381966,.12)--(-1.05,1.05)
 node[above] {$-\ee^{-u}$};
 \node at (2.15,2.15) {$a<|t|<a^{-1}$};
 \draw (.2,-.17)--(.95,-1.1) node[right] {$a=\tfrac14$};
 \node at (2.9,-3.2) {$a^{-1}=4$};
 \node at (0,-5.1) {$h_4(t)=1+3t+t^2$};
\end{scope}
\draw[->,thick] (5.1,1.4)--(7.4,1.4);
\node[align=center,font=\footnotesize] at (6.25,2.15)
 {$z=-\log(-t)$\\one branch};
\begin{scope}[shift={(12,0)}]
 \node[font=\normalsize] at (0,5.55) {Logarithmic coordinate $z$};
 \fill[black!8] (-1.386294,-4.3) rectangle (1.386294,4.3);
 \draw[dashed] (-1.386294,-4.3)--(-1.386294,4.3);
 \draw[dashed] (1.386294,-4.3)--(1.386294,4.3);
 \draw[densely dotted,thick] (-4.2,3.141593)--(4.2,3.141593);
 \draw[densely dotted,thick] (-4.2,-3.141593)--(4.2,-3.141593);
 \draw[->] (-4.35,0)--(4.5,0) node[right] {$\Re z$};
 \draw[->] (0,-4.65)--(0,4.65) node[above] {$\Im z$};
 \fill (-.962424,0) circle (.065);
 \fill (.962424,0) circle (.065);
 \node[below=4pt] at (-.962424,0) {$-u$};
 \node[below=4pt] at (.962424,0) {$u$};
 \node[above,font=\footnotesize] at (-1.386294,4.3) {$-\log4$};
 \node[above,font=\footnotesize] at (1.386294,4.3) {$\log4$};
 \foreach \y in {-3.141593,3.141593}{
  \draw[thick] (-.10,\y-.10)--(.10,\y+.10);
  \draw[thick] (-.10,\y+.10)--(.10,\y-.10);
 }
 \node[above left] at (0,3.141593) {$\ii\pi$};
 \node[below left] at (0,-3.141593) {$-\ii\pi$};
 \node[align=center,font=\footnotesize] at (2.9,1.8)
 {fundamental strip\\$|\Im z|<\pi$};
 \node[align=center,font=\footnotesize] at (0,-5.1)
 {$u=\log\!\bigl((3+\sqrt5)/2\bigr)$};
\end{scope}
\end{tikzpicture}
\caption{The change of coordinates $t=-\ee^{-z}$, illustrated by
$h_4(t)=t^2+3t+1$. Its roots are
$-(3\pm\sqrt5)/2$, and their real logarithmic lifts are
$\mp u$, where $u=\log((3+\sqrt5)/2)$. Circles become vertical lines;
the horizontal lines $\Im z=\pm\pi$ bound the fundamental strip and
contain the nearest kernel poles. The illustrative annulus uses $a=1/4$;
it is not an application of the large-size edge theorem to $h_4$.}
\label{fig:log-geometry}
\end{figure}

The root-counting integral \eqref{eq:argument-principle} has poles at the
\emph{unknown zeros}. It records information but does not yet give the
missing signs. The second contour integral used below,
\[
 \frac1{2\pi\ii}\oint\frac{f(z)}{(z-s)^{\ell+1}}\dif z,
 \qquad f(z)=\frac1{2\cosh(z/2)},
\]
has a different task. Its poles are the evaluation point $s$ and the
\emph{known poles of the kernel}. Their residues reconstruct derivatives
and therefore the transformed polynomial itself. This is the new explicit
information obtained by changing representations.

This is the same coefficient-extraction mechanism as the meromorphic
expansion theorem of \cite[Theorem~IV.10, proof~(ii)]{FlajoletSedgewick}:
move a coefficient contour past known singularities and retain their
residues. Here the coefficient is taken at the variable real center $s$,
the derivative order also varies, and the poles form an infinite lattice.
The square-contour estimate below supplies the explicit remainder for
these parameters. Thus the citation provides the analytic framework, while
the required uniform and finite-size conclusions are proved here.

\subsection{The analytic destination}

For $\Dg\ge2$ and $N=\Dg+1$, define
\begin{equation}\label{eq:H-definition}
 \Hcal_{\Dg}(z)=
 \frac{\ee^{-3z/2}h_{\Dg}(-\ee^{-z})}{(1+\ee^{-z})^N},
 \qquad
 D_\ell(z)=(-\partial_z)^\ell f(z),\quad \ell\ge0.
\end{equation}
On the real axis, the factor relating $\Hcal_{\Dg}(s)$ to
$h_{\Dg}(-\ee^{-s})$ is strictly positive. Hence the two functions have
the same signs and corresponding zeros. The function $\Hcal_{\Dg}$ is
holomorphic in $|\Im z|<\pi$ and meromorphic on the plane. Since
$h_{\Dg}(1)=E_{\Dg}>0$, it has poles of exact order $N$ at every
$(2j+1)\pi\ii$. Its reflection and translation identities are
\begin{equation}\label{eq:H-symmetries}
 \Hcal_{\Dg}(-z)=(-1)^{\Dg}\Hcal_{\Dg}(z),
 \qquad
 \Hcal_{\Dg}(z+2\pi\ii)=-\Hcal_{\Dg}(z).
\end{equation}
The first follows by substituting \eqref{eq:h-reciprocity} into \eqref{eq:H-definition}; the powers of $\ee^z$ cancel because $N-(\Dg-2)=3$. The second follows from $\ee^{-3(z+2\pi\ii)/2}=-\ee^{-3z/2}$. In particular, this transformed function is antiperiodic, although its zero set repeats with period $2\pi\ii$. The exact representation we seek is
\begin{equation}\label{eq:differential-destination}
 \Hcal_{\Dg}(z)
 =-\sum_{k=0}^{\lfloor\Dg/2\rfloor}
 \frac{c_k(\Dg)E_{\Dg-2k}}{4^k(\Dg-2k)!}
 D_{\Dg-2k}(z).
\end{equation}
Here the coefficients are defined, without an asymptotic truncation, by
\begin{equation}\label{eq:c-definition}
 c_k(x)=[u^{2k}]\cos^2u\left(\frac{u}{\sin u}\right)^{x+1}.
\end{equation}
The power is defined by the analytic logarithm equal to zero at $u=0$;
equivalently it is a formal power series. For each fixed $k$, $c_k(x)$ is
a polynomial in $x$ with rational coefficients. The next two subsections
derive \eqref{eq:differential-destination}. Stating the destination first
makes the role of the centered expansion clear: powers of $M+1/2$ become
derivatives of the half-shifted geometric-series kernel $f$.

\subsection{The fixed-alphabet recurrence and a centered expansion}\label{sec:centered}

Fix $M$ and let
\begin{equation}\label{eq:alphabet-series}
 V_M(y)=\sum_{q\ge0}\Omega_q(M)y^q,
 \qquad \Omega_0(M)=1,
 \qquad V_0(y)=1.
\end{equation}
The convention at $M=0$ uses $\Omega_q(0)=0$ for $q\ge1$.
Petersen and Zhuang \cite[Theorem~3.15, Eqs.~(14)--(15)]{PZ} prove
\begin{equation}\label{eq:alphabet-recurrence}
 V_M(y)=\frac1{V_{M-1}(-y)-y}
       =\frac{y+2V_{M-1}(y)}{2-y^2-yV_{M-1}(y)}
 \quad(M\ge2).
\end{equation}
Their $F_m(y)$ is our $V_M(y)$, with their alphabet index $m$ replaced by
$M$ and their size index by $q$. Their coefficients count the same order
maps by \eqref{eq:gap-order-map}. At $M=1$, both formulas extend directly:
substituting $V_0=1$ gives $1/(1-y)$, which is $V_1$ since there is one
word of every length. These checks supply the normalization and the base
case needed to use the recurrence here.

The second fraction in \eqref{eq:alphabet-recurrence} is a fractional
linear transformation. A sine substitution turns it into the tangent
addition formula. Work near $y=0$ and put
\[
 y=2\sin u,\qquad u=\arcsin(y/2),\qquad
 W_M(u)=\frac{V_M(2\sin u)+\sin u}{\cos u}.
\]
Write $a=\sin u$ and $b=\cos u$. If $V_{M-1}=bW-a$, then
\[
 \frac{2a+2V_{M-1}}{2-4a^2-2aV_{M-1}}
 =\frac{2bW}{2b(b-aW)}=\frac{W}{b-aW},
\]
and consequently
\begin{equation}\label{eq:tangent-conjugacy}
 W_M=\frac{bW_{M-1}+a}{b-aW_{M-1}}
     =\frac{W_{M-1}+\tan u}{1-W_{M-1}\tan u}.
\end{equation}
All denominators are nonzero at the origin. Since
$W_0=(1+\sin u)/\cos u=\tan(\pi/4+u/2)$, the tangent addition formula
\cite[Eq.~4.21.4]{DLMF} gives
\[
 W_M=\tan\left(\frac\pi4+\left(M+\frac12\right)u\right)
     =\sec((2M+1)u)+\tan((2M+1)u).
\]
The last identity follows by applying the same addition formula and the
double-angle identities, or by multiplying
$(1+\sin v)/\cos v$ by $\cos(v/2)-\sin(v/2)$.
We obtain the exact local identity
\begin{equation}\label{eq:alphabet-trigonometric}
 V_M(y)=\cos u\{\sec((2M+1)u)+\tan((2M+1)u)\}-\sin u,
 \qquad y=2\sin u.
\end{equation}
The occurrence of $2M+1$ is consistent with the center $M=-1/2$ already
visible in \eqref{eq:order-reciprocity}.

\begin{remark}[Where the sine substitution comes from]\label{rem:fixed-points}
The substitution above is found by comparing fixed points; the derivation
does not depend on this remark. Write the second fraction of
\eqref{eq:alphabet-recurrence} as
$\Phi_y(V)=(y+2V)/(2-y^2-yV)$. For real $y\in(-2,2)$ with $y\ne0$ (at $y=0$
the map is the identity), its fixed points solve $V^2+yV+1=0$, that is,
$V=-\sin u\pm\ii\cos u$ when $y=2\sin u$. The tangent addition map
$R_\gamma(\tau)=(\tau+\tan\gamma)/(1-\tau\tan\gamma)$ fixes $\tau=\pm\ii$ for every $\gamma$,
and the affine map $\tau\mapsto\cos u\,\tau-\sin u$ carries $\pm\ii$ to these
two fixed points. So if $\Phi_y$ is to become a tangent addition,
the natural unknown is $\tau=(V+\sin u)/\cos u$, which is exactly $W$ above.
The fixed points cannot reveal the angle $\gamma$; that is supplied by the
computation leading to \eqref{eq:tangent-conjugacy}.
\end{remark}

\begin{proposition}[Exact centered order-polynomial expansion]\label{prop:centered}
For $\Dg\ge2$, with $c_k$ defined by \eqref{eq:c-definition},
\begin{equation}\label{eq:centered-order}
 \Omega_{\Dg}(M)=
 \sum_{k=0}^{\lfloor\Dg/2\rfloor}
 \frac{c_k(\Dg)E_{\Dg-2k}}{4^k(\Dg-2k)!}
 \left(M+\frac12\right)^{\Dg-2k}.
\end{equation}
This is a polynomial identity in $M$, not a large-$M$ approximation.
The first three coefficients are
\begin{equation}\label{eq:first-c}
 c_0(x)=1,\qquad c_1(x)=\frac{x-5}{6},\qquad
 c_2(x)=\frac{5x^2-48x+67}{360}.
\end{equation}
\end{proposition}

\begin{proof}
The term $-\sin u=-y/2$ in \eqref{eq:alphabet-trigonometric} contributes
nothing at degrees $\Dg\ge2$. To track every factor in the coefficient
change, abbreviate
$A_M(u)=\sec((2M+1)u)+\tan((2M+1)u)$ and choose a sufficiently
small positive circle around $y=0$. The local inverse $u=\arcsin(y/2)$
takes it to a simple closed contour around $u=0$. Since
$\dif y=2\cos u\dif u$, the coefficient integral becomes
\begin{align}
 \Omega_{\Dg}(M)
 &=\frac{1}{2\pi\ii}\oint
       \frac{\cos u\,A_M(u)}{(2\sin u)^{\Dg+1}}
                      2\cos u\dif u\notag\\
 &=\frac{2^{-\Dg}}{2\pi\ii}\oint
       \frac{\cos^2u\,(u/\sin u)^{\Dg+1}A_M(u)}{u^{\Dg+1}}
                      \dif u\notag\\
 &=2^{-\Dg}[u^{\Dg}]
       \cos^2u\left(\frac{u}{\sin u}\right)^{\Dg+1}A_M(u).
 \label{eq:centered-cauchy-change}
\end{align}
The second cosine is therefore the Jacobian, not another assumption about the recurrence. The factor $u/\sin u$ has value one at zero and is analytic there, so the last line is an ordinary Taylor coefficient. Shrinking the contour if necessary avoids every pole of $A_M$; no global inverse for sine is required. The even series in front contributes $c_k(\Dg)u^{2k}$. By \eqref{eq:euler-normalization}, the remaining coefficient at degree $\Dg-2k$ is $E_{\Dg-2k}(2M+1)^{\Dg-2k}/(\Dg-2k)!$. Since $2^{-\Dg}(2M+1)^{\Dg-2k}=4^{-k}(M+1/2)^{\Dg-2k}$, the result is
\eqref{eq:centered-order} for every integer $M\ge0$. Both sides are polynomials of degree at most $\Dg$, by \eqref{eq:order-binomial}, so the identity holds for all complex $M$. For the low coefficients, the Taylor series of sine and cosine \cite[Eqs.~4.19.1--4.19.2]{DLMF} give
\[
 \log\frac{u}{\sin u}=\frac{u^2}{6}+\frac{u^4}{180}+O(u^6),
 \qquad \cos^2u=1-u^2+\frac{u^4}{3}+O(u^6).
\]
Expanding the exponential with multiplier $x+1$, the coefficient of $u^2$
is $(x+1)/6-1=(x-5)/6$. The coefficient of $u^4$ is
\[
 \frac{x+1}{180}+\frac{(x+1)^2}{72}-\frac{x+1}{6}+\frac13
 =\frac{5x^2-48x+67}{360}.
\]
This proves \eqref{eq:first-c}. More generally, coefficient extraction in
the same exponential shows directly that $c_k$ has degree at most $k$
and rational coefficients.
\end{proof}

\subsection{A polynomial-to-kernel transform}\label{sec:kernel-transform}

The half-integer center now selects the kernel. For $\Re z>0$,
\begin{equation}\label{eq:sech-geometric}
 f(z)=\frac{\ee^{-z/2}}{1+\ee^{-z}}
     =\sum_{M\ge0}(-1)^M\ee^{-(M+1/2)z}.
\end{equation}
Differentiation multiplies the $M$th term by $M+1/2$. A polynomial in that
quantity can therefore be summed exactly by applying a finite differential
operator to $f$.

\begin{lemma}[Centered polynomial transform]\label{lem:transform}
Let $P(X)=\sum_{\ell=0}^q A_\ell(X+1/2)^\ell$ be a polynomial. Then
for $\Re z>0$,
\begin{equation}\label{eq:general-transform}
 \sum_{M\ge1}(-1)^{M-1}P(M)\ee^{-(M+1/2)z}
 =-\sum_{\ell=0}^q A_\ell D_\ell(z)+P(0)\ee^{-z/2}.
\end{equation}
In particular, the correction term vanishes when $P(0)=0$.
\end{lemma}

\begin{proof}
On any compact subset of $\Re z>0$, choose $\delta>0$ with
$\Re z\ge\delta$. Every derivative of order at most $q$ of the $M$th term
of \eqref{eq:sech-geometric} is bounded in modulus by
$(M+1/2)^q\ee^{-(M+1/2)\delta}$, up to a constant for the finitely many
small $M$. The majorant is summable. Hence the series and these derivatives
converge locally uniformly and can be differentiated term by term.
Multiplying the resulting identity for $D_\ell$ by $A_\ell$ and summing
proves the version of \eqref{eq:general-transform} with $M\ge0$ and a
minus sign. Removing its $M=0$ term adds precisely $P(0)\ee^{-z/2}$.
\end{proof}

Apply the lemma to $P(M)=\Omega_{\Dg}(M)$. Its value at zero is zero, and
its coefficients are those of \eqref{eq:centered-order}. On the other hand,
\eqref{eq:numerator} yields
\[
 \Hcal_{\Dg}(z)=\ee^{-3z/2}
 \sum_{M\ge1}\Omega_{\Dg}(M)(-\ee^{-z})^{M-1}
 =\sum_{M\ge1}(-1)^{M-1}\Omega_{\Dg}(M)\ee^{-(M+1/2)z}.
\]
This proves \eqref{eq:differential-destination} for $\Re z>0$. It explains
the previously chosen factor $\ee^{-3z/2}$: it changes the exponent $M-1$
into $M+1/2$. Both sides of the identity are holomorphic on the connected
strip $|\Im z|<\pi$, and agree on its right half. The identity theorem
extends it to that strip, in particular to the whole real axis. Their
meromorphic continuations then agree wherever defined.

\begin{remark}[Relation to the derivative-polynomial literature]
\label{rem:derivative-polynomials}
The derivatives just obtained belong to a classical polynomial family.
Hoffman defines $Q_\ell$ by
$\partial_u^\ell\sec u=\sec u\,Q_\ell(\tan u)$
\cite{Hoffman1995}; see also \cite[Section~1, Eq.~(1)]{Hoffman1999}.
Replacing $u$ by $\ii u$ gives
\[
 \partial_u^\ell\sech u
 =\sech u\,\widetilde Q_\ell(\tanh u),
 \qquad \widetilde Q_\ell(x)=\ii^\ell Q_\ell(\ii x).
\]
This is the hyperbolic convention of
\cite[Section~4, Eq.~(4.1) and Proposition~3]{Boyadzhiev} and
\cite[Section~2]{Ma}. The chain rule gives the real recurrence
\begin{equation}\label{eq:derivative-polynomial-recurrence}
 \widetilde Q_0(x)=1,\qquad
 \widetilde Q_{\ell+1}(x)=(1-x^2)\widetilde Q_\ell'(x)
                            -x\widetilde Q_\ell(x),
\end{equation}
so $\widetilde Q_1=-x$ and $\widetilde Q_2=2x^2-1$. In our scaling,
\begin{equation}\label{eq:derivative-polynomial-dictionary}
 D_\ell(s)=\frac{(-1)^\ell}{2^{\ell+1}}
       \sech(s/2)\,\widetilde Q_\ell(\tanh(s/2)).
\end{equation}
The factors $2^{-\ell}$ and $(-1)^\ell$ come respectively from the
half-argument and the operator $-\partial_s$.

The simple real-zero and interlacing properties of these individual
polynomials are attributed to Hetyei \cite{Hetyei} in the explicit
formulation \cite[Lemma~3.1]{Ma}. Ma's proof of Theorem~3.2 then uses
controlled signs at reference roots for a particular Eulerian combination.
Our combination \eqref{eq:differential-destination} involves all orders of
one parity, with different size-dependent coefficients. Real-rootedness
of each summand does not establish real-rootedness of their sum, and we
do not invoke an unproved closure property. Instead we obtain common
reference points from the pole geometry and bound the combined corrections.
This is why the derivative-polynomial connection motivates the kernel but
does not replace the central-sign theorem.
\end{remark}

\subsection{What a conjugate pole pair contributes}\label{sec:pole-geometry}

A simple calculation isolates the geometry before the infinite pole set is
introduced. Let $\zeta=a+\ii b$ with $b>0$, let $c=|c|\ee^{\ii\varphi}$,
and consider the real-valued principal-part pair
\[
 g(z)=\frac{c}{z-\zeta}+\frac{\overline c}{z-\overline\zeta}.
\]
For real $s$, write $s-a+\ii b=R\ee^{\ii\theta}$ with
$0<\theta<\pi$. Direct differentiation gives
\begin{equation}\label{eq:general-pole-pair}
 (-\partial_s)^\ell g(s)
 =2\ell!|c|R^{-\ell-1}
       \cos((\ell+1)\theta+\varphi).
\end{equation}
Thus the distance supplies an amplitude and differentiation multiplies the angle. The residue contributes a fixed phase offset. A sine is the special case of a purely imaginary residue; it is not a universal phase for every meromorphic function. For $f(z)=1/(2\cosh(z/2))$, the poles and residues are
\begin{equation}\label{eq:kernel-poles}
 \zeta_j^\pm=\pm(2j+1)\pi\ii,
 \qquad
 \Res(f;\zeta_j^+)=-\ii(-1)^j,
 \qquad
 \Res(f;\zeta_j^-)=\ii(-1)^j,
 \quad j\ge0.
\end{equation}
Indeed, the derivative of $2\cosh(z/2)$ is $\sinh(z/2)$, which is $\ii(-1)^j$ at $\zeta_j^+$. All poles are simple. This lattice is also consistent with the classical hyperbolic partial
fractions \cite[Eq.~4.36.5]{DLMF}: the identity $f(z)=\frac{\ii}{2}\csch\left(\frac{z+\ii\pi}{2}\right)$ shifts the integer imaginary lattice to the odd imaginary lattice. We derive the paired expansion by a contour argument below, so its summation order and the absence of an additional entire remainder are explicit. For real $s$, put
\begin{equation}\label{eq:pole-distances}
 R_j(s)=\sqrt{s^2+(2j+1)^2\pi^2},\qquad
 \theta_j(s)=\arg(s+(2j+1)\pi\ii)\in(0,\pi).
\end{equation}
The vector defining $\theta_j$ runs from the lower pole $\zeta_j^-$ to
$s$. In particular,
\begin{equation}\label{eq:angle-coordinate}
 s=\pi\cot\theta_0,
 \qquad R_0=\frac\pi{\sin\theta_0},
 \qquad R_j^2=R_0^2+4\pi^2j(j+1).
\end{equation}
The nearest-pair contribution in \eqref{eq:general-pole-pair} is
$2\ell!R_0^{-\ell-1}\sin((\ell+1)\theta_0)$.
Its zeros become equally spaced in $\theta_0$, although their real
$s$-coordinates, and hence their negative $t$-coordinates, are not equally
spaced.

\begin{figure}[tbp]
\centering
\begin{tikzpicture}[x=.7cm,y=.7cm,>=Latex,font=\small]
 \draw[->] (-1.3,0)--(6,0) node[right] {$\Re z$};
 \draw[->] (0,-4.05)--(0,4.05) node[above] {$\Im z$};
 \coordinate (L) at (0,-3);
 \coordinate (U) at (0,3);
 \coordinate (S) at (4,0);
 \draw[dashed] (U)--(S) node[midway,above right] {$R_0(s)$};
 \draw[->,thick] (L)--(S) node[midway,below right] {$s+\ii\pi$};
 \draw[dotted] (L)--(2.4,-3);
 \draw[->] (.95,-3) arc[start angle=0,end angle=36.869898,radius=.95];
 \node at (1.38,-2.57) {$\theta_0$};
 \foreach \P in {L,U}{
  \draw[thick] ($ (\P)+(-.09,-.09)$)--($ (\P)+(.09,.09)$);
  \draw[thick] ($ (\P)+(-.09,.09)$)--($ (\P)+(.09,-.09)$);
 }
 \fill (S) circle (.06);
 \node[above right] at (S) {$s>0$};
 \node[left=5pt] at (L) {$-\ii\pi$};
 \node[left=5pt] at (U) {$\ii\pi$};
 \node[align=left,anchor=west] at (6.8,2)
 {$R_0=\sqrt{s^2+\pi^2}$\\[5pt]
  $\theta_0=\arg(s+\ii\pi)$\\[5pt]
  $s=\pi\cot\theta_0$};
 \node[align=left,anchor=west] at (6.8,-1.9)
 {nearest pair in $D_\ell(s)$:\\[6pt]
  $\displaystyle 2\ell!R_0^{-\ell-1}
       \sin\bigl((\ell+1)\theta_0\bigr)$};
\end{tikzpicture}
\caption{The nearest-pole angle is based at $-\ii\pi$, not at the real
point $s$. The vector $s+\ii\pi$ has length $R_0$ and argument
$\theta_0$. Pairing the two residues converts their complex phases into
$\sin((\ell+1)\theta_0)$ after $\ell$ differentiations.
Reflection $s\mapsto-s$ replaces $\theta_0$ by $\pi-\theta_0$.}
\label{fig:pole-angle}
\end{figure}
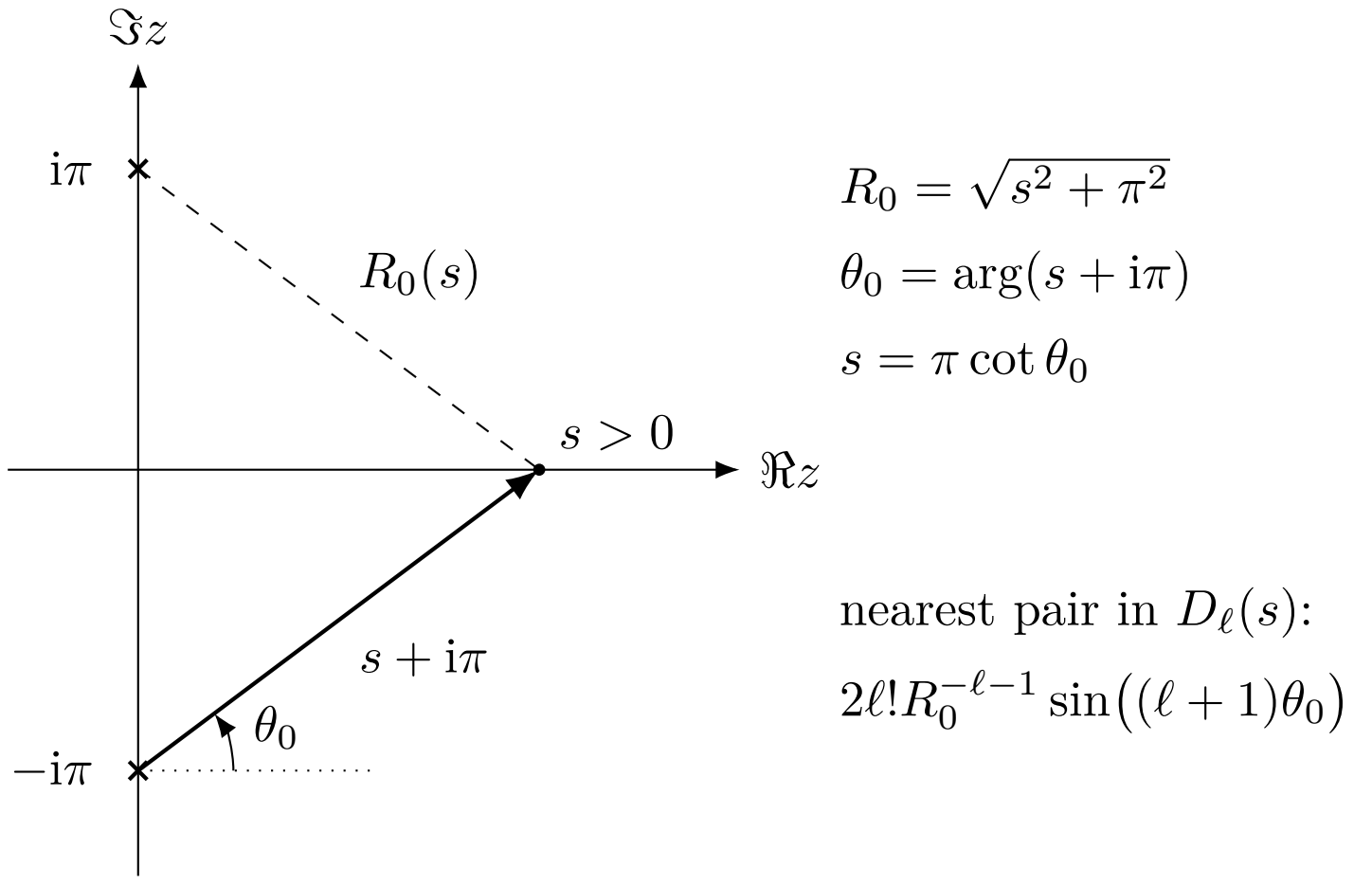

\begin{proposition}[Exact differentiated pole expansion]\label{prop:poles}
For every real $s$ and every integer $\ell\ge0$,
\begin{equation}\label{eq:pole-expansion}
 D_\ell(s)=2\ell!\sum_{j\ge0}(-1)^j
 R_j(s)^{-\ell-1}\sin((\ell+1)\theta_j(s)).
\end{equation}
The series converges absolutely for $\ell\ge1$. For $\ell=0$ it converges
in the displayed conjugate-pair order, but not absolutely.
\end{proposition}

\begin{proof}
Let $L_m=2m\pi$ and let $\mathcal Q_m$ be the open square
$|\Re z|<L_m$, $|\Im z|<L_m$. Its positively oriented boundary contains
no pole of $f$. For $L_m>|s|$, the poles inside are $s$ and
$\zeta_j^\pm$ for $0\le j<m$. The residue theorem gives
\begin{equation}\label{eq:residue-reconstruction}
 \frac1{2\pi\ii}\oint_{\partial\mathcal Q_m}
 \frac{f(z)}{(z-s)^{\ell+1}}\dif z
 =\frac{f^{(\ell)}(s)}{\ell!}
 +\sum_{j=0}^{m-1}(-1)^j
 \left(\frac{-\ii}{(\zeta_j^+-s)^{\ell+1}}
       +\frac{\ii}{(\zeta_j^--s)^{\ell+1}}\right).
\end{equation}
The first residue is the Taylor coefficient of order $\ell$ at $s$: division by $(z-s)^{\ell+1}$ makes it the coefficient of $(z-s)^{-1}$. Explicitly, on a disk about $s$ avoiding the kernel poles,
\[
 \frac{f(z)}{(z-s)^{\ell+1}}
 =\sum_{a\ge0}\frac{f^{(a)}(s)}{a!}(z-s)^{a-\ell-1}.
\]
Only $a=\ell$ has exponent $-1$. At a kernel pole $\zeta\ne s$,
$(z-s)^{-\ell-1}$ is analytic, so its value at $\zeta$ multiplies the
residue of $f$. This gives each other summand of
\eqref{eq:residue-reconstruction}. The list is exhaustive: the horizontal
sides have imaginary parts $\pm2m\pi$, strictly between consecutive odd
multiples of $\pi$, and the vertical sides have nonzero real parts.
Exactly the $m$ conjugate pairs with $0\le j<m$ lie inside.

The boundary integral vanishes as $m\to\infty$. On a horizontal side,
$f(x\pm\ii L_m)=(-1)^m/[2\cosh(x/2)]$, so the integral of its absolute
value along that side is at most
$\int_{-\infty}^{\infty}\ee^{-|x|/2}\dif x=4$. On a vertical side,
$|f(z)|\le[2\sinh(L_m/2)]^{-1}$ by the reverse triangle inequality
applied to $\ee^{z/2}+\ee^{-z/2}$. The two vertical sides have total length
$4L_m$. Finally $|z-s|\ge L_m-|s|$ on the whole boundary. Hence
\begin{equation}\label{eq:contour-vanishing}
 \left|\frac1{2\pi\ii}\oint_{\partial\mathcal Q_m}
 \frac{f(z)}{(z-s)^{\ell+1}}\dif z\right|
 \le\frac{8+2L_m/\sinh(L_m/2)}{2\pi(L_m-|s|)^{\ell+1}}
 \longrightarrow0.
\end{equation}
The numerator stays bounded and the denominator tends to infinity, even
for $\ell=0$.

The distance bound treats the sides separately. On a vertical side,
$|z-s|\ge|\Re z-s|\ge L_m-|s|$; on a horizontal side,
$|z-s|\ge|\Im z|=L_m$. For the kernel bound on the right vertical
side,
\[
 |\ee^{z/2}+\ee^{-z/2}|
 \ge\bigl|\ee^{L_m/2}-\ee^{-L_m/2}\bigr|
 =2\sinh(L_m/2),
\]
and the same argument with the two exponentials interchanged treats the
left side. Thus the perimeter growth of the square has not been ignored:
it appears as $2L_m/\sinh(L_m/2)$, which actually tends to zero.
The horizontal contribution remains bounded independently of the side
length because the kernel decays in the real direction.

Letting $m\to\infty$ in \eqref{eq:residue-reconstruction} and multiplying
by $(-1)^\ell\ell!$ yields
\[
 D_\ell(s)=\ell!\lim_{m\to\infty}\sum_{j=0}^{m-1}(-1)^j
 \left(\frac{-\ii}{(s-\zeta_j^+)^{\ell+1}}
       +\frac{\ii}{(s-\zeta_j^-)^{\ell+1}}\right).
\]
Since $s-\zeta_j^\pm=R_j\ee^{\mp\ii\theta_j}$, the bracket equals $2R_j^{-\ell-1}\sin((\ell+1)\theta_j)$, proving \eqref{eq:pole-expansion}. The contour fixes the ordering by conjugate pairs and accounts for the entire function, not just its principal parts. In particular the finite version, before the contour is sent to infinity, is worth retaining. If $I_{m,\ell}(s)$ denotes the integral on the left of \eqref{eq:residue-reconstruction}, then
\begin{equation}\label{eq:finite-pole-remainder}
 D_\ell(s)=2\ell!\sum_{j=0}^{m-1}(-1)^jR_j^{-\ell-1}
       \sin((\ell+1)\theta_j)+(-1)^\ell\ell! I_{m,\ell}(s).
\end{equation}
The signs follow from $(\zeta-s)^{-\ell-1}=(-1)^{\ell+1}(s-\zeta)^{-\ell-1}$: the minus sign when the pole sum is moved to the other side cancels $(-1)^\ell(-1)^{\ell+1}$. Pairing the two resulting terms gives $\ii(\ee^{-\ii(\ell+1)\theta_j}-\ee^{\ii(\ell+1)\theta_j}) =2\sin((\ell+1)\theta_j)$. Moreover, for $|s|\le S<L_m$, the absolute error in \eqref{eq:finite-pole-remainder} is at most
\begin{equation}\label{eq:uniform-finite-contour-error}
 \frac{\ell!\{8+2L_m/\sinh(L_m/2)\}}
      {2\pi(L_m-S)^{\ell+1}}.
\end{equation}
Hence convergence is uniform on every compact real interval for each fixed derivative order. We establish the formula separately at each order; we have not differentiated a conditionally convergent order-zero series. The vanishing explicit remainder is also what excludes an unspecified entire-function term in a principal-part expansion. For $\ell\ge1$, the absolute values are bounded by a constant times $(2j+1)^{-\ell-1}$, which is summable. For $\ell=0$, the absolute values are $2(2j+1)\pi/[s^2+(2j+1)^2\pi^2]$. Once $(2j+1)\pi\ge|s|$, this is at least $1/[(2j+1)\pi]$, so the absolute sum diverges. Convergence in the paired order has already been established by the contour limit.
\end{proof}

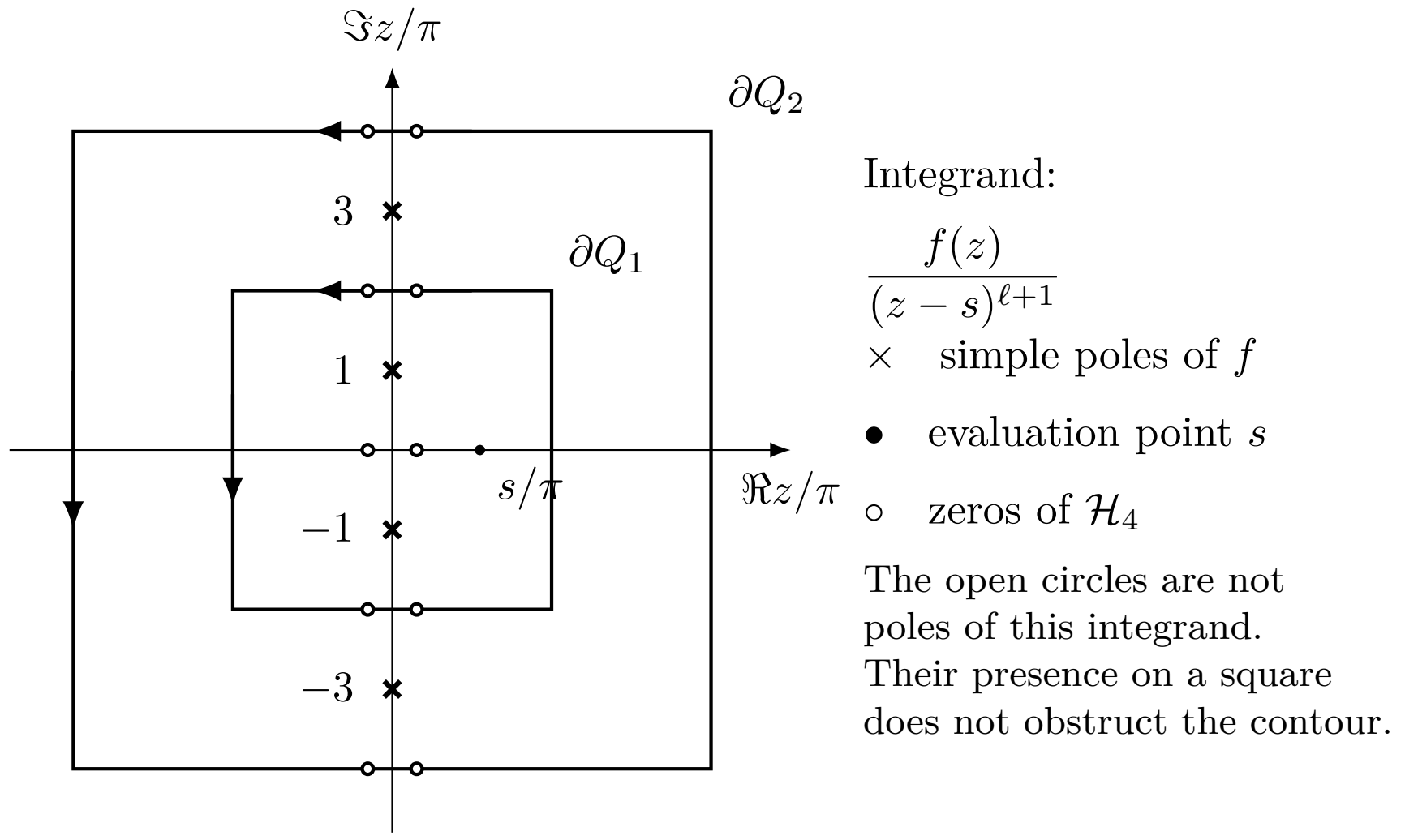
\begin{figure}[tbp]
\centering
\begin{tikzpicture}[x=.65cm,y=.65cm,>=Latex,font=\small]
 \draw[->] (-4.8,0)--(5,0)
  node[below=1pt] {$\Re z/\pi$};
 \draw[->] (0,-4.8)--(0,4.8) node[above] {$\Im z/\pi$};
 \draw[thick] (-4,-4) rectangle (4,4);
 \draw[thick] (-2,-2) rectangle (2,2);
 \draw[->,thick] (1,4)--(-1,4);
 \draw[->,thick] (1,2)--(-1,2);
 \draw[->,thick] (-4,1)--(-4,-1);
 \draw[->,thick] (-2,.7)--(-2,-.7);
 \node[above right] at (4,4) {$\partial Q_2$};
 \node[above right] at (2,2) {$\partial Q_1$};
 \foreach \y in {-3,-1,1,3}{
  \draw[very thick] (-.10,\y-.10)--(.10,\y+.10);
  \draw[very thick] (-.10,\y+.10)--(.10,\y-.10);
  \node[left=5pt] at (0,\y) {$\y$};
 }
 \foreach \y in {-4,-2,0,2,4}{
  \draw[fill=white,thick] (-.306349,\y) circle (.07);
  \draw[fill=white,thick] (.306349,\y) circle (.07);
 }
 \fill (1.1,0) circle (.065) node[below right] {$s/\pi$};
 \node[anchor=west,align=left] at (5.7,2.6)
 {Integrand:\\[6pt]$\displaystyle\frac{f(z)}{(z-s)^{\ell+1}}$};
 \node[anchor=west,align=left] at (5.7,.2)
 {$\times$\quad simple poles of $f$\\[6pt]
  $\bullet$\quad evaluation point $s$\\[6pt]
  $\circ$\quad zeros of $\Hcal_4$};
 \node[anchor=west,align=left,font=\footnotesize] at (5.7,-2.5)
 {The open circles are not\\poles of this integrand.\\Their presence on a square\\does not obstruct the contour.};
\end{tikzpicture}
\caption{Kernel contours and lifted polynomial zeros, again using \(h_4\) as an exact geometric example. The squares are \(\mathcal Q_1\) and \(\mathcal Q_2\) in the proof of Proposition~\ref{prop:poles}. Crosses are simple poles of \(f\) and order-five poles of \(\Hcal_4\). Open circles are zeros of \(\Hcal_4\), repeated by \(2\pi\ii\). Some circles lie on the square boundaries; this causes no obstruction because the integrand is \(f(z)/(z-s)^{\ell+1}\), not \(\Hcal_4'/\Hcal_4\). The real evaluation point \(s\) is chosen strictly inside the square used for that integral. The nearest-pole geometry of Figure~\ref{fig:pole-angle} may be superimposed on the present diagram at this same real evaluation point \(s\); its angle \(\theta_0\) and distance \(R_0\) are measured from \(s\) to the nearest pole pair shown here.}
\label{fig:kernel-contours}
\end{figure}

\subsection{Normalization, extrema, and the endpoint discrepancy}

Normalize each derivative to the amplitude of its nearest pole pair:
\begin{equation}\label{eq:normalized-D}
 \widehat D_\ell(s)=\frac{R_0(s)^{\ell+1}}{2\ell!}D_\ell(s).
\end{equation}
For $\ell\ge1$, Proposition~\ref{prop:poles} becomes
\begin{equation}\label{eq:normalized-poles}
 \widehat D_\ell(s)=\sin((\ell+1)\theta_0)
 +\sum_{j\ge1}(-1)^j\left(\frac{R_0}{R_j}\right)^{\ell+1}
                      \sin((\ell+1)\theta_j).
\end{equation}
The ratios are strictly below one, but approach one as $|s|\to\infty$.
Thus nearest-pair dominance becomes useful in a central window, not
uniformly along the entire real axis. For $\ell=0$, we retain the closed
form of $f$ rather than take absolute values in a conditionally convergent
pole series. Normalize the full transformed numerator in the same way:
\begin{equation}\label{eq:N-definition}
 \Ncal_{\Dg}(s)=-\frac{R_0(s)^N}{2E_{\Dg}}\Hcal_{\Dg}(s),
 \qquad
 w_{\Dg,k}(s)=c_k(\Dg)\frac{E_{\Dg-2k}}{E_{\Dg}}
                         \left(\frac{R_0(s)^2}{4}\right)^k.
\end{equation}
Substitution of \eqref{eq:differential-destination} gives the exact finite
sum
\begin{equation}\label{eq:N-sum}
 \Ncal_{\Dg}(s)=\sum_{k=0}^{\lfloor\Dg/2\rfloor}
             w_{\Dg,k}(s)\widehat D_{\Dg-2k}(s),
 \qquad w_{\Dg,0}=1.
\end{equation}

To verify the cancellation of normalizations, put $\ell=\Dg-2k$ in
one summand of \eqref{eq:differential-destination}. After multiplication
by $-R_0^N/(2E_{\Dg})$ and substitution of
$D_\ell=2\ell!R_0^{-\ell-1}\widehat D_\ell$, it becomes
\[
 \frac{R_0^N}{2E_{\Dg}}
 \frac{c_k(\Dg)E_\ell}{4^k\ell!}
 \frac{2\ell!}{R_0^{\ell+1}}\widehat D_\ell
 =c_k(\Dg)\frac{E_\ell}{E_{\Dg}}
     \frac{R_0^{2k}}{4^k}\widehat D_\ell,
\]
because $N-\ell-1=2k$. Thus the factorial in the pole residue and the
factorial in the centered expansion cancel exactly. The factorial ratio
that later makes the corrections small comes from the Euler-number ratio,
not from an additional normalization of $D_\ell$.
All normalizing factors are explicit. In particular,
\begin{equation}\label{eq:N-sign}
 \sgn h_{\Dg}(-\ee^{-s})=-\sgn\Ncal_{\Dg}(s),
 \qquad
 \Ncal_{\Dg}(-s)=(-1)^{\Dg}\Ncal_{\Dg}(s).
\end{equation}
The first term in \eqref{eq:N-sum} has principal sine
$\sin(N\theta_0)$. Its extrema occur at
\begin{equation}\label{eq:extrema}
 \theta_p=\frac{(p+1/2)\pi}{N},\qquad
 s_p=\pi\cot\theta_p,\qquad
 \sigma_p=(-1)^p,\qquad 0\le p\le\Dg.
\end{equation}
At these points the principal term has absolute value one, giving the largest available sign margin. One cannot certify its sign at every extremum. That would give $\Dg$ distinct negative roots, while $\deg h_{\Dg}=\Dg-2$. The discrepancy has an exact endpoint explanation. At fixed $\Dg$, as $s\to+\infty$,
\begin{equation}\label{eq:endpoint-asymptotic}
 \Ncal_{\Dg}(s)\sim-\frac{s^N\ee^{-3s/2}}{2E_{\Dg}}<0,
 \qquad
 \sin(N\theta_0(s))\sim\frac{N\pi}{s}>0.
\end{equation}
Here $h_{\Dg}(-\ee^{-s})\to h_{\Dg}(0)=1$, $R_0\sim s$, and
$\theta_0\sim\pi/s$. The full function approaches zero exponentially
from the opposite side to the isolated sine. Reflection gives the analogous
mismatch at the other endpoint. The omitted poles and lower derivative
terms must jointly account for this difference; they cannot all be ignored
there.

\begin{figure}[tbp]
\centering
\includegraphics[width=\linewidth]{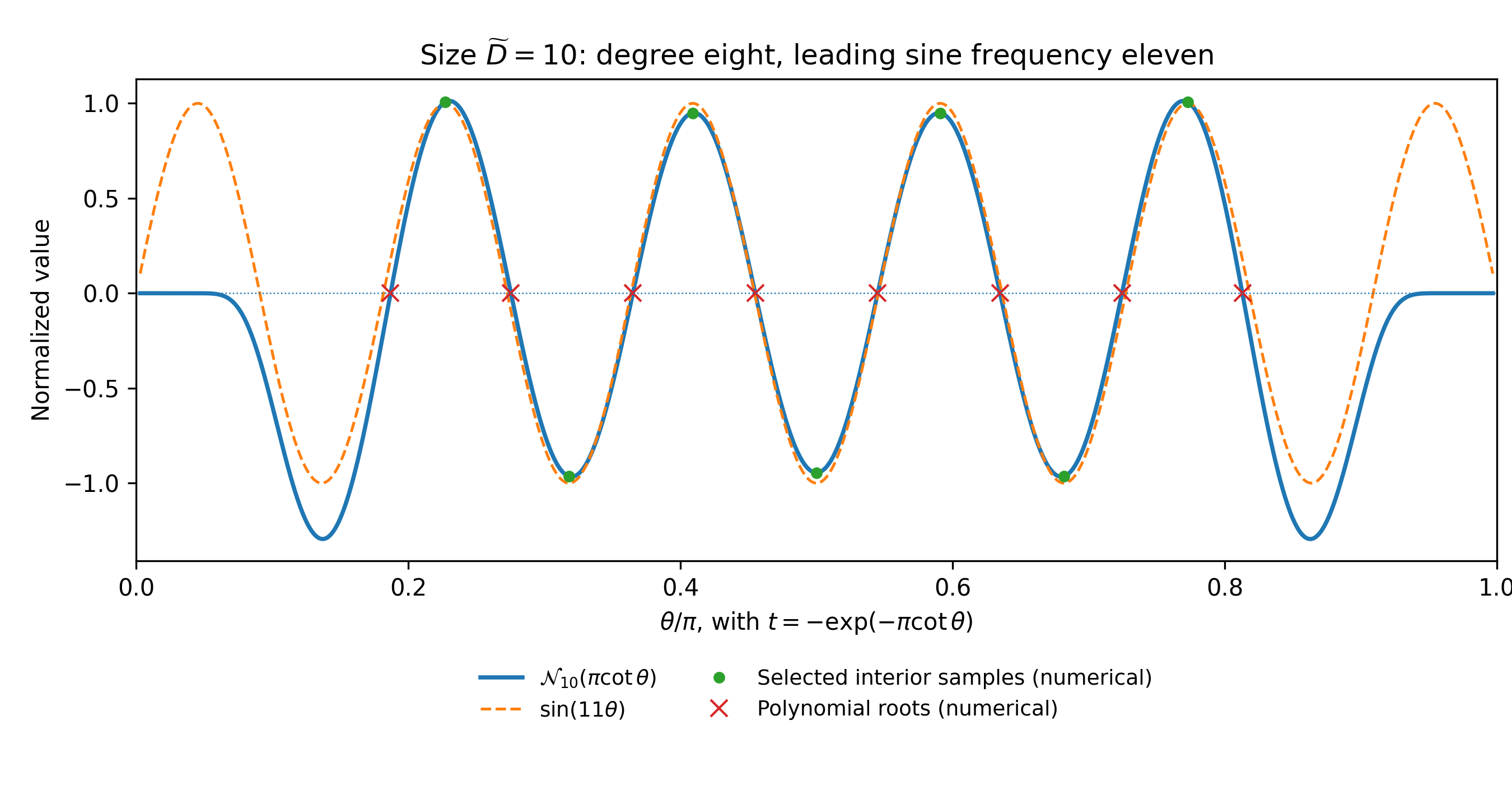}
\caption{A diagnostic at size $\Dg=10$ (degree eight). The full normalized
function $\Ncal_{10}(\pi\cot\theta)$ is compared with $\sin(11\theta)$.
Their interior oscillations are close, but the endpoint signs disagree as
in \eqref{eq:endpoint-asymptotic}; the two outermost sine zeros have no
polynomial counterparts. The plotted sample points and roots are numerical
illustrations, not sign certificates. Exact integer brackets for all eight
roots are given in Section~\ref{sec:finite}.}
\label{fig:size10}
\end{figure}

This phenomenon is not restricted to the small sizes in the figure. What
is special about those sizes is that a short reciprocal reduction certifies
all their roots without a uniform pole estimate. For larger sizes the
spectral edge count specifies which interior extrema are needed, and the
next calculation places them in a common window.

\subsection{The selected grid lies in a uniform pole-distance window}\label{sec:window}

Look again at Figure~\ref{fig:size10}. The dashed curve is the sine
suggested by the nearest pole pair, and the solid curve is the true
normalized function. Through the middle of the picture the two are almost
indistinguishable, and each sine zero sits almost on top of a polynomial
root. At
the two ends something else happens. The sine has $\Dg+1$ extrema and
$\Dg$ interior zeros, but the polynomial has only $\Dg-2$ roots, so two
sine zeros must go unmatched. The figure shows which two: at each end the
sine completes one more half-oscillation, while the true function declines
to follow it and instead flattens to zero from the opposite side, exactly
as \eqref{eq:endpoint-asymptotic} predicts. One sine zero at each end is
left with no root beneath it.

So the sine model is right in the middle and wrong at the ends, and the
proof has to say where the middle stops. This subsection answers that
question in three stages. We first ask the pole geometry why the model
fails, and find that a single ratio controls the failure. That ratio
dictates the \emph{shape} of the region in which the sine can be trusted.
We then ask the root count which extrema are actually needed, and the
answer fixes the \emph{size} of the region. Only then do we prove that the
two requirements are compatible. The constant $15/4$ and the starting size
$\Dg=11$ will both come out of that proof instead of being put into it.

\paragraph{Why the model fails at the ends.}
Figure~\ref{fig:pole-angle} gives the geometry. A real point $s$ is seen
from the pole $-\ii\pi$ at angle $\theta_0$ and distance $R_0$, and by
\eqref{eq:angle-coordinate} $R_0=\frac{\pi}{\sin\theta_0}$ and $R_j^2=R_0^2+4\pi^2j(j+1).$
In the middle of the axis $\theta_0$ is near $\pi/2$, the point $s$ is as
close to the nearest pair as it can be, and the next pair is three times
as far away. As $s$ moves out toward either end, $\theta_0$ tends to $0$ or
to $\pi$, the vector in the figure flattens, and $R_0$ grows without bound.
The second formula then shows the difficulty. The squared distances to
different pole pairs differ by \emph{fixed} amounts, so from far away all
the poles look equally distant. The nearest pair is nearest only by an
additive margin, and at the ends that margin is negligible.

Two terms of the exact expansion feel this, and it is worth seeing that
they feel it through the same quantity. The first is the contribution of
the farther poles in \eqref{eq:normalized-poles}. For derivative order
$\ell$ its size is governed by
\begin{equation}\label{eq:far-factor-heuristic}
 \left(\frac{R_0}{R_j}\right)^{\ell+1}
 =\left(1+\frac{4\pi^2j(j+1)}{R_0^2}\right)^{-(\ell+1)/2}
 \ge\exp\left(-\frac{2\pi^2j(j+1)(\ell+1)}{R_0^2}\right),
\end{equation}
where the inequality is $\log(1+y)\le y$. The orders that matter are close
to $\Dg$. If $R_0^2$ is much larger than $\Dg$, the exponent is close to
zero, the factor is close to one, and the farther pairs are as large as
the nearest one. If instead $R_0^2$ is at most a fixed multiple of $\Dg$,
the exponent stays bounded away from zero and the farther pairs are
suppressed by a fixed factor at every size. The second term is the family of weights in \eqref{eq:N-definition},
\[
 w_{\Dg,k}(s)=c_k(\Dg)\frac{E_{\Dg-2k}}{E_{\Dg}}
              \left(\frac{R_0(s)^2}{4}\right)^k,
\]
which multiply the lower derivatives and grow like $R_0^{2k}$. Their size
can be read off at leading order for a fixed $k$ as $\Dg\to\infty$, and the
computation is short enough to carry out in full. A weight is a product of
three factors, and each of them will contribute one number to the final
constant. The first factor is the coefficient $c_k(\Dg)$. In the variable $v=u^2$,
the definition \eqref{eq:c-definition} reads
\begin{align}
 c_k(\Dg)
 &=[v^k]\,\cos^2\!\sqrt v
   \left(\frac{\sqrt v}{\sin\sqrt v}\right)^{\Dg+1}\notag\\
 &=[v^k]\,\bigl(1-v+O(v^2)\bigr)
   \left(1+\frac v6+\frac{7v^2}{360}+O(v^3)\right)^{\Dg+1}\notag\\
 &=\frac1{6^k}\binom{\Dg+1}{k}+O_k(\Dg^{k-1})
  =\frac{\Dg^k}{6^kk!}\bigl(1+O_k(\Dg^{-1})\bigr).
  \label{eq:ck-leading}
\end{align}
Here and below, $O_k$ means that the implied constant depends on $k$ alone.
To see the last line, expand the power. A contribution to $v^k$ that takes
a nonconstant term from $m$ of the $\Dg+1$ equal factors can be chosen in
$\binom{\Dg+1}{m}$ ways, which is a polynomial of degree $m$ in $\Dg$. The
largest possible $m$ is $k$, and it occurs only when $k$ factors each
contribute $v/6$ and the cosine contributes its constant term. Every other
contribution uses fewer than $k$ nonconstant factors and so has degree at
most $k-1$ in $\Dg$. The result agrees with the leading coefficients $1/6$
and $5/360=1/(6^2\,2!)$ in \eqref{eq:first-c}.

The second factor is the ratio of Euler numbers, and it is worth noticing
that its size comes from a nearest pole, the same mechanism that this
section applies to the kernel. By \eqref{eq:euler-normalization},
$E_q/q!$ is the coefficient of $u^q$ in
$\sec u+\tan u=\tan(\pi/4+u/2)$. The only singularities of this function
are simple poles at $u=\pi/2+2\pi m$ for integers $m$. The one nearest to
the origin is $u=\pi/2$, and the next is $u=-3\pi/2$. Near the nearest
pole,
\begin{equation}\label{eq:euler-pole}
 \sec u+\tan u=\frac{1+\sin u}{\cos u}
 =\frac{2}{\pi/2-u}+O\!\left(\frac\pi2-u\right),
 \qquad
 \frac{2}{\pi/2-u}=\frac4\pi\sum_{q\ge0}\left(\frac2\pi\right)^qu^q.
\end{equation}
Subtracting the principal part leaves a function that is analytic in
$|u|<3\pi/2$, so by Cauchy's estimate on the circle $|u|=\pi$ its
coefficients are $O(\pi^{-q})$. Compared with $(2/\pi)^q$ this is a
relative error $O(2^{-q})$, and therefore
\begin{equation}\label{eq:euler-leading}
 E_q=q!\,\frac4\pi\left(\frac2\pi\right)^q\bigl(1+O(2^{-q})\bigr).
\end{equation}
The third factor, $(R_0^2/4)^k$, is explicit. Now multiply the three:
\begin{align}
 w_{\Dg,k}(s)
 &=c_k(\Dg)\frac{E_{\Dg-2k}}{E_{\Dg}}
   \left(\frac{R_0^2}{4}\right)^k\notag\\
 &=c_k(\Dg)\,
   \frac{(\Dg-2k)!\,\frac4\pi\left(\frac2\pi\right)^{\Dg-2k}
         \bigl(1+O(2^{-(\Dg-2k)})\bigr)}
        {\Dg!\,\frac4\pi\left(\frac2\pi\right)^{\Dg}
         \bigl(1+O(2^{-\Dg})\bigr)}
   \left(\frac{R_0^2}{4}\right)^k\notag\\
 &=c_k(\Dg)\,
   \frac{1}{\Dg^{2k}\prod_{j=0}^{2k-1}(1-j/\Dg)}
   \left(\frac{\pi^2}{4}\right)^k
   \bigl(1+O_k(2^{-\Dg})\bigr)
   \left(\frac{R_0^2}{4}\right)^k\notag\\
 &=\frac{\Dg^k}{6^kk!}\left(\frac{\pi^2}{4}\right)^k\Dg^{-2k}
   \left(\frac{R_0^2}{4}\right)^k\bigl(1+O_k(\Dg^{-1})\bigr)\notag\\
 &=\frac1{k!}\left(\frac{\pi^2R_0^2}{6\cdot4\cdot4\cdot\Dg}\right)^k
   \bigl(1+O_k(\Dg^{-1})\bigr)
  =\frac{\chi^k}{k!}\bigl(1+O_k(\Dg^{-1})\bigr),
 \label{eq:weight-heuristic}
\end{align}
where the whole dependence on the position $s$ and on the size has been
collected into the single ratio
\begin{equation}\label{eq:chi-definition}
 \chi\Def\frac{\pi^2}{96}\,\frac{R_0(s)^2}{\Dg}.
\end{equation}
The second line of \eqref{eq:weight-heuristic} inserts \eqref{eq:euler-leading} at $q=\Dg-2k$ and at
$q=\Dg$. In the third, the constants $4/\pi$ cancel, the powers of $2/\pi$
leave $(2/\pi)^{-2k}=(\pi^2/4)^k$, and the factorials are rewritten exactly
as
\[
 \frac{(\Dg-2k)!}{\Dg!}
 =\frac1{\Dg(\Dg-1)\cdots(\Dg-2k+1)}
 =\frac1{\Dg^{2k}\prod_{j=0}^{2k-1}(1-j/\Dg)},
\]
where the finite product is $1+O_k(\Dg^{-1})$ for fixed $k$. The fourth
line inserts \eqref{eq:ck-leading}. The number $96=6\cdot4\cdot4$ in the
last line keeps a record of where everything came from: the $6$ is the
$v/6$ in $c_k$, one $4$ is the $(\pi/2)^2$ of the Euler ratio, and the
other $4$ is the normalization $R_0^2/4$ in the weight. Because $R_0$
enters exactly, the error term does not depend on $s$; the approximation
is uniform along the real axis.

The corrections therefore behave like the terms of an exponential series in $\chi$, that is, in $R_0^2/\Dg$. If that ratio is bounded, they are a bounded and rapidly decaying family, and the principal sine has a chance to dominate them. If it is unbounded, no finite group of terms dominates. Both obstructions are thus controlled by the one dimensionless ratio $R_0(s)^2/\Dg$. This is what dictates the shape of the hypothesis. The natural condition under which to trust the sine model at a point $s$ is
\begin{equation}\label{eq:window-shape}
 R_0(s)^2\le\Lambda\Dg
 \qquad\Longleftrightarrow\qquad
 \sin\theta_0(s)\ge\frac{\pi}{\sqrt{\Lambda\Dg}}
\end{equation}
for a fixed constant $\Lambda$, small enough for the estimates of
Section~\ref{sec:central} to close. We call \eqref{eq:window-shape} a
\emph{pole-distance window}. In the angular coordinate of
Figure~\ref{fig:size10} it simply removes a strip of width about
$\pi/\sqrt{\Lambda\Dg}$ from each end of the axis. The extrema sit at
$\theta_p=(p+\tfrac12)\pi/N$, so the removed strips contain those with
$p+\tfrac12$ less than about $\sqrt{N/\Lambda}$, and their mirror images.
The sine model gives up roughly $\sqrt{N/\Lambda}$ extrema at each end.

\paragraph{Which extrema the root count needs.}
Giving up extrema is affordable only if they are not needed, and here the
two halves of the proof meet. The spectral argument has already certified
the zeros at each end: by Theorem~\ref{thm:edges} there are $K-1$ of them
on each side, where
\[
 K=K_N=\left\lfloor\sqrt{\frac{8N}{25}}-\frac12\right\rfloor
\]
is the edge index \eqref{eq:edge-index}. So each method loses or covers
about $\sqrt N$ points at each end, and the question becomes a comparison
of constants. The central method gives up about $\sqrt{N/\Lambda}$ extrema
and the edge method covers about $\sqrt{8N/25}$ zeros. The first number is
the smaller one precisely when
\begin{equation}\label{eq:window-leading-constant}
\frac1\Lambda<\frac8{25},\implies\Lambda>\frac{25}{8}.
\end{equation}
This is the first constant of the window, and it is inherited from the
edge cutoff. A large $\Lambda$ would satisfy it comfortably, but by
\eqref{eq:far-factor-heuristic} and \eqref{eq:weight-heuristic} every
increase in $\Lambda$ makes the farther poles and the weights larger. We
therefore want $\Lambda$ only slightly above $25/8$, and we take
$\Lambda=15/4=25/8+5/8$. The proof below shows what the extra $5/8$ pays
for.

To make the comparison exact we must say which extrema are needed, and
that is decided by signs. Let $\Dg\ge11$, $N=\Dg+1$, and $K=K_N$. Then
$K\ge1$, with $K=1$ exactly when $12\le N\le19$. By
\eqref{eq:central-annulus-count}, the central annulus contains $\Dg-2K$
zeros; when $K=1$ no edge circle is used and this is the whole degree
$\Dg-2$. To place them by Lemma~\ref{lem:sign-count} we need $\Dg-2K$
intervals with a sign change, hence $\Dg-2K-1$ interior points whose signs
alternate and whose two outermost signs are opposite to the signs at the
edge boundaries. If the principal sine prevails at $s_p$, then by
\eqref{eq:N-sign} and \eqref{eq:extrema}
\[
 \sgn h_{\Dg}(-\ee^{-s_p})=-\sigma_p=(-1)^{p+1}.
\]
At the near boundary the sign is $(-1)^{K-1}$: this is
\eqref{eq:edge-circle-data} at $j=K$ when $K\ge2$, and $h_{\Dg}(0)=1$ when
$K=1$. The extremum $p=K$ would have sign $(-1)^{K+1}=(-1)^{K-1}$, the
same as the boundary, and so cannot bound a sign change with it. The
extremum $p=K+1$ has the opposite sign. At the far boundary, reciprocity
\eqref{eq:h-reciprocity} gives the sign $(-1)^{\Dg-2}(-1)^{K-1}
=(-1)^{\Dg+K-3}$; when $K=1$ this is the sign of $h_{\Dg}(t)$ for very
negative $t$, because $h_{\Dg}$ is monic of degree $\Dg-2$. The extremum $p=\Dg-K$ has sign $(-1)^{\Dg-K+1}$, which
agrees with it, whereas $p=\Dg-K-1$ has sign $(-1)^{\Dg-K}$, which is
opposite because the two exponents differ by the odd number $2K-3$. The
extrema that the root count needs are therefore exactly
\begin{equation}\label{eq:selected-indices}
 K+1\le p\le\Dg-K-1.
\end{equation}
There are $\Dg-2K-1$ of them, and with the two edge boundaries they delimit
$\Dg-2K$ intervals, exactly the missing count. Figure~\ref{fig:size10}
shows the pattern at $\Dg=10$, where $K=1$: the seven marked extrema are
$p=2,\ldots,8$, they cut out the eight intervals that hold the eight roots,
and the unmarked extremum $p=1$ has the same sign as $h_{10}(0)=1$, so it
is not needed. That the selected points really lie between the two
boundaries is a separate matter of ordering, proved in
Section~\ref{sec:joining}. Here we show that they lie in the window.

\paragraph{The two requirements are compatible.}
By \eqref{eq:extrema}, $\theta_{\Dg-p}=\pi-\theta_p$ and
$s_{\Dg-p}=-s_p$, so the selected set is symmetric and the window
\eqref{eq:window-shape} depends only on $\sin\theta_p$. The selected
extremum closest to an end is the first one, $p=K+1$. Everything therefore
reduces to one number, the angle
$\theta_{K+1}=(K+\tfrac32)\pi/N$, and to the claim that it is not too
small.

\begin{lemma}[Selected-point window]\label{lem:window}
For every selected index in \eqref{eq:selected-indices},
\begin{equation}\label{eq:central-window}
 R_0(s_p)^2=s_p^2+\pi^2\le\frac{15}{4}\Dg.
\end{equation}
The range \eqref{eq:selected-indices} contains at least two indices, and
$\theta_{K+1}<\pi/2$.\end{lemma}

\begin{proof}
We begin with the first selected angle, and claim that it is small:
\begin{equation}\label{eq:endpoint-angle-bound}
 0<\theta_{K+1}=\frac{(K+3/2)\pi}{N}\le\frac23<\frac\pi2.
\end{equation}
The bound $\pi/2$ is what the geometry needs; the sharper bound $2/3$ is there only to make the cosecant estimate \eqref{eq:cosecant-small} below available. If $K=1$, then $N\ge12$ and $\pi<22/7$ give $\theta_2=\frac{5\pi}{2N}\le\frac{5\pi}{24}<\frac{55}{84}<\frac23.$ If $K\ge2$, the floor in \eqref{eq:edge-index} gives $K+\tfrac12\le\sqrt{8N/25}$, that is, $N\ge(25/8)(K+1/2)^2$. Hence
\[
 \theta_{K+1}\le\frac{8\pi}{25}\,
 \frac{K+3/2}{(K+1/2)^2}
 \le\frac{8\pi}{25}\cdot\frac{14}{25}
 =\frac{112\pi}{625}<\frac{2464}{4375}<\frac23.
\]
For the middle inequality put $x=K+\tfrac12\ge\tfrac52$. The function $(x+1)/x^2$ has derivative $-(x+2)/x^3<0$, so its largest value is at $x=5/2$, where it equals $(7/2)/(25/4)=14/25$. The strict inequality $\theta_{K+1}<\pi/2$ says $K+\tfrac32<N/2$, that is, $2K+2<\Dg$, or $K+1<\Dg-K-1$. Thus the two integer endpoints of \eqref{eq:selected-indices} are distinct and properly ordered. Next we check that the first selected point is the worst one, so that nothing more than \eqref{eq:endpoint-angle-bound} will be needed about the others. All selected angles lie between $\theta_{K+1}$ and $\theta_{\Dg-K-1}=\pi-\theta_{K+1}$. Since $\theta_{K+1}<\pi/2$ and the
sine is increasing on $(0,\pi/2)$ and symmetric about $\pi/2$, its smallest value over the selected angles is attained at the two ends. By \eqref{eq:angle-coordinate},
\begin{equation}\label{eq:window-reduction}
 \sin\theta_p\ge\sin\theta_{K+1},
 \qquad
 R_0(s_p)^2=\frac{\pi^2}{\sin^2\theta_p}
 \le\frac{\pi^2}{\sin^2\theta_{K+1}}
 =\pi^2\csc^2\theta_{K+1}.
\end{equation}

It remains to bound the right-hand side of \eqref{eq:window-reduction},
and there is a small mismatch to repair first. The floor in the definition of
$K$ gives information about the angle $\theta_{K+1}$, not about its sine,
so we need to pass from $\csc^2\theta$ to $\theta^{-2}$. The Laurent
expansion $\csc^2x=x^{-2}+\tfrac13+\tfrac1{15}x^2+\cdots$ shows what to
expect: the two differ by a bounded amount, starting at $1/3$ and growing
slowly. The explicit form we use is
\begin{equation}\label{eq:cosecant-small}
 \csc^2x\le x^{-2}+\frac38,\qquad 0<x\le\frac23.
\end{equation}
To prove it, write $u=x^2\le4/9$. The alternating Taylor series gives
$\sin x\ge x-x^3/6=x(1-u/6)>0$, so
$\csc^2x\le x^{-2}(1-u/6)^{-2}$, and it is enough to show
$(1-u/6)^{-2}\le1+3u/8$, that is, $(1+3u/8)(1-u/6)^2\ge1$. Direct
multiplication gives
\[
 \left(1+\frac{3u}{8}\right)\left(1-\frac u6\right)^2-1
 =\frac u{24}-\frac{7u^2}{72}+\frac{u^3}{96}
 =\frac{u(12-28u+3u^2)}{288}.
\]
The quadratic $12-28u+3u^2$ has derivative $-28+6u<0$ on $[0,4/9]$, so it
decreases there, and its value at $u=4/9$ is $4/27>0$. It is therefore
positive throughout, which proves \eqref{eq:cosecant-small}. With the sine replaced by the angle, the edge cutoff can finally enter. The index $K$ is the \emph{largest}
integer allowed by \eqref{eq:edge-index}, so $K+1$ violates the floor
condition: $K+1>\sqrt{8N/25}-\tfrac12$. Squaring,
\begin{equation}\label{eq:window-maximality}
 \left(K+\frac32\right)^2>\frac{8N}{25},
 \qquad\text{and therefore}\qquad
 \frac{\pi^2}{\theta_{K+1}^2}
 =\frac{N^2}{(K+3/2)^2}<\frac{25N}{8}.
\end{equation}
This is the exact form of the comparison \eqref{eq:window-leading-constant}: the first extremum that the edge argument does not cover is automatically close enough to the center. All the pieces are now in hand. We combine \eqref{eq:window-reduction}, \eqref{eq:cosecant-small} at $x=\theta_{K+1}$, which is permitted because $\theta_{K+1}\le2/3$ by \eqref{eq:endpoint-angle-bound}, \eqref{eq:window-maximality}, and $\pi^2<10$:
\begin{align}
 R_0(s_p)^2
 &\le\pi^2\csc^2\theta_{K+1}
 \le\frac{\pi^2}{\theta_{K+1}^2}+\frac{3\pi^2}{8}
  =\frac{N^2}{(K+3/2)^2}+\frac{3\pi^2}{8}\notag\\
 &<\frac{25N}{8}+\frac{15}{4}
 =\frac{25(\Dg+1)+30}{8}
 =\frac{25\Dg+55}{8}\le\frac{15\Dg}{4}.\label{eq:window-calculation}
\end{align}
The final comparison is $25\Dg+55\le30\Dg$, which is exactly $\Dg\ge11$.
\end{proof}

The constants can now be read backwards, and each has a single source. The
leading term $25N/8$ is the edge cutoff: it is the reciprocal of the
$8/25$ in \eqref{eq:edge-index}, which was itself chosen in
Section~\ref{sec:edges} to make the Rouch\'e margin close. The additive
$15/4$ in the second line of \eqref{eq:window-calculation} is the price of
replacing a sine by its angle, namely $3\pi^2/8$ rounded up with
$\pi^2<10$. Writing $N=\Dg+1$ turns the two into $(25\Dg+55)/8$. A window
$R_0^2\le\Lambda\Dg$ therefore contains every selected extremum as soon as
\begin{equation}\label{eq:window-threshold}
 \frac{25\Dg+55}{8}\le\Lambda\Dg \implies
 \Dg\ge\frac{55}{8\Lambda-25}.
\end{equation}
This makes the trade visible. At $\Lambda=25/8$ the condition is never
met, in agreement with \eqref{eq:window-leading-constant}. Every
$\Lambda>25/8$ works from some size on, and a smaller $\Lambda$ means a
later start. With $\Lambda=15/4$ the denominator is $5$ and the threshold
is $\Dg\ge11$. That is where the cutoff at eleven comes from: the extra
$5/8$ above $25/8$ is exactly what absorbs the constant $55/8$ from size
eleven onward. The sizes below eleven are not lost; they are the ones
certified exactly in Section~\ref{sec:finite}.

It is also worth recording what the choice buys, as orientation for
Section~\ref{sec:central}. Inside the window, the ratio $\chi$ of
\eqref{eq:chi-definition} is at most $15\pi^2/384<2/5$, so the weights
are roughly $(2/5)^k/k!$: already small at $k=3$ and decaying factorially.
The rigorous counterpart of $\chi$ is the factor $75\Dg/32$ in
\eqref{eq:weight-scale-constant}; with the same leading factors
$\Dg^k/(6^kk!)$ and $\Dg^{-2k}$ it becomes $(25/64)^k/k!$, and $25/64$ is
$15\pi^2/384$ with $\pi^2$ rounded up to $10$.
For the nearest of the farther pairs, on the boundary of the window,
\eqref{eq:far-factor-heuristic} with $\ell+1=N$ is
$(1+8\pi^2/(\Lambda\Dg))^{-N/2}$, which tends to
$\ee^{-4\pi^2/\Lambda}=\ee^{-16\pi^2/15}<3\cdot10^{-5}$ as $\Dg\to\infty$
and is larger, but still small, at the first sizes. These figures are
leading-order statements and are not used in any proof; the Euler step
is made exact, with explicit constants, in Lemma~\ref{lem:euler-ratios}. Near $\Dg=11$ the exact
coefficients differ noticeably from their leading terms; for instance
$c_1(11)=1$, whereas $11/6\approx1.83$. The rigorous versions are the
three bounds of Sections~\ref{sec:phase-bound}, \ref{sec:far-bound},
and~\ref{sec:high-bound}, each of which uses \eqref{eq:central-window} and
nothing else about the position of $s_p$.

Finally, the window should not be mistaken for a new boundary of the
unresolved annulus. It is a sufficient pole-distance condition that covers
the particular extrema required by the root count, and it is not sharp:
the chain \eqref{eq:window-calculation} fails at $\Dg=10$, yet in
Figure~\ref{fig:size10} the condition $R_0^2\le\tfrac{15}{4}\cdot10$ reads
$0.171\le\theta/\pi\le0.829$, which still contains all seven marked
extrema and excludes the four unmarked ones, where the two curves visibly
part. The cutoff at eleven is where this uniform chain of inequalities
closes, not a change in the polynomial's qualitative geometry.

\subsection{The exact three-loss decomposition}\label{sec:loss-partition}

Fix an extremum $s_p$ with $\Dg\ge11$ and
\eqref{eq:central-window}; suppress its argument in $R_j,\theta_j,w_{\Dg,k}$.
By \eqref{eq:extrema},
\begin{equation}\label{eq:phase-lock}
 \sigma_p\sin((N-2k)\theta_p)=\cos(2k\theta_p).
\end{equation}

Indeed, $\sin(N\theta_p)=\sigma_p$ and $\cos(N\theta_p)=0$.
The subtraction formula therefore reads
\[
 \sin(N\theta_p-2k\theta_p)
 =\sigma_p\cos(2k\theta_p)-0\sin(2k\theta_p).
\]
Multiplying by $\sigma_p$, whose square is one, gives \eqref{eq:phase-lock}. The reference points belong to the leading oscillation, but this identity controls the phases of every lower
same-parity derivative there. This is the useful phase relation: the lower derivatives do not contribute arbitrary phases at the principal extrema. Their nearest-pair phases are locked to explicit cosines. The first three weights satisfy $w_{\Dg,0}=1$ and $w_{\Dg,1},w_{\Dg,2}>0$. Indeed, \eqref{eq:first-c} gives $c_1>0$ and $c_2>0$ for $\Dg\ge11$; the quadratic numerator of $c_2$ is increasing there and equals $144$ at eleven. The other factors in the weights are
positive. We therefore retain $k=0,1,2$ together, rather than estimating $k=1,2$ as purely adverse errors. Define
\begin{align}
 \Lphase&=w_{\Dg,1}(1-\cos2\theta_p)
          +w_{\Dg,2}(1-\cos4\theta_p),\label{eq:phase-loss}\\
 X_{\mathrm{far}}
 &=\sum_{k=0}^2w_{\Dg,k}\sum_{j\ge1}(-1)^j
       (R_0/R_j)^{N-2k}\sigma_p\sin((N-2k)\theta_j),\label{eq:Xfar}\\
 X_{\mathrm{high}}
 &=\sigma_p\sum_{k=3}^{\lfloor\Dg/2\rfloor}
              w_{\Dg,k}\widehat D_{\Dg-2k}(s_p),\label{eq:Xhigh}\\
 \Lfar&=\sum_{k=0}^2w_{\Dg,k}\sum_{j\ge1}(R_0/R_j)^{N-2k},\label{eq:far-loss}\\
 \Lhigh&=\sum_{k=3}^{\lfloor\Dg/2\rfloor}
              |w_{\Dg,k}|\,|\widehat D_{\Dg-2k}(s_p)|.\label{eq:high-loss}
\end{align}
The farther-pole series in the retained terms are absolutely convergent:
the smallest derivative order is $\Dg-4\ge7$. The final sum is finite
and leaves any derivative of order zero intact.

\begin{proposition}[Sign accounting at an extremum]\label{prop:loss-partition}
With the definitions above,
\begin{equation}\label{eq:exact-loss-identity}
 \sigma_p\Ncal_{\Dg}(s_p)
 =1+w_{\Dg,1}+w_{\Dg,2}-\Lphase
       +X_{\mathrm{far}}+X_{\mathrm{high}},
\end{equation}
where $\Lphase,\Lfar,\Lhigh\ge0$,
$|X_{\mathrm{far}}|\le\Lfar$, and
$|X_{\mathrm{high}}|\le\Lhigh$. Consequently,
\begin{equation}\label{eq:loss-lower-bound}
 \sigma_p\Ncal_{\Dg}(s_p)
 \ge1-\Lphase-\Lfar-\Lhigh.
\end{equation}
\end{proposition}

\begin{proof}
For $k=0,1,2$, insert \eqref{eq:normalized-poles} into
\eqref{eq:N-sum} and multiply by $\sigma_p$. The nearest-pair terms sum,
by \eqref{eq:phase-lock}, to
$1+w_{\Dg,1}\cos2\theta_p+w_{\Dg,2}\cos4\theta_p$.
Subtracting and adding $w_{\Dg,1}+w_{\Dg,2}$ gives the first four terms
of \eqref{eq:exact-loss-identity}; all remaining terms are exactly
\eqref{eq:Xfar} and \eqref{eq:Xhigh}. Positivity of the retained weights
and $1-\cos u\ge0$ give nonnegativity. The stated absolute bounds follow
from the triangle inequality and $|\sin u|\le1$. Dropping the two
nonnegative reinforcement terms proves \eqref{eq:loss-lower-bound}.
\end{proof}

This is a partition, not a list of partially overlapping error estimates.
The first loss concerns only the nearest pair in the retained corrections;
the second concerns only farther pairs in those same retained terms; the
third contains every remaining correction as a complete derivative.
The adjective ``high'' refers to the correction index $k$, not to the
derivative order, which decreases as $k$ increases.

One can check the accounting by classifying a term by both of its indices.
For $k\le2$, the pair $j=0$ is assigned to the principal sine and the
retained reinforcement, with only its cosine deficit charged to
$\Lphase$; all $j\ge1$ at those same $k$ are charged to $\Lfar$.
For $k\ge3$ the whole derivative, including its nearest pair, is charged
to $\Lhigh$. Thus the three domains are disjoint and cover the finite
differential combination. In particular, taking absolute values of the
last group does not duplicate the farther-pole error in the first group.

The unit baseline in \eqref{eq:loss-lower-bound} is deliberately
conservative. Even if the retained cosines are not both positive,
$w_1\cos(2\theta_p)+w_2\cos(4\theta_p)
=w_1+w_2-\Lphase$ is an exact equality. We retain its nonnegative
$w_1+w_2$ on the identity side and may then discard it in a lower bound.
No claim that a lower derivative is globally small, or that its sign agrees
with the leading derivative away from the chosen points, is used.

The organization also identifies three different mechanisms of smallness.
Phase loss is controlled by the trigonometric deficit and the first two
exact coefficient formulas. Farther poles are controlled by their relative
distances raised to powers comparable with $\Dg$. The remaining corrections
are controlled by Euler-number ratios, a factorial quotient, and a uniform
coefficient majorant. Section~\ref{sec:central} proves
\eqref{eq:intro-budget}; then \eqref{eq:loss-lower-bound} gives
\[
 \sigma_p\Ncal_{\Dg}(s_p)>1-\frac{19}{100}-\frac{13}{1000}-\frac34
 =\frac{47}{1000}.
\]
The resulting sign transfers to $h_{\Dg}$ through
\eqref{eq:N-sign}. Before proving those uniform estimates, we dispose of
the smaller sizes by exact algebra, so that no exceptional polynomial
remains outside the final root count.

% ===== 06_finite =====
\section{Exact finite-size certificates}\label{sec:finite}

The endpoint discrepancy in \eqref{eq:endpoint-asymptotic} explains why the
principal sine should not be imposed globally. It does not explain a sudden
change at size eleven: no such change is asserted. Eleven is the lower
threshold of our uniform window calculation. Moreover, at some smaller
sizes the retained coefficients themselves have the wrong signs for the
nonnegative-loss decomposition: $c_1(\Dg)<0$ for $2\le\Dg\le4$, and
$c_2(\Dg)<0$ for $4\le\Dg\le7$. The differential identity remains exact;
it is the chosen inequality that no longer applies unchanged.

Reciprocity provides a simpler finite treatment. It pairs the negative
roots and reduces the required degree by roughly one half. For sizes up to
ten the resulting polynomials have degree at most four, and integer sign
brackets certify every root. No approximate root finder is needed.

\subsection{A reciprocal coordinate}

\begin{lemma}[Reciprocal reduction]\label{lem:reciprocal-reduction}
Let $h$ be a real monic reciprocal polynomial of degree $2m+\epsilon$,
where $\epsilon\in\{0,1\}$ and $h(0)=1$. There is a unique real monic
polynomial $\Upsilon$ of degree $m$ such that
\begin{equation}\label{eq:reciprocal-reduction}
 h(t)=(1+t)^\epsilon(-t)^m \Upsilon(x),
 \qquad x=-\frac{(1+t)^2}{t}=-t-2-t^{-1}.
\end{equation}
If $\Upsilon$ has $m$ distinct positive roots, then $h$ has $2m+\epsilon$ distinct
strictly negative roots.
\end{lemma}

\begin{proof}
A reciprocal polynomial of odd degree vanishes at $-1$, since
$h(-1)=-h(-1)$. Dividing by $1+t$ gives a monic reciprocal polynomial of
even degree and constant coefficient one. Thus it suffices to consider
$h(t)/(1+t)^\epsilon$. Multiplication by $t^{-m}$ gives a Laurent polynomial
invariant under $t\mapsto t^{-1}$. Every such Laurent polynomial is a
polynomial in $u=t+t^{-1}$: the symmetric terms $t^j+t^{-j}$ satisfy
$C_0=2$, $C_1=u$, and $C_{j+1}=uC_j-C_{j-1}$. Their leading coefficients
are one for $j\ge1$, so the resulting polynomial in $u$ is monic of degree
$m$ and unique. Replacing $u$ by $-x-2$ and multiplying by $(-1)^m$ gives
the asserted monic $\Upsilon$. If its roots are $\xi_1,\ldots,\xi_m>0$, then substitution gives
\begin{equation}\label{eq:quadratic-factors}
 h(t)=(1+t)^\epsilon\prod_{j=1}^m
       \bigl(t^2+(\xi_j+2)t+1\bigr).
\end{equation}
Each quadratic has discriminant $\xi_j(\xi_j+4)>0$, positive root product, and negative root sum. It therefore has two distinct negative roots, $\frac{-(\xi_j+2)\pm\sqrt{\xi_j(\xi_j+4)}}2.$
Quadratics with distinct $\xi_j$ cannot share a nonzero root, because subtracting their equations would give $(\xi_i-\xi_j)t=0$. Neither has root $-1$, where its value is $-\xi_j$. Thus all factors in \eqref{eq:quadratic-factors} supply distinct negative roots.
\end{proof}
This coordinate is consistent with the earlier logarithmic geometry.
For $t=-\ee^{-s}$,
\begin{equation}\label{eq:reciprocal-log-coordinate}
 x=\ee^s+\ee^{-s}-2=4\sinh^2(s/2).
\end{equation}
A positive root of $\Upsilon$ gives two opposite real values of $s$, hence a
reciprocal pair of negative values of $t$. The optional factor $1+t$
corresponds to $s=0$.

\subsection{Integer coefficient extraction}

For $\Dg\ge2$, multiplication of \eqref{eq:numerator} by $(1-t)^{\Dg+1}$
gives
\begin{equation}\label{eq:finite-extraction}
 [t^k]h_{\Dg}(t)=\sum_{j=0}^k(-1)^j\binom{\Dg+1}{j}
                      \Omega_{\Dg}(k+1-j),
 \qquad 0\le k\le\Dg-2.
\end{equation}
The gap counts can be generated entirely with integers. If $b_q^{(M)}(a)$
counts words of length $q$ ending at gap $a$, then
\begin{equation}\label{eq:integer-gap-recurrence}
 b_1^{(M)}(a)=1,\qquad
 b_{q+1}^{(M)}(a)=\sum_{u=0}^{M-1-a}b_q^{(M)}(u),\qquad
 \Omega_{\Dg}(M)=\sum_{a=0}^{M-1}b_{\Dg}^{(M)}(a).
\end{equation}
This is the same compatibility test as \eqref{eq:threshold}. The zero
boundary gaps add no further restriction to the first or last coordinate.
\begin{table}[tbp]
\centering
\begin{tabular}{@{}c rrrrr@{}}\toprule
$\Dg$ & $[t^0]h$ & $[t^1]h$ & $[t^2]h$ & $[t^3]h$ & $[t^4]h$\\\midrule
2&1&--&--&--&--\\
3&1&--&--&--&--\\
4&1&3&--&--&--\\
5&1&7&--&--&--\\
6&1&14&31&--&--\\
7&1&26&109&--&--\\
8&1&46&334&623&--\\
9&1&79&937&2951&--\\
10&1&133&2475&12331&20641\\\bottomrule
\end{tabular}
\caption{The first half of the exact numerator coefficients. Reciprocity
supplies the complementary half. A dash means that the index lies beyond
the half listed, not that the coefficient vanishes.}
\label{tab:finite-coefficients}
\end{table}
Equations~\eqref{eq:finite-extraction}--\eqref{eq:integer-gap-recurrence}
give Table~\ref{tab:finite-coefficients}. For example,
\[
 (\Omega_{10}(1),\ldots,\Omega_{10}(5))=(1,144,4004,48620,358671).
\]
The first five coefficients are obtained as follows:
\begin{align*}
 [t^0]h_{10}&=1,\\
 [t^1]h_{10}&=144-11=133,\\
 [t^2]h_{10}&=4004-11\cdot144+55=2475,\\
 [t^3]h_{10}&=48620-11\cdot4004+55\cdot144-165=12331,\\
 [t^4]h_{10}&=358671-11\cdot48620+55\cdot4004-165\cdot144+330
             =20641.
\end{align*}
These determine the entire reciprocal polynomial. This example also
specifies all operations needed to reproduce the smaller rows.

\subsection{The complete finite certificate}

For each size, put $m=\lfloor(\Dg-2)/2\rfloor$ and
$\epsilon=\Dg-2-2m$. The polynomials from
Lemma~\ref{lem:reciprocal-reduction} are
\begin{equation}\label{eq:finite-Q}
\begin{aligned}
 \Upsilon_2(x)=\Upsilon_3(x)&=1,\\
 \Upsilon_4(x)&=x-1,\\
 \Upsilon_5(x)&=x-4,\\
 \Upsilon_6(x)&=x^2-10x+5,\\
 \Upsilon_7(x)&=x^2-21x+36,\\
 \Upsilon_8(x)&=x^3-40x^2+159x-45,\\
 \Upsilon_9(x)&=x^3-72x^2+556x-528,\\
 \Upsilon_{10}(x)&=x^4-125x^3+1697x^2-3612x+665.
\end{aligned}
\end{equation}
The reductions are exact polynomial identities: a term $q_jx^j$ on the
right side of \eqref{eq:reciprocal-reduction} contributes
$q_j(-1)^{m+j}t^{m-j}(1+t)^{2j+\epsilon}$. Since $0\le j\le m$, there are
no negative powers, and ordinary expansion reproduces
Table~\ref{tab:finite-coefficients} and its reciprocal half.

\begin{table}[tbp]
\centering
\begin{tabular}{@{}c l l@{}}\toprule
$\Dg$ & Evaluation points & Exact values of $\Upsilon_{\Dg}$\\\midrule
4 & $(0,2)$ & $(-1,1)$\\
5 & $(0,5)$ & $(-4,1)$\\
6 & $(0,1,10)$ & $(5,-4,5)$\\
7 & $(0,2,21)$ & $(36,-2,36)$\\
8 & $(0,1,5,40)$ & $(-45,75,-125,6315)$\\
9 & $(0,2,8,72)$ & $(-528,304,-176,39504)$\\
10 & $(0,1,3,13,110)$ & $(665,-1374,1808,-5562,172045)$\\\bottomrule
\end{tabular}
\caption{Sign certificates for \cref{prop:finite}. At the listed integer
points the reduced polynomial $\Upsilon_{\Dg}$ of \eqref{eq:finite-Q}
alternates in sign, giving $\deg\Upsilon_{\Dg}$ disjoint positive intervals,
each containing exactly one simple root.}
\label{tab:finite-signs}
\end{table}

\begin{proposition}[Finite-size root certification]\label{prop:finite}
Theorem~\ref{thm:main} holds for $1\le\Dg\le10$.
\end{proposition}

\begin{proof}
The constant numerators $h_1,h_2$ and the polynomial $h_3=1+t$ are
immediate. For $4\le\Dg\le10$, Table~\ref{tab:finite-signs} gives $m$
disjoint open intervals in $(0,\infty)$ with opposite endpoint signs.
A polynomial of degree $m$ has $m$ zeros in the plane, counted with
multiplicity. Lemma~\ref{lem:sign-count}, with $U=\CC$ and $n=m$, therefore
gives exactly one root of $\Upsilon_{\Dg}$ in each interval, shows that it is
simple, and excludes any other root.
Lemma~\ref{lem:reciprocal-reduction} proves the desired conclusion for
$h_{\Dg}$. Finally $h_{\Dg}(0)=1$, so multiplying by $t$ adds exactly one
simple zero at zero.
\end{proof}

The argument uses the figures only for motivation. The certificate consists
of integer counts, polynomial identities, and sign evaluations. The same
operations are included in the verification file described in
Appendix~\ref{app:certificates}. We now turn to the estimates that cover
all remaining sizes at once.

% ===== 07_central =====
\section{The uniform central-sign theorem}\label{sec:central}

The decomposition in \cref{prop:loss-partition} has reduced a question about
zeros to three estimates at a prescribed set of real points. The estimates
have different purposes. The phase estimate measures how much of the first
two reinforcing corrections is lost when their oscillations are not exactly
aligned. The farther-pole estimate measures the effect of replacing three
complete derivative expansions by their nearest pole pairs. The last estimate
controls every correction not retained in those three expansions. In
particular, it must include derivatives of orders zero and one; these are not
uniformly controlled by an absolute sum of the order-zero pole expansion.
Here is the statement that the estimates will prove. Its window hypothesis
makes it useful independently of the particular index set selected in
\cref{eq:selected-indices}.

\begin{theorem}[Central signs]\label{thm:central-sign}
Let $\Dg\ge11$, $N=\Dg+1$, and let $p\in\{0,\ldots,\Dg\}$.
At the point $s_p$ of \cref{eq:extrema}, suppose that
$R_0(s_p)^2\le15\Dg/4$. Then
\begin{equation}\label{eq:central-sign-conclusion}
 (-1)^p\Ncal_{\Dg}(s_p)>\frac{47}{1000}>0,
 \qquad
 \sgn h_{\Dg}(-\ee^{-s_p})=(-1)^{p+1}.
\end{equation}
In particular, these conclusions hold for every
$K_N+1\le p\le\Dg-K_N-1$.
\end{theorem}

Throughout this section, $s=s_p$, $R_0=R_0(s)$, and
$\theta=\theta_0(s)$ satisfy these hypotheses. We occasionally write
$w_k=w_{\Dg,k}(s)$ when the size and evaluation point are fixed. The relation
$\sin\theta=\pi/R_0$ from \cref{eq:angle-coordinate} will repeatedly convert
phase errors into powers of the same distance that occurs in the coefficient
weights. This cancellation is the main reason the estimates remain uniform
as the polynomial size grows.

\subsection{Two estimates shared by the three losses}\label{sec:shared-bounds}

We first separate the combinatorial normalization from the analytic
normalization. The Euler numbers control the coefficients of the differential
combination; the normalized derivatives control the pole geometry. Neither
estimate depends on the selected index $p$.

\begin{lemma}[Euler-number ratios]\label{lem:euler-ratios}
For $\Dg\ge11$ and $0\le k\le\lfloor\Dg/2\rfloor$, put
$\ell=\Dg-2k$. Then
\begin{equation}\label{eq:euler-ratio-bound}
 \frac{E_\ell}{E_{\Dg}}
 \le \varkappa_\ell\left(\frac{\pi^2}{4}\right)^k
          \frac{\ell!}{\Dg!},
 \qquad
 \varkappa_\ell=\begin{cases}
 79/100,&\ell=0,\\
 5/4,&\ell=1,\\
 51/50,&\ell=2,3,\\
 501/500,&\ell=4,5,\\
 1001/1000,&\ell\ge6.
 \end{cases}
\end{equation}
\end{lemma}
\begin{proof}
We use the Euler-number expansion recorded in
\cite[Section~1, immediately after (1.10)]{StanleyAlt}:
\begin{equation}\label{eq:euler-dirichlet}
 \frac{E_\ell}{\ell!}
 =\frac4\pi\left(\frac2\pi\right)^\ell S_\ell,
 \qquad
 S_\ell=\sum_{j\ge0}
 \frac{(-1)^{j(\ell+1)}}{(2j+1)^{\ell+1}},\qquad \ell\ge1.
\end{equation}
This formula uses the positive Euler zigzag numbers of
\cref{eq:euler-normalization}, not the signed coefficients of the hyperbolic secant; the distinction is
spelled out after \cref{eq:euler-normalization}.
If $\ell$ is even, the alternating-series bounds give
$1-3^{-\ell-1}\le S_\ell\le1$. If $\ell$ is odd, all terms are positive,
so $S_\ell\ge1$. Consequently
\begin{equation}\label{eq:euler-denominator}
 S_{\Dg}\ge b_*\Def1-3^{-12}=\frac{531440}{531441}
 \qquad(\Dg\ge11).
\end{equation}
The exponent twelve is valid for odd $\Dg$ because $S_{\Dg}\ge1$, and
for even $\Dg$ because then $\Dg\ge12$.

For $\ell\ge1$, dividing the two exact formulas before estimating gives
\begin{equation}\label{eq:euler-ratio-exact}
 \frac{E_\ell}{E_{\Dg}}
 =\frac{\ell!}{\Dg!}\left(\frac2\pi\right)^{\ell-\Dg}
                       \frac{S_\ell}{S_{\Dg}}
 =\left(\frac{\pi^2}{4}\right)^k
        \frac{\ell!}{\Dg!}\frac{S_\ell}{S_{\Dg}}.
\end{equation}
Thus $\varkappa_\ell$ only bounds a quotient of two normalized reciprocal-power
sums. It does not hide a factorial, a power of $\pi$, or a dependence on
$k$. The near-unit constants for $\ell\ge2$ reflect the rapid decay of
the tails of these sums. For $q>1$, comparison with a decreasing integral yields
\[
 \sum_{j\ge1}(2j+1)^{-q}
 \le3^{-q}+\int_1^\infty(2x+1)^{-q}\dif x
 =3^{-q}+\frac{3^{1-q}}{2(q-1)}.
\]
Using $q=4,6,8$ for the first possible odd indices in the indicated ranges,
and the bound $S_\ell\le1$ for even indices, gives
\[
 S_\ell\le\frac{55}{54}\ (\ell\ge2),\qquad
 S_\ell\le\frac{7303}{7290}\ (\ell\ge4),\qquad
 S_\ell\le\frac{91871}{91854}\ (\ell\ge6).
\]

For example, at the first odd order $\ell=3$,
$1+3^{-4}+3^{-3}/6=55/54$; at $\ell=5$ the same calculation is
$1+3^{-6}+3^{-5}/10=7303/7290$; and at $\ell=7$ it is
$1+3^{-8}+3^{-7}/14=91871/91854$.
Increasing the exponent decreases every positive tail term. Even orders
need only the smaller bound one, so each displayed bound controls its
whole stated range. This is the reason the piecewise constants occur in
pairs of consecutive derivative orders.
For $\ell\ge2$, divide \cref{eq:euler-dirichlet} at indices $\ell$ and
$\Dg$. The powers of $2/\pi$ become $(\pi^2/4)^k$, and the quotient
$S_\ell/S_{\Dg}$ is at most the relevant displayed rational number divided
by $b_*$. These quotients satisfy
\begin{equation}\label{eq:euler-roundings}
 \frac{55}{54b_*}<\frac{51}{50},\qquad
 \frac{7303}{7290b_*}<\frac{501}{500},\qquad
 \frac{91871}{91854b_*}<\frac{1001}{1000}.
\end{equation}
For completeness, the positive differences in these three comparisons are,
respectively, $7863/5314400$, $5709/26572000$, and $11633/14308000$.

The cases $\ell=0,1$ should not be obtained by extending an absolutely
convergent formula beyond its stated range. Instead, use $E_0=E_1=1$
and \cref{eq:euler-dirichlet} only in the denominator. After removing the
factor $(\pi^2/4)^k\ell!/\Dg!$, the remaining factors are
$\pi/(4S_{\Dg})$ and $\pi^2/(8S_{\Dg})$, respectively.
The elementary bound $\pi<22/7$ gives
\[
 \frac{11}{14b_*}<\frac{79}{100},\qquad
 \frac{121}{98b_*}<\frac54;
\]
the differences are $159377/37200800$ and $797039/52081120$.
This proves every case of \cref{eq:euler-ratio-bound}.
\end{proof}

\begin{lemma}[Normalized derivative bounds]\label{lem:derivative-envelope}
For real $s$ and the normalization in \cref{eq:normalized-D},
\begin{equation}\label{eq:derivative-low-orders}
 |\widehat D_0(s)|\le\frac\pi4<\frac{79}{100},\qquad
 |\widehat D_1(s)|\le\frac{\pi^2}{8}<\frac54.
\end{equation}
For every integer $\ell\ge1$,
\begin{equation}\label{eq:derivative-envelope}
 |\widehat D_\ell(s)|\le1+\frac{R_0(s)}{\pi\sqrt\ell}.
\end{equation}
\end{lemma}
\begin{proof}
Direct differentiation of the kernel gives
\[
 |\widehat D_0(s)|=\frac{\sqrt{s^2+\pi^2}}{4\cosh(s/2)},\qquad
 |\widehat D_1(s)|=
 \frac{(s^2+\pi^2)|\tanh(s/2)|}{8\cosh(s/2)}.
\]
The Taylor series \cite[(4.33.2)]{DLMF} gives $\cosh(s/2)\ge1+s^2/8$. Squaring this inequality shows $\cosh^2(s/2)\ge1+s^2/4\ge1+s^2/\pi^2$, proving the first bound. For the second, $|\tanh(s/2)|\le1$ gives $|\widehat D_1(s)|\le\frac{s^2+\pi^2}{s^2+8}\le\frac{\pi^2}{8}.$
Both strict rational comparisons follow from $\pi<22/7$. For $\ell\ge1$, the pole series is absolutely convergent, so
\cref{eq:normalized-poles,eq:angle-coordinate} imply
\[
 |\widehat D_\ell(s)|
 \le1+\sum_{j\ge1}
 \left(1+\frac{4\pi^2j^2}{R_0^2}\right)^{-(\ell+1)/2}
 \le1+\frac{R_0}{2\pi}J_\ell,
 \quad
 J_\ell=\int_0^\infty(1+y^2)^{-(\ell+1)/2}\dif y.
\]
For the last comparison, if $G$ is the positive decreasing summand viewed
as a function of a real variable, then
$G(j)\le\int_{j-1}^{j}G(x)\dif x$. Summing over $j\ge1$ and
setting $y=2\pi x/R_0$ gives the factor $R_0/(2\pi)$.
We have $J_1=\pi/2<2$ and $J_2=1<2/\sqrt2$.
The recurrence for higher orders can be checked without a special-function
evaluation. For $\ell\ge1$,
\[
 \frac{\dif}{\dif y}\{y(1+y^2)^{-(\ell+1)/2}\}
 =(1+y^2)^{-(\ell+1)/2}
   -(\ell+1)y^2(1+y^2)^{-(\ell+3)/2}.
\]
The bracket vanishes at both endpoints of the integral. Also the integral of the second factor $y^2(1+y^2)^{-(\ell+3)/2}$ is $J_\ell-J_{\ell+2}$. Hence $0=J_\ell-(\ell+1)(J_\ell-J_{\ell+2})$, or $J_{\ell+2}=\ell J_\ell/(\ell+1)$. If $J_\ell\le2/\sqrt\ell$, this gives $J_{\ell+2}\le\frac{2\sqrt\ell}{\ell+1}\le\frac2{\sqrt{\ell+2}},$ because $\ell(\ell+2)\le(\ell+1)^2$. Induction from the two base orders proves the bound for both parities. This proves \cref{eq:derivative-envelope} without applying absolute values to the
conditionally convergent order-zero expansion.
\end{proof}

\subsection{Phase loss: preserving the reinforcing corrections}\label{sec:phase-bound}

The factors $R_0^{2k}$ in the weights appear to grow near an angular endpoint.
For the retained corrections, however, the phase discrepancy supplies a
factor $\sin^2\theta=\pi^2/R_0^2$. The first correction therefore loses its
distance dependence entirely; the second leaves only the single power
controlled by the window.
\begin{proposition}\label{prop:phase-bound}
Under the hypotheses of \cref{thm:central-sign},
\begin{equation}\label{eq:phase-budget}
 \Lphase<\frac{38311}{204800}<\frac{19}{100}.
\end{equation}
\end{proposition}
\begin{proof}
The identities
$1-\cos(2\theta)=2\sin^2\theta$ and
$1-\cos(4\theta)=8\sin^2\theta\cos^2\theta\le8\sin^2\theta$
are sufficient. In \cref{eq:euler-ratio-bound}, the orders of the first two
corrections are at least seven, so both use $\varkappa=1001/1000$. The definition of
the weights and \cref{eq:first-c} give
\begin{align}
 w_1(1-\cos2\theta)
 &\le\frac{\varkappa\pi^4}{8}
       \frac{\Dg-5}{6\Dg(\Dg-1)},\label{eq:phase-first}\\
 w_2(1-\cos4\theta)
 &\le\frac{\varkappa\pi^6}{32}\frac{15}{4}
 \frac{5\Dg^2-48\Dg+67}{360(\Dg-1)(\Dg-2)(\Dg-3)}.
 \label{eq:phase-second}
\end{align}
Here are the cancellations producing these powers and numerical factors. For $k=1$, \cref{eq:euler-ratio-bound} and $1-\cos2\theta=2\pi^2/R_0^2$ give $c_1(\Dg)\frac{\varkappa\pi^2}{4\Dg(\Dg-1)}\frac{R_0^2}{4}\frac{2\pi^2}{R_0^2}=\frac{\varkappa\pi^4}{8}\frac{c_1(\Dg)}{\Dg(\Dg-1)}.$ For $k=2$, the corresponding calculation, writing the four consecutive factors in the factorial quotient explicitly, is
\begin{align*}
c_2(\Dg)\frac{\varkappa\pi^4}{16\Dg(\Dg-1)(\Dg-2)(\Dg-3)}\frac{R_0^4}{16}\frac{8\pi^2}{R_0^2}
=\frac{\varkappa\pi^6}{32}\frac{c_2(\Dg)R_0^2}{\Dg(\Dg-1)(\Dg-2)(\Dg-3)}.
\end{align*}
The window replaces $R_0^2/\Dg$ by $15/4$. Thus the two rational
quotients in \cref{eq:phase-first,eq:phase-second} are the complete
remaining size dependence; neither estimate has taken a maximum over an
uncontrolled angular endpoint.
The first quotient decreases for $\Dg\ge11$: its derivative has the sign
of $-x^2+10x-5=20-(x-5)^2<0$. Its value at eleven is $1/110$.
The second quotient is less than $1/1600$ throughout the same range.
Indeed, after multiplying by the positive denominators, this is equivalent
to positivity of
\[
 360(x-1)(x-2)(x-3)-1600(5x^2-48x+67)
 =40(9x^3-254x^2+2019x-2734).
\]
With $v=x-13\ge-2$, the cubic on the right is $9v^2(v+2)+79\left(v-\frac{11}{79}\right)^2+\frac{28319}{79}>0.$ Finally, $\pi^2<10$ in \cref{eq:phase-first,eq:phase-second} yields
\[
 \Lphase<\frac{1001}{1000}
 \left(\frac{100}{8\cdot110}
       +\frac{1000}{32}\frac{15}{4\cdot1600}\right)
 =\frac{38311}{204800}.
\]
The difference between $19/100$ and this last rational number is
$601/204800$.
\end{proof}

\begin{remark}[A general phase-deficit bound]
\label{rem:markov-phase}
The two elementary identities used above are instances of a reusable
inequality. Write $T_k(\cos u)=\cos(ku)$ for the Chebyshev polynomial
\cite[Eq.~18.5.1]{DLMF}. The classical Markov inequality, in the form
stated in \cite[Theorem~1.1]{Pierzchala}, is
$\|P'\|_{[-1,1]}\le(\deg P)^2\|P\|_{[-1,1]}$.
Since $T_k(1)=1$ and $\|T_k\|_{[-1,1]}=1$, integration from
$x=\cos2\theta$ to one yields
\[
 0\le1-\cos(2k\theta)=1-T_k(\cos2\theta)
 \le k^2(1-\cos2\theta)=2k^2\sin^2\theta.
\]
This explains the natural quadratic cost of retaining additional harmonics.
For the two harmonics actually needed here, the direct identities are
shorter and give the same bounds. The general inequality is recorded to
make its role and provenance explicit, not to introduce another hypothesis
into the proof.
\end{remark}

\subsection{Farther poles: distance defeats the retained powers}\label{sec:far-bound}

The preceding argument used the phases of the nearest pair. For farther
pairs it is more efficient to discard the phase and use their larger
distances. The issue is uniformity: as the size increases, the allowed
window also expands. We therefore prove that the complete distance-power
bound decreases with the size, rather than checking it at a collection of
sample sizes.

\begin{proposition}\label{prop:far-bound}
Under the hypotheses of \cref{thm:central-sign},
\begin{equation}\label{eq:far-budget}
 \Lfar<\frac{10625479}{1638400000}<\frac7{1000}<\frac{13}{1000}.
\end{equation}
\end{proposition}
\begin{proof}
First, the retained weights have size-independent bounds
\begin{equation}\label{eq:retained-weight-bounds}
 w_0=1,\qquad w_1<\frac{1001}{2560},\qquad
 w_2<\frac{5005}{65536}.
\end{equation}
To see this, use \cref{eq:euler-ratio-bound}, the window, and $\pi^2<10$:
\[
 w_k\le\frac{1001}{1000}
 \left(\frac{75\Dg}{32}\right)^k
 c_k(\Dg)\frac{(\Dg-2k)!}{\Dg!},\qquad k=1,2.
\]
The factor $75/32$ has a single source:
\begin{equation}\label{eq:weight-scale-constant}
 \left(\frac{\pi^2}{4}\right)^k
 \left(\frac{R_0^2}{4}\right)^k
 =\left(\frac{\pi^2R_0^2}{16}\right)^k
 <\left(\frac{10(15\Dg/4)}{16}\right)^k
 =\left(\frac{75\Dg}{32}\right)^k.
\end{equation}
The Euler ratio supplies the first factor and the weight definition the second. We will reuse exactly this constant in the higher-correction estimate, rather than introduce a second scale for the same product. For $k=1$, the coefficient factor after taking out $(75/32)$ is
$(\Dg-5)/(6(\Dg-1))<1/6$. For $k=2$, it is $\frac{\Dg(5\Dg^2-48\Dg+67)}{360(\Dg-1)(\Dg-2)(\Dg-3)}<\frac1{72},$ because $5(x-1)(x-2)(x-3)-x(5x^2-48x+67)=6(3x-5)(x+1)>0$. Substitution gives \cref{eq:retained-weight-bounds}. In particular,
\[
 \frac{1001}{1000}\frac{75}{32}\frac16=\frac{1001}{2560},
 \qquad
 \frac{1001}{1000}\left(\frac{75}{32}\right)^2\frac1{72}
 =\frac{5005}{65536}.
\]
These estimates bound the corrections themselves. The previous phase estimate was smaller because it also used their cosine deficits. For $j\ge1$, the distance identity and window imply
\begin{equation}\label{eq:far-distance-window}
 \left(\frac{R_0}{R_j}\right)^{\Dg+1-2k}
 \le\left(1+\frac{48j(j+1)}{5\Dg}\right)^{-(\Dg+1-2k)/2},
 \qquad k=0,1,2.
\end{equation}

Indeed,
$R_j^2/R_0^2=1+4\pi^2j(j+1)/R_0^2$ is at least
$1+16\pi^2j(j+1)/(15\Dg)$, and $\pi^2>9$ makes the last
coefficient larger than $48/(5\Dg)$. Raising its reciprocal to the
positive power $(\Dg+1-2k)/2$ gives
\cref{eq:far-distance-window}. Notice that as $\Dg$ grows the base of
this bound decreases while its exponent increases. Neither change alone
settles monotonicity; that is why the logarithmic derivative is needed.
Fix $A=48j(j+1)/5$ and define
$F_k(x)=(1+A/x)^{-(x+1-2k)/2}$ for $x\ge11$.
Writing $u=A/x$, differentiation gives
\begin{equation}\label{eq:far-log-derivative}
 \frac{F_k'(x)}{F_k(x)}
 =-\frac12\left[
 \log(1+u)-\left(1+\frac{1-2k}{x}\right)\frac{u}{1+u}
 \right].
\end{equation}
For $k=1,2$, the bracket is positive because $\log(1+u)>u/(1+u)$ for $u>0$. For $k=0$, use the slightly stronger inequality $\log(1+u)\ge2u/(2+u)$: the difference vanishes at zero and has
derivative $u^2/((1+u)(2+u)^2)$. The bracket in \cref{eq:far-log-derivative} is then at least $\frac{u\{A(1-1/x)-2\}}{x(1+u)(2+u)}>0,$ since $A(1-1/x)-2\ge(96/5)(10/11)-2=170/11>0$. For $k=1,2$, the simpler difference can be verified directly as well:
\[
 \frac{\dif}{\dif u}\left(\log(1+u)-\frac{u}{1+u}\right)
 =\frac{u}{(1+u)^2}>0,
\]
with value zero at $u=0$. These two calculations prove strict decrease for all three retained orders and every pole index, not only for the nearest omitted pole. Thus every $F_k$ is decreasing, and \cref{eq:far-distance-window} is at most $\left(1+\frac{48j(j+1)}{55}\right)^{-(6-k)}.$ This is the desired reduction from every size to the smallest size eleven. Put $q=6-k\in\{4,5,6\}$. Retain the terms $j=1,2$ exactly. For $j\ge3$, replace $j(j+1)$ by $j^2$ and compare the decreasing summand with its integral on $[2,\infty)$. With $a=48/55$, $\int_2^\infty(1+ax^2)^{-q}\dif x\le\frac{(1+4a)^{1-q}}{4a(q-1)}.$ Indeed, $-\frac{\dif}{\dif x}(1+ax^2)^{1-q}=2a(q-1)x(1+ax^2)^{-q}\ge4a(q-1)(1+ax^2)^{-q}$ on this interval. We have proved the entirely rational bound
\begin{equation}\label{eq:Hq-bound}
 \sum_{j\ge1}\left(\frac{R_0}{R_j}\right)^{\Dg+1-2k}
 \le H_q\Def
 \left(\frac{55}{151}\right)^q+
 \left(\frac{55}{343}\right)^q+
 \frac{55}{192(q-1)}\left(\frac{55}{247}\right)^{q-1}.
\end{equation}
Direct integer comparisons give
\begin{equation}\label{eq:Hq-roundings}
 H_6<\frac{239}{100000},\qquad
 H_5<\frac{67}{10000},\qquad
 H_4<\frac{483}{25000}.
\end{equation}
The three summands in each comparison are certified separately in
\cref{tab:far-cross-products}; no numerical evaluation of a pole series is
involved.

The denominators in \cref{eq:Hq-bound} are also determined by the
comparison, rather than selected numerically. The two explicit poles have
$55+48\cdot1\cdot2=151$ and $55+48\cdot2\cdot3=343$;
the integral starts at $x=2$, so its boundary value uses
$55+48\cdot2^2=247$. Its prefactor is
$1/(4a(q-1))=55/[192(q-1)]$.
For example, for $q=6$ the three summands are
$55^6/151^6$, $55^6/343^6$, and $55^6/(960\cdot247^5)$.
The first row of the appendix verifies
\[
 2336\cdot151^6-10^6\cdot55^6=10096845504736>0;
\]
the next two rows give upper numerators $17$ and $32$ with the same
scale $10^6$. Their sum $2385$ is less than $2390$, proving
$H_6<239/100000$. The $q=5$ and $q=4$ rows repeat the same
integer comparison. This is a finite rational certificate for the integral
bound, not a truncation of the infinite pole sum without a remainder. Combining \cref{eq:far-loss,eq:retained-weight-bounds,eq:Hq-bound,eq:Hq-roundings}
now gives
\[
 \Lfar<\frac{239}{100000}
 +\frac{1001}{2560}\frac{67}{10000}
 +\frac{5005}{65536}\frac{483}{25000}
 =\frac{10625479}{1638400000}.
\]
The gap to $7/1000$ is $843321/1638400000$. We retain the larger budget
$13/1000$ in the theorem so that the three allocations have a simple stated
margin.
\end{proof}

\subsection{Remaining corrections: coefficients, endpoint orders, and a tail}\label{sec:high-bound}

The final loss is the only part of the proof in which a substantial finite
arithmetic certificate is useful. There are two different infinite ranges:
the size $\Dg$ varies for each fixed correction index $k$, and the correction
index itself is unbounded. We handle them differently. A monotonicity lemma
reduces the first range to one rational calculation per $k$, with two exact
endpoint formulas for derivative orders zero and one. A positive coefficient
majorant then controls the second range by a factorial tail.

The distinction is important. A table of sampled sizes could not prove a
uniform bound. Conversely, expanding every coefficient polynomial for every
$k$ would obscure why a finite table suffices. The following construction
makes both reductions explicit.

\subsubsection{An exact coefficient recurrence}\label{sec:coefficient-recurrence}

Let
\begin{equation}\label{eq:S-eta}
 S(v)\Def\frac{\sqrt v}{\sin\sqrt v},\qquad
 \log S(v)=\sum_{j\ge1}\eta_jv^j,
 \qquad \tau_j\Def\frac{(-1)^j}{(2j+1)!}.
\end{equation}
These identities are power-series identities at zero, so the apparent square
roots introduce no branch ambiguity. Since $S(v)^{-1}=\sum_{j\ge0}\tau_jv^j$,
logarithmic differentiation gives the triangular rational recurrence
\begin{equation}\label{eq:eta-recurrence}
 \eta_j=-\tau_j-\frac1j\sum_{i=1}^{j-1}i\eta_i\tau_{j-i}.
\end{equation}

For example, multiply
$\sum_{j\ge1}j\eta_jv^{j-1}=-\tau'(v)/\tau(v)$ by
$\tau(v)=\sum_{j\ge0}\tau_jv^j$, with $\tau_0=1$.
At degree $j-1$ this says
$j\eta_j+\sum_{i=1}^{j-1}i\eta_i\tau_{j-i}=-j\tau_j$.
Solving for $\eta_j$ proves the displayed recurrence and shows why only
already computed coefficients are needed. For $j=2$, for instance,
$\eta_2=-1/120+(1/6)^2/2=1/180$.
For example, $\eta_1=1/6$, $\eta_2=1/180$, $\eta_3=1/2835$, and
$\eta_4=1/37800$. The double-angle identity gives
$\cos\sqrt v=S(v)/S(4v)$. Therefore the defining series
\cref{eq:c-definition} becomes
\begin{equation}\label{eq:c-exponential}
 \sum_{k\ge0}c_k(x)v^k
 =\frac{S(v)^{x+3}}{S(4v)^2}
 =\exp\left(\sum_{j\ge1}\eta_j
       (x+3-2^{2j+1})v^j\right).
\end{equation}
If the series on the left is denoted by $C_x(v)$, differentiation says $C_x'(v)=C_x(v)\sum_{j\ge1}j\eta_j(x+3-2^{2j+1})v^{j-1}.$ The coefficient of $v^{k-1}$ on the left is $kc_k(x)$; on the right it is a finite convolution with $j=1,\ldots,k$. Therefore
\begin{equation}\label{eq:c-recurrence}
 c_0(x)=1,\qquad
 k c_k(x)=\sum_{j=1}^k j\eta_j
                  (x+3-2^{2j+1})c_{k-j}(x).
\end{equation}
In particular, $c_k$ is a rational polynomial of degree $k$, with leading
coefficient $1/(6^k k!)$. For the leading term, only $j=1$ can contribute
degree $k$ in the recurrence, which proves this assertion inductively. Write $c_k(x)=\sum_{i=0}^k a_{k,i}x^i$, with missing coefficients interpreted
as zero. The completely explicit coefficient form is
\begin{equation}\label{eq:a-recurrence}
 a_{0,0}=1,\qquad
 a_{k,i}=\frac1k\sum_{j=1}^k j\eta_j
 \bigl(a_{k-j,i-1}+(3-2^{2j+1})a_{k-j,i}\bigr).
\end{equation}
Thus every rational number used below is specified by factorials, finite
sums, and integer powers. As examples beyond \cref{eq:first-c},
\begin{align}
 c_3(x)&=\frac{35x^3-483x^2+1213x-285}{45360},\label{eq:c3}\\
 c_4(x)&=\frac{175x^4-3080x^3+10274x^2-1528x+2223}{5443200},\label{eq:c4}\\
 c_5(x)&=\frac{385x^5-8085x^4+31174x^3+10626x^2+28345x+6675}
 {359251200}.\label{eq:c5}
\end{align}
We use the absolute-coefficient majorant
\begin{equation}\label{eq:c-majorant}
 \widehat c_k(x)\Def\sum_{i=0}^k|a_{k,i}|x^i,
 \qquad |c_k(x)|\le\widehat c_k(x)\quad(x\ge0).
\end{equation}
The hat in this notation indicates coefficientwise majorization, not a
numerical approximation.

\subsubsection{Why one size controls each regular range}

The next elementary lemma is not specific to the present coefficient
polynomials. It is the mechanism that turns the finite cap table into a
uniform statement.

\begin{lemma}[Polynomial--factorial monotonicity]\label{lem:factorial-monotonicity}
Let $p(x)=\sum_{i=0}^k b_i x^i$, and put
$\widehat p(x)=\sum_{i=0}^k|b_i|x^i$. The function
\begin{equation}\label{eq:factorial-monotonicity}
 \frac{\widehat p(x)x^k}{\prod_{j=0}^{2k-1}(x-j)}
 =\left(\sum_{i=0}^k|b_i|x^{i-k}\right)
  \prod_{j=0}^{2k-1}(1-j/x)^{-1}
\end{equation}
is nonnegative and nonincreasing for $x>2k-1$. At an integer $x=\Dg\ge2k$,
its left side is $\widehat p(\Dg)\Dg^k(\Dg-2k)!/\Dg!$.
\end{lemma}
\begin{proof}
Each power $x^{i-k}$ is nonincreasing on this interval. Each factor $(1-j/x)^{-1}$ is positive and nonincreasing there as well. A finite sum or product of nonnegative nonincreasing functions is nonincreasing. The sign condition is doing real work here: $p(x)$ itself may have
cancellation and need not be monotone. After replacing coefficients by absolute values, the factorization isolates only decreasing factors. For example, $\frac{\dif}{\dif x}(1-j/x)^{-1}
=-\frac{j}{x^2}(1-j/x)^{-2}\le0,$ while $(i-k)x^{i-k-1}\le0$ for $0\le i\le k$. The restriction $x>2k-1$ makes every denominator positive, including the last one; no monotonicity claim across a pole is used. Factoring $x$ from each of the $2k$ denominator factors gives the stated identity.
\end{proof}
By \cref{eq:euler-ratio-bound}, the window, and $\pi^2<10$, each term of
\cref{eq:high-loss} is bounded by
\begin{equation}\label{eq:high-basic-bound}
 |w_{\Dg,k}|\,|\widehat D_{\Dg-2k}(s)|
 \le \varkappa_{\Dg-2k}
 \left(\frac{75\Dg}{32}\right)^k
 \frac{(\Dg-2k)!}{\Dg!}|c_k(\Dg)|\,
 |\widehat D_{\Dg-2k}(s)|.
\end{equation}
For derivative orders at least two, define
\begin{equation}\label{eq:regular-cap-data}
 \Dg^{\star}_k\Def\max(11,2k+2),\quad \ell_k\Def\Dg^{\star}_k-2k,
 \quad x_k\Def\frac{5\Dg^{\star}_k}{12\ell_k},\qquad
 P_k\Def \varkappa_{\ell_k}\left(\frac{75\Dg^{\star}_k}{32}\right)^k
             \frac{\ell_k!}{\Dg^{\star}_k!}\widehat c_k(\Dg^{\star}_k).
\end{equation}
For all integers $\Dg\ge\Dg^{\star}_k$, put $\ell\Def\Dg-2k.$ By the definition of $\Dg^{\star}_k$ in \cref{eq:regular-cap-data}, one has $\ell=\Dg-2k\ge \Dg^{\star}_k-2k=\ell_k \ge2.$ Thus the general normalized derivative estimate \cref{eq:derivative-envelope} applies to the derivative factor that remains on the right-hand side of \cref{eq:high-basic-bound}. We now derive the regular cap explicitly.  Starting from \cref{eq:high-basic-bound}, first use \cref{eq:c-majorant} to replace the coefficient by its absolute-coefficient majorant:
\[
\begin{aligned}
 |w_{\Dg,k}|\,|\widehat D_{\Dg-2k}(s)|
 &\le
 \varkappa_{\Dg-2k}
 \left(\frac{75\Dg}{32}\right)^k
 \frac{(\Dg-2k)!}{\Dg!}
 |c_k(\Dg)|\,
 |\widehat D_{\Dg-2k}(s)|
\\
 &\le
 \varkappa_{\ell}
 \left(\frac{75\Dg}{32}\right)^k
 \frac{\ell!}{\Dg!}
 \widehat c_k(\Dg)\,
 |\widehat D_{\ell}(s)|.
\end{aligned}
\]
The coefficient majorant alone is not yet uniform in $\Dg$; the factorial quotient and the derivative factor must still be incorporated. Since $\ell\ge2$, \cref{eq:derivative-envelope} gives $|\widehat D_{\ell}(s)|\le1+\frac{R_0}{\pi\sqrt{\ell}}.$ The window hypothesis in \cref{thm:central-sign} gives $R_0^2\le15\Dg/4$.  Hence $\frac{R_0}{\pi\sqrt{\ell}}\le
\sqrt{\frac{15\Dg}{4\pi^2\ell}}.$ Using $\pi^2>9$, $\frac{15}{4\pi^2}<\frac{15}{36}=\frac5{12},$
and therefore
\[
 |\widehat D_{\ell}(s)|
 \le
 1+\sqrt{\frac{5\Dg}{12\ell}}
 =
 1+\sqrt{\frac{5\Dg}{12(\Dg-2k)}}.
\]
Substituting this estimate into the preceding bound yields
\[
\begin{aligned}
 |w_{\Dg,k}|\,|\widehat D_{\Dg-2k}(s)|
 &\le
 \varkappa_{\ell}
 \left(\frac{75\Dg}{32}\right)^k
 \frac{\ell!}{\Dg!}
 \widehat c_k(\Dg)
 \left(
 1+\sqrt{\frac{5\Dg}{12\ell}}
 \right)
\\
 &=
 \varkappa_{\ell}
 \left(\frac{75}{32}\right)^k
 \left[
   \widehat c_k(\Dg)\Dg^k
   \frac{(\Dg-2k)!}{\Dg!}
 \right]
 \left(
 1+\sqrt{\frac{5\Dg}{12(\Dg-2k)}}
 \right).
\end{aligned}
\]
This factorization is the point at which \cref{lem:factorial-monotonicity} becomes applicable.  Indeed, apply that lemma with $p(x)=c_k(x)$ and $\widehat p(x)=\widehat c_k(x).$ Since $\Dg\ge\Dg^{\star}_k\ge2k+2>2k-1$, the lemma gives $\widehat c_k(\Dg)\Dg^k\frac{(\Dg-2k)!}{\Dg!}
\le\widehat c_k(\Dg^{\star}_k)(\Dg^{\star}_k)^k\frac{(\Dg^{\star}_k-2k)!}{\Dg^{\star}_k!}.$
Because $\Dg^{\star}_k-2k=\ell_k,$ this becomes $\widehat c_k(\Dg)\Dg^k\frac{(\Dg-2k)!}{\Dg!}
\le\widehat c_k(\Dg^{\star}_k)(\Dg^{\star}_k)^k\frac{\ell_k!}{\Dg^{\star}_k!}.$ The other two size-dependent factors are controlled at the same first admissible size.  First, the constants in
\cref{eq:euler-ratio-bound} satisfy
\[
 \varkappa_2=\varkappa_3=\frac{51}{50},
 \qquad
 \varkappa_4=\varkappa_5=\frac{501}{500},
 \qquad
 \varkappa_\ell=\frac{1001}{1000}
 \quad(\ell\ge6),
\]
and hence $\varkappa_\ell$ is nonincreasing for $\ell\ge2$.  Since $\ell\ge\ell_k$, $\varkappa_\ell\le\varkappa_{\ell_k}.$ Second, for real $x>2k$, $\frac{\dif}{\dif x}\left(\frac{x}{x-2k}\right)=-\frac{2k}{(x-2k)^2}<0.$ Thus $x/(x-2k)$ is strictly decreasing on this interval.  Since $\Dg\ge\Dg^{\star}_k>2k$, $\frac{\Dg}{\Dg-2k}\le\frac{\Dg^{\star}_k}{\Dg^{\star}_k-2k}=
\frac{\Dg^{\star}_k}{\ell_k}.$ Multiplication by $5/12$, followed by the increasing square-root map, therefore gives $\sqrt{\frac{5\Dg}{12(\Dg-2k)}}\le\sqrt{\frac{5\Dg^{\star}_k}{12\ell_k}}=
\sqrt{x_k}.$ All size-dependent factors in the regular estimate have therefore been bounded simultaneously by their values at the single first admissible size $\Dg^{\star}_k$.  Combining the three preceding comparisons gives
\begin{align}
 |w_{\Dg,k}|\,|\widehat D_{\Dg-2k}(s)|
 &\le
 \varkappa_{\ell_k}
 \left(\frac{75}{32}\right)^k
 \left[
   \widehat c_k(\Dg^{\star}_k)
   (\Dg^{\star}_k)^k
   \frac{\ell_k!}{\Dg^{\star}_k!}
 \right]
 (1+\sqrt{x_k})
 \nonumber\\
 &=
 P_k(1+\sqrt{x_k}),
 \qquad
 \Dg\ge\Dg^{\star}_k.
 \label{eq:regular-cap}
\end{align}
The last equality is exactly the definition of $P_k$ in \cref{eq:regular-cap-data}, since $\left(\frac{75}{32}\right)^k(\Dg^{\star}_k)^k=\left(\frac{75\Dg^{\star}_k}{32}\right)^k.$
Thus one calculation at $\Dg^{\star}_k$ controls the entire regular range $\Dg\ge\Dg^{\star}_k$.
The choice $\Dg^{\star}_k=\max(11,2k+2)$ encodes two independent requirements.  The uniform central-sign theorem begins at $\Dg=11$, while the argument above uses the regular derivative
envelope only when $\Dg-2k\ge2.$ Consequently, $k=3,4$ have first regular size $\Dg^{\star}_k=11$,
whereas $k\ge5$ has first regular size $\Dg^{\star}_k=2k+2$.  For a fixed $k$, the only admissible sizes not covered by \cref{eq:regular-cap} are therefore $\Dg=2k$ and $\Dg=2k+1$ provided those sizes are at least eleven.  These correspond exactly to derivative orders zero and one and are treated separately in the next subsection. It is essential that \cref{lem:factorial-monotonicity} is applied to the combined quantity $\widehat c_k(\Dg)\Dg^k\frac{(\Dg-2k)!}{\Dg!},$
rather than to $\widehat c_k(\Dg)$ alone.  The coefficient polynomial may grow with $\Dg$; the uniform decrease appears only after that growth is paired with the factorial quotient.

\subsubsection{The two exceptional derivative orders}\label{sec:exceptional-orders}

The majorant $\widehat c_k$ sacrifices cancellation among coefficients. At
$\Dg=2k$ and $2k+1$, that cancellation is particularly useful, so we evaluate
the coefficient exactly before taking an absolute value.

\begin{lemma}[Endpoint coefficient identities]\label{lem:endpoint-coefficients}
For $k\ge1$,
\begin{align}
 c_k(2k)&=-\frac{(2k-2)!}{2^{2k-1}k!(k-1)!},\label{eq:coefficient-order0}\\
 c_k(2k+1)&=-\frac{2^{2k-1}k!(k-1)!}{(2k+1)!}.\label{eq:coefficient-order1}
\end{align}
\end{lemma}
\begin{proof}
The coefficient definition, expressed as a residue at zero, gives
\[
 c_k(\Dg)
 =\Res_{u=0}\frac{\cos^2u\,u^{\Dg-2k}}{\sin^{\Dg+1}u}.
\]
Use the local analytic change $x=\sin u$, with $u=\arcsin x$ and
$\cos u=\sqrt{1-x^2}$, both on their branches at zero. One factor of
$\cos u$ is absorbed by $\dif x$, so
\begin{equation}\label{eq:arcsine-coefficient}
 c_k(\Dg)=[x^{\Dg}](\arcsin x)^{\Dg-2k}\sqrt{1-x^2}.
\end{equation}
For $\Dg=2k$, the binomial series of $\sqrt{1-x^2}$ gives
\cref{eq:coefficient-order0} directly.

In factorial form, the coefficient is
\[
 (-1)^k\binom{1/2}{k}
 =-\frac{1\cdot3\cdots(2k-3)}{2^k k!}
 =-\frac{(2k-2)!}{2^{2k-1}k!(k-1)!},
\]
where the product is empty when $k=1$. This also fixes its sign before
an absolute value is taken.

For $\Dg=2k+1$, put $B(x)=\arcsin x\sqrt{1-x^2}
=\sum_{k\ge0}b_kx^{2k+1}$. Differentiating gives
$(1-x^2)B'(x)=1-x^2-xB(x)$. Comparison of the constant and quadratic
coefficients gives $b_0=1$, $b_1=-1/3$. For $k\ge2$, comparison of the
coefficient of $x^{2k}$ gives
$(2k+1)b_k=(2k-2)b_{k-1}$. Multiplying these ratios yields
$b_k=-2^{2k-1}k!(k-1)!/(2k+1)!$, as required.
\end{proof}

Combining these identities with \cref{eq:high-basic-bound,eq:derivative-low-orders}
produces the rational exceptional bounds
\begin{align}
 T_{k,0}&\Def\left(\frac{79}{100}\right)^2
 \left(\frac{75(2k)}{32}\right)^k
 \frac{(2k-2)!}{(2k)!\,2^{2k-1}k!(k-1)!},\label{eq:exception0-cap}\\
 T_{k,1}&\Def\left(\frac54\right)^2
 \left(\frac{75(2k+1)}{32}\right)^k
 \frac{2^{2k-1}k!(k-1)!}{((2k+1)!)^2}.\label{eq:exception1-cap}
\end{align}
The first is needed only for $k\ge6$, and the second only for $k\ge5$,
because $\Dg\ge11$.

The squares in these formulas have distinct origins. At order zero,
$\varkappa_{0}=79/100$ bounds the Euler-number ratio and a second $79/100$
bounds $|\widehat D_0|$. At order one the corresponding two factors
are both $5/4$. The exact negative values of $c_k(2k)$ and
$c_k(2k+1)$ are first replaced by their positive magnitudes, and the
remaining factorials come from \cref{eq:high-basic-bound}.
In particular, the order-one case contains two copies of $(2k+1)!$
in its denominator: one from the Euler ratio and one from
\cref{eq:coefficient-order1}. For any fixed $\Dg$ and $k$, the regular case and
the two exceptional cases are mutually exclusive. Accordingly, a cap must
bound their maximum, not their sum.

\subsubsection{Sixteen finite caps}\label{sec:finite-caps}

\Cref{tab:correction-caps} specifies rational numbers $q_k$ and $U_k$ for
$3\le k\le18$. They satisfy
\begin{equation}\label{eq:cap-certificates}
 q_k^2>x_k,\qquad P_k(1+q_k)<U_k,
 \qquad T_{k,0}<U_k\ (k\ge6),\qquad T_{k,1}<U_k\ (k\ge5).
\end{equation}
All quantities here are rational; the square roots in
\cref{eq:regular-cap} have been eliminated by the first comparison.
The recurrence \cref{eq:eta-recurrence,eq:a-recurrence} specifies $P_k$
exactly. \Cref{app:certificates} records positive rational lower bounds
for all the differences in \cref{eq:cap-certificates}, together with the
finite algorithm that reproduces them. The accompanying coefficient ledger
provides the expanded polynomials, but they are not extra assumptions in the
certificate.

The numbers $q_k$ and $U_k$ are chosen rational upper bounds, not
intrinsic constants of the polynomial family. Their admissibility consists
exactly of the finitely many strict inequalities in
\cref{eq:cap-certificates}. For a proposed row, compute $\eta_j$ up to
$j=k$ by \cref{eq:eta-recurrence}, then compute $a_{k,i}$ by
\cref{eq:a-recurrence}, evaluate $\widehat c_k(\Dg^{\star}_k)$, and form
$P_k$. Checking $q_k^2>x_k$ uses integers after clearing denominators;
so does checking $U_k-P_k(1+q_k)>0$. The endpoint formulas supply the
other two checks with no additional coefficient computation. Because
\cref{eq:regular-cap} already controls every larger size, a successful
row certifies an infinite range rather than a list of experiments.

\begin{table}[tbp]
\centering\small
\caption{Uniform caps for the remaining correction indices. For each row,
$\Dg^{\star}_k$ and $\ell_k$ are the first regular size and derivative order,
$q_k>\sqrt{x_k}$ is rational, and $U_k$ bounds the regular term and each
applicable exceptional term.}\label{tab:correction-caps}
\begin{tabular}{rrrrrr}
\toprule
$k$&$\Dg^{\star}_k$&$\ell_k$&$x_k$&$q_k$&$U_k$\\\midrule
3&11&5&$11/12$&$23/24$&$53/200$\\
4&11&3&$55/36$&$26/21$&$221/1000$\\
5&12&2&$5/2$&$68/43$&$43/250$\\
6&14&2&$35/12$&$41/24$&$29/500$\\
7&16&2&$10/3$&$11/6$&$87/5000$\\
8&18&2&$15/4$&$31/16$&$47/10000$\\
9&20&2&$25/6$&$37/18$&$57/50000$\\
10&22&2&$55/12$&$11/5$&$13/50000$\\
11&24&2&$5$&$9/4$&$27/500000$\\
12&26&2&$65/12$&$5/2$&$3/250000$\\
13&28&2&$35/6$&$3$&$3/1000000$\\
14&30&2&$25/4$&$3$&$1/2000000$\\
15&32&2&$20/3$&$3$&$1/10000000$\\
16&34&2&$85/12$&$3$&$1/50000000$\\
17&36&2&$15/2$&$3$&$3/1000000000$\\
18&38&2&$95/12$&$3$&$1/1000000000$\\\bottomrule
\end{tabular}
\end{table}

Here is one complete regular calculation. From \cref{eq:c3},
$\widehat c_3(11)=824/315$, hence
\[
 P_3=\frac{501}{500}\left(\frac{825}{32}\right)^3
       \frac{5!}{11!}\frac{824}{315}
 =\frac{10406605}{77070336}.
\]
The two differences are
\[
 \left(\frac{23}{24}\right)^2-\frac{11}{12}=\frac1{576},\qquad
 \frac{53}{200}-\frac{10406605}{77070336}\left(1+\frac{23}{24}\right)
 =\frac{26422549}{46242201600}>0.
\]
There is no exceptional case at this index because $2k+1<11$.
For an index where an exceptional case does occur, $k=5$ gives
\[
 P_5=\frac{3391462929375}{50929722195968},\qquad
 \left(\frac{68}{43}\right)^2-\frac52=\frac3{3698},
\]
\[
 \frac{43}{250}-P_5\left(1+\frac{68}{43}\right)
 =\frac{27980025094291}{273747256803328000}>0,
\]
and
\[
 T_{5,1}=\frac{2599609375}{157840048128},\qquad
 \frac{43}{250}-T_{5,1}=\frac{3068609862877}{19730006016000}>0.
\]
Thus the same cap applies both at $\Dg=11$, where the derivative order is
one, and at every $\Dg\ge12$, where the monotonicity reduction applies.

The first index for which all three branches occur is $k=6$.
Here $\Dg^{\star}_6=14$, $\ell_6=2$, and $x_6=35/12$. The same recurrence
gives
\[
 \widehat c_6(14)=\frac{287383564043}{392302310400},\qquad
 P_6=\frac{3740476700747171875}{174834909637219713024}.
\]
With $q_6=41/24$ and $U_6=29/500$, the regular comparisons are
\begin{align*}
 q_6^2-x_6&=\frac1{576},\\
 U_6-P_6(1+q_6)
 &=\frac{2300083331189121677}{40346517608589164544000}>0.
\end{align*}
At sizes twelve and thirteen the exact coefficients are
$c_6(12)=-21/1024$ and $c_6(13)=-256/9009$, respectively.
Substituting into the two exceptional formulas gives
\begin{align*}
 U_6-T_{6,0}&=\frac{270769207307419}{6047313952768000}>0,\\
 U_6-T_{6,1}&=\frac{7990159113270307}{152789166587904000}>0.
\end{align*}
Thus one cap covers $\Dg=12$, $\Dg=13$, and every $\Dg\ge14$,
using a different justified bound in each branch. There is no $k=6$
term when $\Dg=11$. The appendix records the analogous positive margins
for every remaining row.

Summing the sixteen caps gives
\begin{equation}\label{eq:finite-cap-sum}
 \sum_{k=3}^{18}U_k=\frac{92446203}{125000000}
 <\frac{73957}{100000}.
\end{equation}
The difference in this last comparison is $47/125000000$. At a size for
which some of these correction indices are absent, summing all sixteen caps
only enlarges the bound.

\subsubsection{A positive majorant for every higher index}\label{sec:infinite-tail}

To finish, we need a bound that does not require computing $c_k$ for
$k\ge19$. We write $A\preceq B$ for coefficientwise inequality between
power series with nonnegative coefficients. The sine product
\cite[(4.22.1)]{DLMF}, expanded logarithmically near zero, gives
\begin{equation}\label{eq:eta-positive}
 \eta_j=\frac1j\sum_{n\ge1}\frac1{(\pi^2n^2)^j}>0.
\end{equation}

To see the coefficient extraction, substitute $u=\sqrt v$ in the product for sine and take the reciprocal: $S(v)=\prod_{n\ge1}\left(1-\frac{v}{\pi^2n^2}\right)^{-1}.$ Expanding $-\log(1-w)=\sum_{j\ge1}w^j/j$ in each factor and collecting $v^j$ gives \cref{eq:eta-positive}. For $|v|\le r<\pi^2$, the double absolute sum is bounded by $\sum_n -\log(1-r/(\pi^2n^2))<\infty$, which justifies the interchange. The positivity is now coefficientwise; a pointwise inequality
between two real analytic functions would not by itself justify the coefficient bound used next.
The convergence is absolute on every closed disk $|v|<\pi^2$, so the logarithm and the product can be expanded termwise. Since $\eta_1=1/6$ from \cref{eq:eta-recurrence}, $\eta_j\le\frac1{6j\pi^{2j-2}}\le\frac3{2j9^j}.$ In the first comparison, each summand of order $j$ is bounded by $\pi^{-2j+2}$ times its order-one summand; the second uses $\pi^2>9$.
Multiplying by $N$ and comparing logarithmic coefficients gives
\[
 N\log S(v)\preceq\sum_{j\ge1}\frac{3N}{2j}(v/9)^j
 =-\frac{3N}{2}\log(1-v/9).
\]
Exponentiation preserves this comparison: every coefficient of $\exp A=\sum_{m\ge0}A^m/m!$ is a sum of nonnegative products of coefficients of $A$. Hence
\begin{equation}\label{eq:S-positive-majorant}
 S(v)^N\preceq(1-v/9)^{-3N/2}.
\end{equation}

The absolute coefficients of $\cos^2\sqrt v$ are those of
$\cosh^2\sqrt v$. For $j\ge1$, their ratio to the coefficient $1/j!$
of $\ee^v$ is $r_j=2^{2j-1}j!/(2j)!$; here $r_1=1$ and
$r_{j+1}/r_j=2/(2j+1)\le1$. Thus
$\cosh^2\sqrt v\preceq\ee^v$.
Also $\ee^v\preceq(1-v/9)^{-9}$, since the logarithm of the latter
series is $v$ plus terms with nonnegative coefficients. Applying these
majorants directly to the coefficient definition, not to the quotient in
\cref{eq:c-exponential}, gives
\begin{equation}\label{eq:coefficient-tail-majorant}
 |c_k(\Dg)|\le\frac{(a)_k}{9^k k!},
 \qquad a=\frac{3\Dg}{2}+\frac{21}{2},
\end{equation}
where $(a)_k=a(a+1)\cdots(a+k-1)$.

In detail, coefficientwise absolute values of a product are bounded by
the product of the coefficientwise absolute values. The nonnegative
majorants for the two factors therefore give
\[
 |[v^k]\cos^2\sqrt v\,S(v)^N|
 \le[v^k](1-v/9)^{-9-3N/2}
 =\frac{(9+3N/2)_k}{9^k k!}.
\]
Since $N=\Dg+1$, the exponent is
$9+3(\Dg+1)/2=3\Dg/2+21/2$. Thus the constants $9$ and $21/2$
come from the cosine majorant and the shift in $N$, respectively.
This argument deliberately does not divide coefficientwise inequalities
in \cref{eq:c-exponential}; such division would not preserve their order. This step keeps the cancellation
question separate from positivity: $S$ has nonnegative coefficients by
\cref{eq:eta-positive}, and the cosine factor was majorized coefficientwise.

For $\Dg\ge2k$,
\begin{equation}\label{eq:pochhammer-bound}
 (a)_k\le(a+k-1)^k
 \le\left[\Dg\left(2+\frac{19}{4k}\right)\right]^k.
\end{equation}

Indeed, $a+k-1=3\Dg/2+k+19/2$. The hypothesis $\Dg\ge2k$
allows the two estimates $k\le\Dg/2$ and
$19/2\le19\Dg/(4k)$. Adding them gives precisely the second
inequality. It remains valid at the endpoint $\Dg=2k$, which is needed
because the factorial-tail argument must also cover derivative order zero.
We also need a factorial bound uniform even when $\Dg=2k$:
\begin{equation}\label{eq:factorial-tail-bound}
 \frac{(\Dg-2k)!}{\Dg!}\le\left(\frac{\ee}{\Dg}\right)^{2k}.
\end{equation}
Indeed, the average of the largest $2k$ terms of
$\log1,\ldots,\log\Dg$ is at least their full average, and
\[
 \log(\Dg!)\ge\int_1^{\Dg}\log x\dif x
 =\Dg\log\Dg-\Dg+1\ge\Dg(\log\Dg-1).
\]
Exponentiating the resulting lower bound for the product of those $2k$
largest integers proves \cref{eq:factorial-tail-bound}. Equivalently, the full chain of inequalities for that product is
\[
 \log\frac{\Dg!}{(\Dg-2k)!}
 =\sum_{j=\Dg-2k+1}^{\Dg}\log j
 \ge\frac{2k}{\Dg}\log(\Dg!)
 \ge2k(\log\Dg-1).
\]
All logarithms have positive integer arguments, including when $\Dg=2k$. No asymptotic Stirling remainder or lower bound away from that endpoint is being assumed. Every $\varkappa_\ell$ in \cref{eq:euler-ratio-bound} is at most $5/4$. Using the window and $\pi^2<10$ in the weight definition, then \cref{eq:coefficient-tail-majorant,eq:pochhammer-bound,eq:factorial-tail-bound},
we obtain $|w_{\Dg,k}|\le\frac5{4k!}\left[\frac{25\ee^2}{96}\left(2+\frac{19}{4k}\right)\right]^k.$ The powers of the size cancel visibly: \cref{eq:high-basic-bound} without
the derivative contributes $(75\Dg/32)^k$, the coefficient majorant contributes $\Dg^k(2+19/(4k))^k/(9^k k!)$, and the factorial quotient contributes $\ee^{2k}/\Dg^{2k}$. Their size-free factor is $(75/32)/9=25/96$. This cancellation is the reason a single tail bound
can hold for all sizes, even though the individual coefficient polynomials grow with $\Dg$.
For $k\ge19$, the second factor in parentheses is at most $9/4$. Furthermore $\ee^2<37/5$: the terms through order four in its exponential series sum to seven, while the tail starts with $4/15$ and its subsequent ratios are at most $1/3$, with a strict inequality after the first ratio. Consequently
\[
 \frac{25\ee^2}{96}\left(2+\frac{19}{4k}\right)
 <\frac{25}{96}\frac{37}{5}\frac94=\frac{555}{128}
 \qquad(k\ge19),
\]
and thus
\begin{equation}\label{eq:weight-tail-bound}
 |w_{\Dg,k}|\le\frac5{4k!}\left(\frac{555}{128}\right)^k,
 \qquad k\ge19,\quad\Dg\ge2k.
\end{equation}

It remains to bound the normalized derivative without introducing a factor
growing with $\Dg$. A coarse bound linear in $k$ is enough:
\begin{equation}\label{eq:derivative-tail-bound}
 |\widehat D_{\Dg-2k}(s)|\le\frac{7k}{3}\qquad(k\ge1,\ \Dg\ge2k).
\end{equation}
At $\Dg=2k$, use \cref{eq:derivative-low-orders}. At $\Dg\ge4k$, \cref{eq:derivative-envelope} gives a bound at most $1+\sqrt{5/6}<2$. In the remaining range $2k<\Dg<4k$, the order is at
least one, and the same lemma gives $1+\sqrt{5k/3}<1+4\sqrt k/3\le7k/3$. This proves
\cref{eq:derivative-tail-bound} in all cases. Combining \cref{eq:weight-tail-bound,eq:derivative-tail-bound}, and extending a finite sum to an infinite positive sum, gives
\begin{align}
 \sum_{k=19}^{\lfloor\Dg/2\rfloor}
 |w_{\Dg,k}|\,|\widehat D_{\Dg-2k}(s)|
 \le\frac{35}{12}\sum_{k\ge19}
       \frac{(555/128)^k}{(k-1)!}
 \le\frac{35}{12}\frac{(555/128)^{19}}{18!}
         \frac{2432}{1877}
 <\frac1{1000}.\label{eq:infinite-tail-budget}
\end{align}
The geometric comparison uses the successive ratio $(555/128)/k\le555/2432$ for $k\ge19$. More explicitly, let $a_k=(555/128)^k/(k-1)!$. Then
$a_{k+1}/a_k=555/(128k)\le555/2432<1$.
Therefore $a_{19+m}\le a_{19}(555/2432)^m$, and summing over
$m\ge0$ gives the factor
$(1-555/2432)^{-1}=2432/1877$.
The split at eighteen is useful because from nineteen onwards both this
geometric comparison and the simple bound $2+19/(4k)\le9/4$ hold.
It is a certified finite/tail split, not a claim that the coefficients
change their behavior at that index.
The last comparison is exact: if its middle expression is $A/B$ in lowest
terms, then
\begin{align*}
 A&=534920176379658889180566641261733446502685546875,\\
 B&=712306073018453567477086497596298503861623423238144,\\
 B-1000A&=177385896638794678296519856334565057358937876363144>0.
\end{align*}
For $\Dg<38$, the sum on the left is empty and the same bound is immediate.

\begin{proposition}\label{prop:high-bound}
Under the hypotheses of \cref{thm:central-sign},
\begin{equation}\label{eq:high-budget}
 \Lhigh<\frac{74057}{100000}<\frac34.
\end{equation}
\end{proposition}
\begin{proof}
For $3\le k\le18$, \cref{eq:regular-cap,eq:exception0-cap,eq:exception1-cap,eq:cap-certificates}
bound every applicable term by $U_k$. Summing
\cref{eq:finite-cap-sum} and \cref{eq:infinite-tail-budget} gives the first
inequality in \cref{eq:high-budget}. The gap to $3/4$ is $943/100000$.
\end{proof}

\subsection{Assembly of the central signs}\label{sec:central-assembly}

The estimates have now returned us to the geometric decomposition with which
we began. No estimate has been made on a quotient by the leading sine, so
there is no singular division near its zeros. Instead, we evaluate at its
extrema, where its magnitude is exactly one and where the first two
corrections have the controlled phase relation \cref{eq:phase-lock}. The three estimates can now be compared on the same normalization:
\begin{center}
\begin{tabular}{@{}lll@{}}\toprule
Term & Bound actually obtained & Allocated budget\\\midrule
$\Lphase$ & $38311/204800$ & $19/100$\\
$\Lfar$ & $10625479/1638400000$ & $13/1000$\\
$\Lhigh$ & $74057/100000$ & $3/4$\\\bottomrule
\end{tabular}
\end{center}
Every entry in the middle column has been derived from the window, a
monotonicity argument, or an exact rational certificate. The right-hand
column is a convenient enlargement of those bounds. Its sum is $953/1000$,
leaving the stated margin $47/1000$ against the unit principal sine.
The margin is not fitted to numerical values of $h_{\Dg}$ and does not
rely on knowing its zeros in advance.

\begin{proof}[Proof of \cref{thm:central-sign}]
At $s=s_p$, \cref{eq:loss-lower-bound} and
\cref{prop:phase-bound,prop:far-bound,prop:high-bound} imply
\[
 (-1)^p\Ncal_{\Dg}(s_p)
 \ge1-\Lphase-\Lfar-\Lhigh
 >1-\frac{19}{100}-\frac{13}{1000}-\frac34
 =\frac{47}{1000}.
\]
The sign relation in \cref{eq:N-sign} gives the asserted sign of
$h_{\Dg}$. Finally, \cref{lem:window} verifies the window hypothesis for
all indices in \cref{eq:selected-indices}.
\end{proof}

\begin{remark}[The margin is a certificate, not a measurement]\label{rem:margin}
The constant $47/1000$ comes from enlarging three rigorous bounds to round
budgets, and it is far smaller than the true slack. A floating-point
evaluation of the exact identity at every selected extremum for sizes
$11\le\Dg\le150$ gave $(-1)^p\Ncal_{\Dg}(s_p)\ge0.947$, with $\Lphase\le0.14$,
$\Lfar\le0.001$ and $\Lhigh\le0.006$ throughout. Most of the proved budget is
spent on $\Lhigh$ because it must hold for all sizes at once, via the
infinite tail; this is an artifact of making the argument uniform, and
it is not a sign of near failure. These numerical values are an illustration
and are not used in the proof.
\end{remark}

The central theorem gives a real zero between every two consecutive selected
points. It does not by itself prove simplicity, nor does it account for the
two intervals between the selected points and the already certified edge
regions. The final section addresses exactly those two issues: it places the
boundary points in the correct order, and then uses the degree to show that
there is no room for extra roots or extra multiplicity.

% ===== 08_joining =====
\section{Joining the intervals and exhausting the degree}\label{sec:joining}

We now have two kinds of certified information. The spectral argument counts
roots inside disks and gives signs at their negative real boundary points.
The central argument gives signs at a grid chosen from the pole geometry.
To combine them, signs alone are not enough: the two grids must occur in the
correct order. Otherwise two apparently different sign intervals could
contain the same root. We first prove a quantitative separation of their
boundary points, and only then perform the global count.

Fix $\Dg\ge11$, let $N=\Dg+1$, and write $K=K_N$. By
\cref{eq:edge-index}, $K\ge1$, with
\begin{equation}\label{eq:joining-size-ranges}
 K=1\quad\Longleftrightarrow\quad11\le\Dg\le18,
 \qquad K\ge2\quad\Longleftrightarrow\quad\Dg\ge19.
\end{equation}
Set
\begin{equation}\label{eq:tp}
 t_p\Def-\ee^{-s_p},\qquad K+1\le p\le\Dg-K-1.
\end{equation}
The selected range contains at least two indices by \cref{lem:window}.
Since $\theta_p$ increases with $p$, $s_p=\pi\cot\theta_p$ decreases,
and $s\mapsto-\ee^{-s}$ increases. Therefore
\begin{equation}\label{eq:tp-order-sign}
 t_{p+1}<t_p<0,\qquad
 \sgn h_{\Dg}(t_p)=(-1)^{p+1}.
\end{equation}
The sign statement is \cref{thm:central-sign}; its window hypothesis has
already been checked in \cref{lem:window} and is not an additional assumption
in the joining argument. Reflection of the angle also gives
\begin{equation}\label{eq:tp-reciprocity}
 \theta_{\Dg-p}=\pi-\theta_p,\qquad s_{\Dg-p}=-s_p,
 \qquad t_{\Dg-p}=t_p^{-1}.
\end{equation}

\subsection{Boundary separation in the logarithmic coordinate}\label{sec:boundary-separation}

Assume first that $K\ge2$, and abbreviate
$a=a_{\Dg,K}$. The near-zero region already counted is $(-a,0)$.
It is more convenient to compare the new boundary point with $-a$ in the
logarithmic coordinate: the spectral radius becomes an integral of the
spectral logarithm, while the central point is a cotangent. Their natural
scales are neighboring half-integers.

\begin{proposition}[Separation of the two grids]\label{prop:joining-separation}
For $K\ge2$,
\begin{equation}\label{eq:joining-order}
 -a^{-1}<t_{\Dg-K-1}<t_{K+1}<-a<0.
\end{equation}
More precisely, with $b=K+1/2$ and $s_{\mathrm{edge}}=-\log a$,
\begin{equation}\label{eq:joining-quantitative}
 s_{\mathrm{edge}}-s_{K+1}
 >\frac{37N}{72b(b+1)}>0.
\end{equation}
\end{proposition}
\begin{proof}
The balanced-radius definition \cref{eq:balanced-radius}, with
$\phi=\log\lambda$ from \cref{eq:continuous-leading}, gives
\begin{equation}\label{eq:edge-integral}
 s_{\mathrm{edge}}
 =\frac N2\bigl(\phi(K+1)-\phi(K-1)\bigr)
 =\frac N2\int_{K-1}^{K+1}\phi'(x)\dif x.
\end{equation}
Here $b\ge5/2$, and $x+1/2$ runs over $[b-1,b+1]$. Put $\delta\Def\frac5{16(b-1)^2}\le\frac5{36}.$
The derivative in \cref{eq:phi-derivatives} is $\phi'(x)=\vartheta(x)\cot\vartheta(x)/(x+1/2)$, where $\vartheta(x)=\pi/(4(x+1/2))$. On this integration interval, $0<\vartheta(x)\le\pi/6$. Since $\sin u\le u$ and $\cos u>0$ there,
\[
 \vartheta(x)\cot\vartheta(x)\ge\cos\vartheta(x)
 \ge1-\frac{\vartheta(x)^2}{2}
 >1-\frac5{16(b-1)^2}=1-\delta.
\]
The strict final bound uses $\pi^2<10$. Thus
\begin{equation}\label{eq:edge-integral-lower}
 s_{\mathrm{edge}}
 >\frac{N(1-\delta)}2\int_{b-1}^{b+1}\frac{\dif y}{y}
 >\frac{N(1-\delta)}b.
\end{equation}
For the last inequality, pair $y=b-v$ and $y=b+v$: for $0<v\le1$, $1/(b-v)+1/(b+v)=2b/(b^2-v^2)>2/b$, and integrate from zero to one. On the other hand, the first selected angle is
$\theta_{K+1}=(b+1)\pi/N\in(0,\pi/2)$ by \cref{lem:window}. Since $\tan u>u$ for $0<u<\pi/2$,
\begin{equation}\label{eq:central-boundary-upper}
 s_{K+1}=\pi\cot\left(\frac{(b+1)\pi}{N}\right)
 <\frac{N}{b+1}.
\end{equation}
The function $(b+1)/(b-1)^2$ decreases for $b>1$: its derivative is
$-(b+3)/(b-1)^3<0$. Hence $\delta(b+1)=\frac{5(b+1)}{16(b-1)^2}\le\frac{35}{72} \implies 1-\delta(b+1)\ge\frac{37}{72}.$ Subtracting \cref{eq:central-boundary-upper} from \cref{eq:edge-integral-lower} now gives $s_{\mathrm{edge}}-s_{K+1}>\frac{N\{1-\delta(b+1)\}}{b(b+1)} \ge\frac{37N}{72b(b+1)}.$ Since $s\mapsto-\ee^{-s}$ increases, this proves $t_{K+1}<-a$.
Applying the decreasing reciprocal map on the negative real axis and using \cref{eq:tp-reciprocity} gives $-a^{-1}<t_{\Dg-K-1}$. The middle inequality in \cref{eq:joining-order} follows from the strict ordering in \cref{eq:tp-order-sign}.
\end{proof}

The separation is stronger than mere nonintersection, but the proof uses
only what is needed for the root count. It also explains why choosing the
central grid one step beyond the spectral cutoff works: the central upper
bound has denominator $b+1$, whereas the spectral lower bound has denominator
$b$ with a relative error of order $(b-1)^{-2}$. That error is smaller than
the separation between the two scales.

\subsection{The two boundary intervals}\label{sec:boundary-intervals}

The spectral sign in \cref{eq:edge-count-sign} gives
\begin{equation}\label{eq:edge-boundary-signs}
 \sgn h_{\Dg}(-a)=(-1)^{K-1},\qquad
 \sgn h_{\Dg}(-a^{-1})=(-1)^{\Dg+K-3}.
\end{equation}
The second equality uses the reciprocal identity
\cref{eq:h-reciprocity}, including the sign of $(-a)^{\Dg-2}$.
The near central boundary has sign
$\sgn h_{\Dg}(t_{K+1})=(-1)^{K+2}$, opposite to the first sign in
\cref{eq:edge-boundary-signs}. The far central boundary has sign
$\sgn h_{\Dg}(t_{\Dg-K-1})=(-1)^{\Dg-K}$, opposite to the second
sign, since the two exponents differ by the odd integer $2K-3$.
Consequently, in addition to the intervals between successive central points,
the two intervals
\begin{equation}\label{eq:boundary-intervals}
 (-a^{-1},t_{\Dg-K-1}),\qquad (t_{K+1},-a)
\end{equation}
each contain a real zero. By \cref{prop:joining-separation}, these intervals
are nonempty and disjoint from all the edge intervals.

\begin{proposition}[Complete interval count]\label{prop:complete-count}
For every $\Dg\ge11$, the polynomial $h_{\Dg}$ has
$\Dg-2$ distinct negative real zeros.
For $\Dg\ge19$, the central zeros lie one in each of
\begin{equation}\label{eq:central-root-intervals}
 \begin{gathered}
 (-a^{-1},t_{\Dg-K-1}),\qquad
 (t_{p+1},t_p)\quad(K+1\le p\le\Dg-K-2),\\
 (t_{K+1},-a),
 \end{gathered}
\end{equation}
and the other zeros lie in the edge intervals of \cref{thm:edges} and their
reciprocals. Every listed interval contains exactly one zero, and that zero
is simple.
\end{proposition}
\begin{proof}
When $K\ge2$, \cref{thm:edges} accounts for the $K-1$ zeros in $|t|<a$
and the $K-1$ zeros in $|t|>a^{-1}$, and by
\cref{eq:central-annulus-count} the annulus $a<|t|<a^{-1}$ contains exactly
$\Dg-2K$ zeros, counted with multiplicity. Between the selected central
points there are $\Dg-2K-2$ intervals, each with opposite endpoint signs by
\cref{eq:tp-order-sign}. The two intervals in \cref{eq:boundary-intervals}
contribute two more. These $\Dg-2K$ intervals are pairwise disjoint and lie
in $(-a^{-1},-a)$, hence in the central annulus, by
\cref{prop:joining-separation}. \Cref{lem:sign-count}, with $U$ the central
annulus and $n=\Dg-2K$, shows that each contains exactly one zero, that it
is simple, and that the annulus contains no other zero. The three contour
counts of \cref{thm:edges} have thus been matched by sign intervals:
\begin{equation}\label{eq:degree-exhaustion}
 (K-1)+(\Dg-2K-2)+2+(K-1)=\Dg-2.
\end{equation}

When $K=1$, no edge circle is required, and we take $U=\CC$, where the
number of zeros is the degree $\Dg-2$. The selected points have indices
$2,\ldots,\Dg-2$, and the intervals between consecutive points number
$\Dg-4$. At the near end, $h_{\Dg}(t_2)<0$ and
$h_{\Dg}(0)=1$, giving a sign change on $(t_2,0)$.
At the far end, the monic degree $\Dg-2$ established in
\cref{prop:reciprocity} gives
$h_{\Dg}(t)/t^{\Dg-2}\longrightarrow1$ as $t\to-\infty$.
Thus a sufficiently negative $t_*<t_{\Dg-2}$ has sign
$\sgn h_{\Dg}(t_*)=(-1)^{\Dg-2}$, opposite to
$\sgn h_{\Dg}(t_{\Dg-2})=(-1)^{\Dg-1}$. These are
$\Dg-4+2=\Dg-2$ disjoint negative intervals with sign changes, and
\cref{lem:sign-count} with $n=\Dg-2$ applies.

In both cases every listed interval contains exactly one zero, that zero is
simple, and there are no other roots, real or complex.
\end{proof}

\subsection{The main theorem and coefficient consequences}\label{sec:consequences}

\begin{proof}[Proof of \cref{thm:main}]
The sizes $\Dg=1,2$ give $h_{\Dg}=1$. The sizes $3\le\Dg\le10$ are
covered by \cref{prop:finite}. All sizes $\Dg\ge11$ are covered by
\cref{prop:complete-count}. Finally, $Z_{\Dg}(t)=t h_{\Dg}(t)$ and
$h_{\Dg}(0)=1$, so multiplication by $t$ adds exactly one simple zero
at zero and leaves the strictly negative simple zeros unchanged.
\end{proof}

\begin{proof}[Proof of \cref{cor:main-lc}]
Write $m=\Dg-2$ for $\Dg\ge2$. By the theorem and monicity,
\[
 h_{\Dg}(t)=\prod_{j=1}^m(t-\rho_j),\qquad \rho_j<0,
\]
with the empty product interpreted as one when $m=0$. Hence every
coefficient $\mathsf h_j$ is positive. Newton's inequalities in the normalized form
\cite[Lemma~1.1]{Branden} state
\[
 \left(\frac{\mathsf h_j}{\binom mj}\right)^2
 \ge\frac{\mathsf h_{j-1}}{\binom m{j-1}}
     \frac{\mathsf h_{j+1}}{\binom m{j+1}},\qquad 1\le j\le m-1.
\]
Their hypotheses are satisfied because $h_{\Dg}$ has only real zeros and
nonnegative coefficients. Multiplying by $\binom mj^2$ gives
\begin{equation}\label{eq:strict-newton}
 \mathsf h_j^2\ge \mathsf h_{j-1}\mathsf h_{j+1}
 \frac{(j+1)(m-j+1)}{j(m-j)}
 >\mathsf h_{j-1}\mathsf h_{j+1}.
\end{equation}
The multiplier exceeds one and the adjacent coefficients are positive, which
explains strictness even though Newton's displayed inequality is non-strict.
For $m<2$ there are no interior coefficient indices. This proves both
normalized log-concavity and strict ordinary log-concavity in every
nonvacuous case. Multiplication by $t$ preserves the inequalities for the
nonzero coefficient sequence of $Z_{\Dg}$.
\end{proof}

The reciprocal symmetry of the coefficients, already established before the
analytic proof, combines with log-concavity to give unimodality. This last
consequence does not replace the earlier combinatorial results about
unimodality and $\gamma$-nonnegativity; rather, it records how the stronger
root statement controls the coefficient sequence.

\subsection{What the two representations contribute}

The proof began with a local constraint on an ordered list of visits. The
polyhedral interpretation identified the polynomial family, while the matrix
interpretation retained enough quantitative information to locate the extreme
zeros. In the logarithmic coordinate, the same family admitted a different
exact representation: a finite combination of derivatives of a kernel with
known poles. The conjugate-pole phases then supplied a central grid that the
coefficient-dominance argument could not provide.

Several steps of this route are useful beyond the particular family. The
annular certification lemma requires a real-coefficient series with suitably
dominant successive monomials. The centered transform requires a polynomial
coefficient sequence and an explicit exponential kernel. The pole-pair
calculation explains an oscillatory phase whenever a conjugate pair of simple
poles controls a derivative expansion. Finally, polynomial--factorial
monotonicity and a positive coefficient majorant provide a way to turn finitely
many rational comparisons into an estimate uniform over the size.
These tools do not assert that every combinatorial family has the requisite
spectral curvature or the same three losses. Here, the locality of the gap
constraints and the exact centered expansion provide those additional
structures. Keeping those inputs explicit is what makes the method both
verifiable in this problem and available for adaptation to another one.

\appendix

% ===== 09_certificates =====
\section{Exact arithmetic certificates}\label{app:certificates}

The arithmetic in the proof has a finite and a uniform part. The uniform
part is proved in the main text: \cref{lem:factorial-monotonicity} controls
all regular sizes, \cref{lem:endpoint-coefficients} handles the two endpoint
orders, and \cref{eq:infinite-tail-budget} controls all correction indices
beyond eighteen. This appendix supplies the finite comparisons to which
those arguments reduce. All displayed decimals are avoided; powers of ten
below denote exact rational numbers.

\subsection{The farther-pole roundings}

For $q=4,5,6$, the three terms of $H_q$ in \cref{eq:Hq-bound} have the form
$55^q/d$. \Cref{tab:far-cross-products} gives an integer $b$ such that
$55^q/d<b/10^6$, together with the positive difference
$b d-10^6 55^q$. Thus each row is independently checked by integer
exponentiation and subtraction.

\begin{table}[htbp]
\centering\small
\caption{Integer certificates for the three terms of each farther-pole bound.
The last column certifies $55^q/d<b/10^6$.}\label{tab:far-cross-products}
\begin{tabular}{rrrl}
\toprule
$q$&$d$&$b$&$b d-10^6 55^q$\\\midrule
6&$151^6$&2336&10096845504736\\
6&$343^6$&17&2390539477633\\
6&$960\cdot247^5$&32&562044077935040\\\addlinespace
5&$151^5$&6412&75102515412\\
5&$343^5$&107&4704706563901\\
5&$768\cdot247^4$&177&2682749738816\\\addlinespace
4&$151^4$&17602&401348802\\
4&$343^4$&662&12307127062\\
4&$576\cdot247^3$&1055&6640432640\\\bottomrule
\end{tabular}
\end{table}
For $q=6,5,4$, respectively, the three integer numerators sum to
$2385$, $6696$, and $19319$. These are less than $2390$, $6700$, and
$19320$, giving exactly the three bounds in \cref{eq:Hq-roundings}.

\subsection{The correction-cap margins}

To specify the finite verification without printing polynomials with hundreds
of digits on every line, define the exact differences
\begin{equation}\label{eq:cap-margin-definitions}
 \begin{aligned}
 \Delta_{k,\mathrm{reg}}&\Def U_k-P_k(1+q_k),\\
 \Delta_{k,0}&\Def U_k-T_{k,0},\qquad
 \Delta_{k,1}\Def U_k-T_{k,1}.
 \end{aligned}
\end{equation}
Here $U_k,q_k$ are in \cref{tab:correction-caps}; $P_k$ is defined by
\cref{eq:regular-cap-data}; its coefficients are specified by
\cref{eq:eta-recurrence,eq:a-recurrence}; and the exceptional terms are
\cref{eq:exception0-cap,eq:exception1-cap}. Thus the differences are fully
specified rational numbers, not floating-point residuals.

\Cref{tab:cap-margins} records the exact square difference and a strictly
positive rational lower bound for each applicable $\Delta$. For example,
the first regular entry means
$\Delta_{3,\mathrm{reg}}>57/10^5$, consistent with its full fraction
computed in \cref{sec:finite-caps}. A dash means that the corresponding
size is below eleven and is not part of the uniform argument.

\begin{table}[htbp]
\centering\small
\caption{Positive margins for all sixteen correction caps. Entries in the
last three columns are strict lower bounds for the differences in
\cref{eq:cap-margin-definitions}.}\label{tab:cap-margins}
\begin{tabular}{r r r r r}
\toprule
$k$&$q_k^2-x_k$&$\Delta_{k,\mathrm{reg}}>$&$\Delta_{k,0}>$&$\Delta_{k,1}>$\\\midrule
3 & $1/576$ & $57\cdot10^{-5}$ & -- & -- \\
4 & $1/196$ & $31\cdot10^{-5}$ & -- & -- \\
5 & $3/3698$ & $10\cdot10^{-5}$ & -- & $15\cdot10^{-2}$ \\
6 & $1/576$ & $57\cdot10^{-6}$ & $44\cdot10^{-3}$ & $52\cdot10^{-3}$ \\
7 & $1/36$ & $28\cdot10^{-6}$ & $12\cdot10^{-3}$ & $15\cdot10^{-3}$ \\
8 & $1/256$ & $44\cdot10^{-6}$ & $31\cdot10^{-4}$ & $41\cdot10^{-4}$ \\
9 & $19/324$ & $56\cdot10^{-8}$ & $68\cdot10^{-5}$ & $10\cdot10^{-4}$ \\
10 & $77/300$ & $26\cdot10^{-7}$ & $13\cdot10^{-5}$ & $22\cdot10^{-5}$ \\
11 & $1/16$ & $70\cdot10^{-8}$ & $23\cdot10^{-6}$ & $46\cdot10^{-6}$ \\
12 & $5/6$ & $47\cdot10^{-8}$ & $49\cdot10^{-7}$ & $10\cdot10^{-6}$ \\
13 & $19/6$ & $53\cdot10^{-8}$ & $14\cdot10^{-7}$ & $26\cdot10^{-7}$ \\
14 & $11/4$ & $68\cdot10^{-9}$ & $18\cdot10^{-8}$ & $43\cdot10^{-8}$ \\
15 & $7/3$ & $29\cdot10^{-9}$ & $40\cdot10^{-9}$ & $88\cdot10^{-9}$ \\
16 & $23/12$ & $90\cdot10^{-10}$ & $92\cdot10^{-10}$ & $18\cdot10^{-9}$ \\
17 & $3/2$ & $14\cdot10^{-10}$ & $11\cdot10^{-10}$ & $26\cdot10^{-10}$ \\
18 & $13/12$ & $77\cdot10^{-11}$ & $70\cdot10^{-11}$ & $95\cdot10^{-11}$ \\
\bottomrule
\end{tabular}
\end{table}

The following is the exact coefficient-generation algorithm behind these
comparisons. Arrays are indexed from zero, and absent coefficients are zero.
Its correctness follows directly from \cref{eq:eta-recurrence,eq:a-recurrence}:
the outer induction computes $\eta_1,\ldots,\eta_{18}$; the second computes
the polynomials in increasing $k$, so every quantity on the right has already
been determined. No root-finding, special-function evaluation, or numerical
optimization occurs in this calculation.
\newpage

\begin{lstlisting}[language=Python,basicstyle=\ttfamily\footnotesize,
frame=single,breaklines=true,columns=fullflexible,
title={Exact rational generation of the coefficient polynomials.},
label={lst:coefficients}]
from fractions import Fraction as F
from math import factorial

tau = [F((-1)**j, factorial(2*j+1)) for j in range(19)]
eta = [F(0)] * 19
for j in range(1, 19):
    eta[j] = -tau[j] - sum(
        (i*eta[i]*tau[j-i] for i in range(1, j)), F(0)
    ) / j

coefficients = [[F(1)]]
for k in range(1, 19):
    p = [F(0)] * (k+1)
    for j in range(1, k+1):
        for i, value in enumerate(coefficients[k-j]):
            scale = F(j, k) * eta[j] * value
            p[i] += (3 - 2**(2*j+1)) * scale
            p[i+1] += scale
    coefficients.append(p)

# c_k(x) = sum_i coefficients[k][i] * x**i
# Its absolute-coefficient majorant uses abs(coefficients[k][i]).
\end{lstlisting}

For each row, evaluate the majorant at $\Dg^{\star}_k$, multiply by the rational
factors defining $P_k$, and subtract from $U_k$. For the exceptional columns,
substitute the factorial expressions directly. Positivity can be tested by
cross multiplication of positive denominators. The complete companion
\texttt{verification/verify\_certificates.py} executes these steps, checks
\cref{eq:cap-certificates}, regenerates the small-size polynomials and sign
values, and checks the rational tail comparison. It uses only exact integer
and rational arithmetic from the Python standard library.
The files \texttt{supplement/coefficient\_ledger.json} and
\texttt{verification/exact\_margins.json} list every coefficient and every
full rational margin. They provide expanded data for inspection; the
recurrences, table, and formulas in the paper specify the same finite
certificate independently of those files. For access to these files, kindly contact the author.

\subsection{Recovering the finite-size certificates from integer counts}

The gap recurrence \cref{eq:integer-gap-recurrence} is a prefix-sum recurrence
on integers. It computes $\Omega_{\Dg}(M)$ without enumerating paths or
assuming real-rootedness. The numerator coefficients then follow from the
finite-difference formula \cref{eq:finite-extraction}. Finally, the identity
\cref{eq:reciprocal-reduction} can be checked as a polynomial identity: a term
$q_jx^j$ of $\Upsilon$ contributes
\[
 q_j(-1)^{m+j}t^{m-j}(1+t)^{2j+\varepsilon}
\]
to $h_{\Dg}(t)$. Its exponents are all nonnegative because $j\le m$.
Expanding by the binomial theorem reproduces every coefficient in
\cref{tab:finite-coefficients}; evaluating $\Upsilon$ at the integer endpoints
reproduces \cref{tab:finite-signs}.

These checks have deliberately different roles. The integer recurrence
verifies the polynomial being studied; reciprocal reduction verifies the
change of variable; the alternating signs and degree give the root
certificate. A plot is not used in any of these three steps.

% ===== 10_figures =====
\section{Additional geometric diagnostics}\label{app:figures}

The figures in this appendix compare the exact polynomial family with the
leading sine model. They are numerical illustrations of identities proved
in the main text, not additional hypotheses in the real-rootedness proof.
The low-size root certificates remain the integer evaluations in
\cref{tab:finite-signs}; the uniform certificates remain the inequalities in
\cref{sec:central}.

\begin{figure}[h]
 \centering
 \includegraphics[width=.9\linewidth]{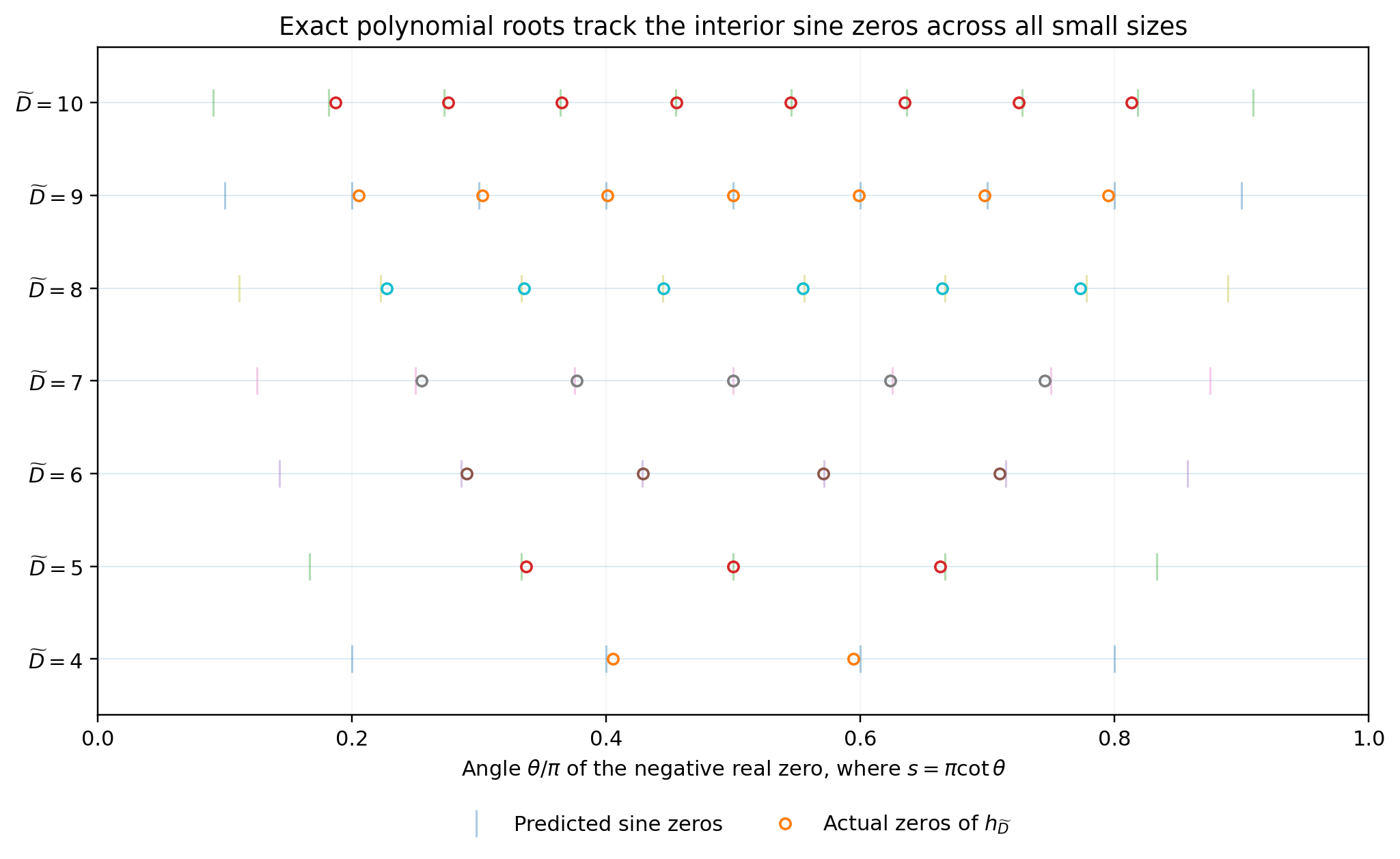}
 \caption{Angular root locations for the small sizes. The horizontal
 coordinate is $\theta/\pi$, where $t=-\exp(-\pi\cot\theta)$.
 Numerical root markers are compared with the zeros suggested
 by the nearest-pole sine. The central correspondence is visible even when
 the uniform window used in the proof is unavailable. Near the endpoints,
 the omitted terms cannot be ignored.}\label{fig:root-atlas}
\end{figure}

\begin{figure}[h]
 \centering
 \includegraphics[width=.8\linewidth,height=.85\textheight,keepaspectratio]{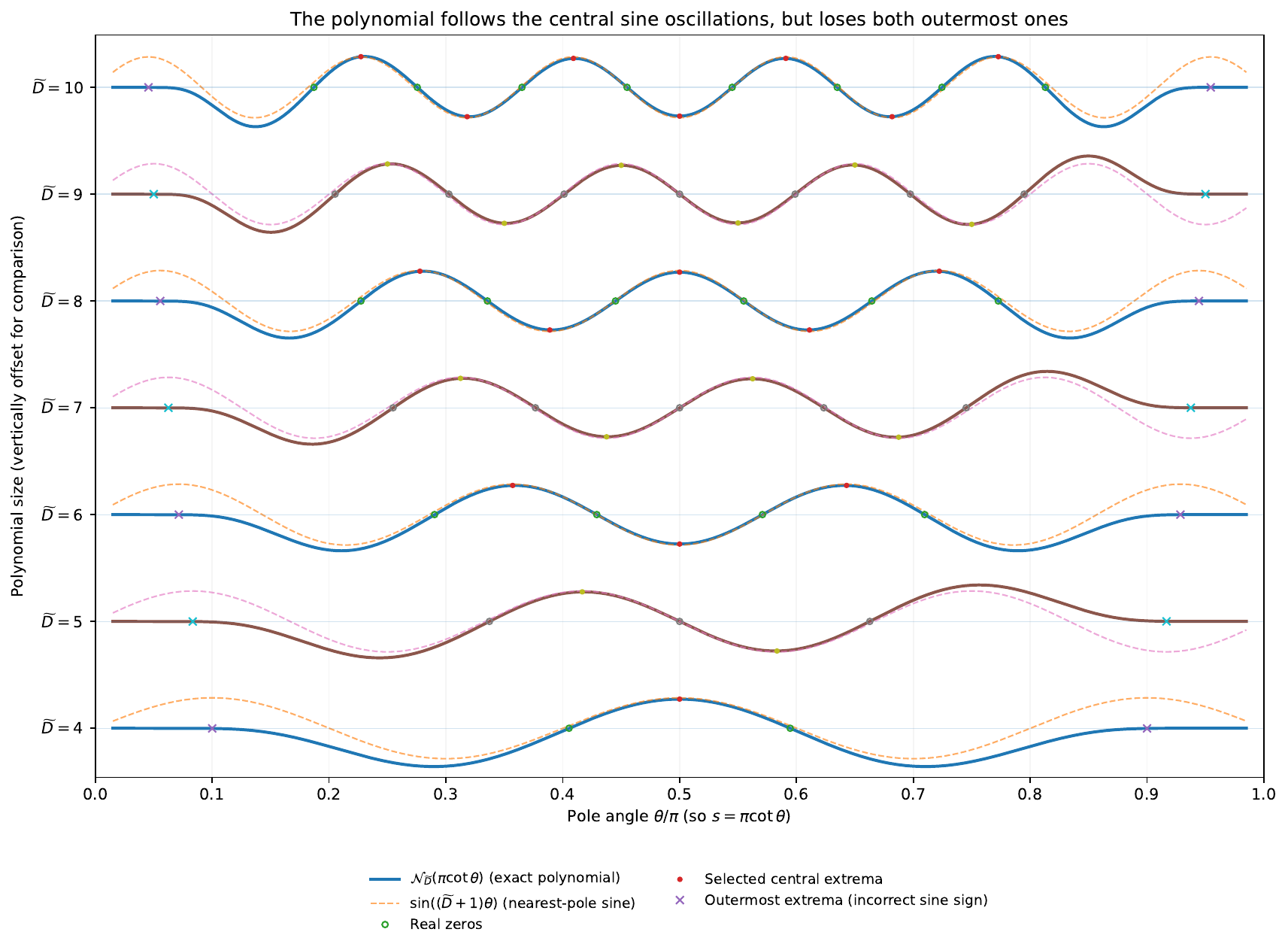}
 \caption{Small-size oscillation plots giving a complementary view of the
 normalized functions. Each curve is a numerical
 evaluation, with the angular coordinate determined by
 $s=\pi\cot\theta$. The curves exhibit the oscillatory interior behavior
 and the different endpoint behavior explained by
 \cref{eq:normalized-poles,eq:endpoint-asymptotic}.}\label{fig:small-oscillations}
\end{figure}
\FloatBarrier

\section*{Acknowledgements}

We want to acknowledge the advice and assistance of Dr. Kyle Petersen of DePaul University, one of the originators of this particular conjecture. His continued advice on this paper, the field, and throughout the process has and will continue to be invaluable. We also recognize another paper that was recently released by Dr. Per Alexandersson of Stockholm University (\cite{Alexandersson}). Our work was developed independently. 

% ===== 11_references =====

\end{document}